\documentclass[12pt,a4paper]{amsart}
\usepackage[pagewise]{lineno}
\usepackage[top=35mm, bottom=32mm, left=30mm, right=30mm]{geometry}
\usepackage{newtxtext,newtxmath}

\usepackage{mathrsfs}
\usepackage{fancyhdr} 
\usepackage{url}
\usepackage[all]{xy}
\usepackage[colorlinks=
true
, linkcolor=red, urlcolor=red, citecolor=red, linktoc=page,citecolor=blue]{hyperref}

\theoremstyle{plain}
\newtheorem{claim}{Claim}
\newtheorem{thm}{Theorem}[section]
\newtheorem{lem}[thm]{Lemma}
\newtheorem{prop}[thm]{Proposition}
\newtheorem{cor}[thm]{Corollary}

\theoremstyle{definition}
\newtheorem{defn}[thm]{Definition}
\newtheorem{rem}[thm]{Remark}

\newcommand{\Z}{{\mathbb{Z}}}
\newcommand{\N}{\mathbb{N}}
\newcommand{\F}{\mathcal F}

\newcommand{\ep}{\varepsilon}
\newcommand{\ra}{\rightarrow}
\newcommand{\lra}{\rightarrow}

\DeclareMathOperator{\diam}{diam}

\newcommand{\Sym}{\mathrm{Sym}}  
\newcommand{\Map}{\mathrm{Map}}
\newcommand{\ord}{\mathrm{ord}}
\newcommand{\mdim}{\mathrm{mdim}}

\def\v{\varepsilon}
\def \d {\delta}

\def \ov {\overline}

\begin{document}

	\title{Sofic conditional mean dimension, relative sofic mean dimension and their localizations}
	
	\author{Xianqiang Li, Zhuowei Liu *, and Xiaofang Luo}
	\date{\today}
	\address[X. Li]{School of Mathematics (Zhuhai), Sun Yat-sen University,
		Zhuhai, Guangdong, 519000, P.R. China}
	\email{lixq233@mail2.sysu.edu.cn}
	\address[Z. Liu]{School of Mathematics (Zhuhai), Sun Yat-sen University,
		Zhuhai, Guangdong, 519000, P.R. China}
	\email{liuzhw55@mail2.sysu.edu.cn}
	\address[X. Luo]{School of Mathematics (Zhuhai), Sun Yat-sen University, Zhuhai, Guangdong, 519000, P.R. China}
	\email{luoxf29@mail2.sysu.edu.cn}
	
	\thanks{* Zhuowei Liu is the corresponding author.}
	
	\begin{abstract}
		Let $\pi:(X,G)\to (Y,G) $ be a factor map between continuous actions of a sofic group $G$, we study sofic conditional mean dimension and relative sofic mean dimension introduced in \cite{LBB2} and \cite{LB}, respectively. We obtain that if $\pi$ has non-negative sofic conditional mean dimension (resp. relative sofic mean dimension), then $\pi$ has the maximal zero sofic conditional mean dimension factor (resp. maximal relative zero  sofic mean dimension factor).
		
		Additionally, the local properties of sofic conditional mean dimension and relative sofic mean dimension are studied. We introduce the sofic conditional mean dimension tuples and relative sofic mean dimension tuples, and show that $\pi$ has positive sofic conditional mean dimension (resp. relative sofic mean dimension) if and only if the set of sofic conditional mean dimension tuples (resp. relative sofic mean dimension tuples) is nonempty.
	\end{abstract}
	
	\maketitle
	%
	%
	
	\tableofcontents
	\section{Introduction}\label{sec1}

	By a \textit{topological dynamical system}, we mean a pair $(X,G)$, where $X$ is a compact metric space and $G$ acts on $X$ as a group of homeomorphsims.	Let $(X,G)$ and $(Y,G)$ be two topological dynamical systems. A continuous map $\pi: X \to Y$ is called  a \textit{factor map} between $(X,G)$ and $(Y,G)$ if it is onto and $\pi(gx) = g\pi(x)$ for all $g \in G$ and $x \in X$.
	
	Mean dimension is a  topological invariant of topological dynamical systems introduced by Gromov in \cite{MG} and was studied systematically by Lindenstrauss and Weiss within the framework of amenable group actions. It plays an important role in the study of embedding topological dynamical systems in certain natural systems (e.g. see \cite{G11,G15,G17,GLT16,GT14}). In \cite{LH}, Li extended the theory to the sofic setting by defining sofic mean dimension, thereby generalising the work of Lindenstrauss–Weiss from amenable to arbitrary sofic groups. Since then, sofic mean dimension has attracted considerable attention (e.g. see \cite{GLT16,LQ,LB}).
	
	The concept of sofic groups was first introduced by Gromov \cite{GM} and explicitly by Weiss in \cite{W00}. In \cite{Bow10b}, Bowen pioneered the study of entropy for measure-preserving actions of sofic groups on probability spaces,  which extends Ornstein-Weiss entropy classification of Bernoulli shifts far beyond the scope of countably infinite amenable groups, covering a large class of nonamenable groups, in particular including all nontorsion countable sofic groups. Extending this framework, Kerr and Li \cite{KL11b,KERRLIHD} introduced the notions of topological entropy  and further established a variational principle bridging the measure-theoretic and topological perspectives. Subsequently, Zhang \cite{Zhangguohua} studied the local properties of sofic entropy, proving a local variational principle for arbitrary finite open covers.

	For a factor map $\pi:(X,G)\to (Y,G)$ between between continuous actions of a sofic group $G$, Luo \cite{Luo} investigated conditional sofic entropy and relative sofic entropy of $Y$ with respect to $X$ for both topological and measurable dynamical systems, extending the variational principle to conditional settings. In particular, when $G$ is amenable, the author in \cite{Luo} showed that relative sofic entropy of $X$ with respect to $X$ is equal to the classical absolute amenable entropy of $Y$.	In \cite{LB}, Li and Liang defined the sofic mean dimension of $(Y,G)$ relative to $(X,G)$ and showed that it is equal to the mean dimension of $(Y,G)$ when the group $G$ is a countably infinite amenable group. Subsequently, building on Luo’s study of conditional entropy and Tsukamoto's \cite{TsukaMa} study on relative mean dimension, in \cite{LBB1} Liang introduced the notions of conditional mean dimension for amenable group actions and showed some applications in the aspect of embedding problems of dynamical systems. Recently, Liang \cite{LBB2} extended the notions of conditional mean dimension from amenable group actions to sofic group actions.

	\medskip
	
	Since mean dimension was introduced later than entropy, the theory of entropy has developed into a relatively more mature framework compared to that of mean dimension.
	In recent years, there has been growing interest in establishing mean dimension analogues of classical entropy concepts. 
	Notable examples include: conditional sofic mean dimension \cite{LBB1} as the counterpart of conditional entropy, relative sofic mean dimension \cite{LB} as the analogue of relative sofic entropy \cite{Luo}, sofic mean dimension pairs and tuples \cite{GRFYG} as the counterparts of entropy pairs and entropy tuples.
	
	Local entropy theory is an important part in entropy theory (see for example \cite{GL24,GY09}), which is introduced by  Blanchard \cite{B1}. The original motivation was to study the system whose non-trivial factors have positive entropy.
	Using the idea from the local entropy theory, García-Ramos and Gutman \cite{{GRFYG}} developed the local mean dimension theory and the notion of sofic mean dimension pairs was introduced. Let  $G$ be a sofic group, $\Sigma$ a sofic approximation sequence for $G$ and $(X,G)$ be a topological dynamical system. A pair $(x,y)\in X\times X$ is said to be a \emph{sofic mean dimension pair} if for every admissible open cover $\alpha$ with respect to $(x,y)$, the sofic mean dimension $\text{mdim}_{\Sigma}(\alpha)$ of the open cover $\alpha$ is positive. Moreover, they introduced the notions of \emph{systems with completely/uniform positive mean dimension}, which is the counterpart of  U.P.E. and C.P.E. (introduced by  Blanchard \cite{B1}) in entropy theory. After that, they give sufficient conditions of when every non-trivial factor of a sofic group action system has positive mean dimension. 
	
	Continuing this programme, it is natural to develop the local theory of sofic conditional mean dimension counterpart of entropy theory. Such as relative entropy tuples, relative U.P.E., C.P.E. extentions \cite{HYZhang}, entropy sets and so on.
	
	Motivated by these developments,
	this paper takes advantage of the techniques from local entropy theory to study mean dimension and advances the theory of sofic conditional mean dimension, relative sofic mean dimension and their localizations. Our main contributions are:
	
	\begin{enumerate}
		\item \textbf{Unification of Definition}:
		We give a definition of sofic conditional  mean dimension using finite open covers in the spirit of \cite{LBB1} and prove equivalency between the finite open cover definition  and Liang’s embedding approach \cite{LBB2} for sofic conditional mean dimension (see Proposition~\ref{xoxoxoxo}).

		\item \textbf{Maximal Zero Factor Theorem}: For a factor map $\pi:(X,G)\to (Y,G)$ between between continuous actions of a sofic group $G$, we establish, for factor maps with non-negative conditional mean dimension, the existence of a maximal zero sofic conditional mean dimension factor—extending Lindenstrauss-Weiss \cite[Proposition 6.12]{LW} to the relative sofic setting (see Theorem \ref{main1}). When $(Y,G)$ is trivial, the maximal zero sofic conditional mean dimension factor recover universal zero sofic mean dimension factor introduced by García-Ramos and Gutman \cite{{GRFYG}}.
		
		\item \textbf{Localization via Tuples}: We introduce the notion of \textbf{sofic conditional mean dimension tuples} and prove that the existence of such a tuple characterizes factor maps with positive sofic conditional mean dimension (see Theorem~\ref{main3}). Moreover, we show that the sofic conditional mean dimension tuples has the similar properties to relative entropy tuples.
		
		\item \textbf{Classification of Factor maps}: Building on local mean dimension theory, we define the notions of sofic $X|\pi$-CPCMD (Completely Positive Conditional Mean Dimension) and sofic $X|\pi$-UPCMD (Uniformly Positive Conditional Mean Dimension), studying their properties and the relationship to sofic conditional mean dimension tuples. In particular, one can characterize sofic $X|\pi$-CPCMD and sofic $X|\pi$-UPCMD via sofic conditional mean dimension tuples.
		
		\item \textbf{Relative Mean Dimension Structures}: We prove analogous theorems for relative sofic mean dimension, including maximal zero relative sofic mean dimension factors (see Theorem~\ref{main3}), relative sofic mean dimension tuples (see Theorem~\ref{main4}), and sofic $Y|X$- CP/UPMD properties.
		
		\item \textbf{Mean Dimension Sets}: Similar to the notions of entropy sets in local entropy theory, we propose the notions of \textbf{sofic conditional  mean dimension sets} and \textbf{relative sofic mean dimension sets}, which is a generalization of sofic conditional mean dimension tuples and relative sofic  mean dimension tuples.

		\item \textbf{Amenable Reconciliation}: We show the reconciliation of sofic conditional mean dimension with amenable case for infinite amenable groups, yielding a relative mean dimension product formula:
		\[
		\mathrm{mdim}(X|\pi)=\frac1n\mathrm{mdim}(X^n|\pi^n)
		\]
		generalizing Jin-Qiao's absolute version \cite{LQ}.
	\end{enumerate}
	
	This work systematically advances the localization theory and relative properties of mean dimension for sofic group actions, unifying prior approaches and enabling new applications in the classification and embedding of dynamical systems.
	
	The paper is organized as follows. In Section \ref{sec2}, we recall some definitions and some related theorems. In 
	Section \ref{sec3},  we use finite open cover to define sofic conditional mean dimension and study maximal zero sofic conditional mean dimension factor. In Section \ref{sec4}, we introduce the notions of sofic conditional mean dimension tuples. In Sections \ref{sec5}
	and \ref{sec6}, we study analogous theorems for relative sofic mean dimension. In Section \ref{sec7}, we introduce the notions of sofic conditional  mean dimension sets and relative sofic mean dimension sets. In Section \ref{sec8} we study what happens when $G$ is an infinite amenable group.

	\medskip
	
	\textit{Acknowledgments}.
	The authors would like to thank their doctoral supervisor, Siming Tu, who provided with much guidance and assistance.

	\section{Preliminaries}\label{sec2}
	In this section, we will present theoretical foundations necessary for this paper. The following notations will be adopted consistently. 
	
	Throughout this paper, $X$ and $Y$ will always represent two compact  metric spaces and G a countable group with identity $e_G$.  Define:
	$
	X^{n} = X \times \cdots \times X \quad \text{($n$-times)}, 
	$
	and let $\Delta_n(X) = \{(x_i)_{i=1}^n \in X^{n} \mid x_1 = \cdots = x_n\}$ denote the $n$-th diagonal of $X$. Let $\pi:X\to Y$ be a  continuous surjective  map. Define:
	$$
	R_\pi^n = \left\{(x_i)_{i=1}^n \in X^n \mid \pi(x_1) = \cdots = \pi(x_n)\right\},
	$$
	abbreviated as $R_\pi$ when $n = 2$. The open ball around $x\in X$ is denoted by $B_\ep(x)$. Its closure is denoted by $\overline{B_\ep(x)}$.
	
	\subsection{Topological dynamical systems} By
	a \emph{$G$-topological dynamical system} ($G$-system for short) $(X, G)$
	we mean that $X$ is a compact metric space and $\Gamma \colon G \times X \to X$, given by $(g, x) \mapsto gx$, is a continuous mapping satisfying $\Gamma(e_G, x) = x$ for every $x \in X$ and $\Gamma(g_1, \Gamma(g_2, x)) = \Gamma(g_1 g_2, x)$ for every $g_1, g_2 \in G$ and $x \in X$. A subset $Y \subset X$ is called $G$\emph{-invariant} if $gy \subset Y$ for all $y \in Y$ and $g \in G$. If $(X, G)$ and $(Y, G)$ are two $G$-systems, their product system is the system $(X \times Y, G)$, where $g(x, y) = (gx, gy)$. A nonempty closed $G$-invariant subset $Y \subset X$ naturally defines a subsystem $(Y, G)$ of $(X, G)$.
	
	\subsection{Factor maps}
	Let $(X,G)$ and $(Y,G)$ be two $G$-systems. A continuous map $\pi: X \to Y$ is called  a \textit{factor map} between $(X,G)$ and $(Y,G)$ if it is onto and $\pi(gx) = g\pi(x)$ for all $g \in G$ and $x \in X$. In this case, we say $(X,G)$ is an extension of $(Y,G)$ or $(Y,G)$ is a factor of $(X,G)$.  The systems are said to be \emph{conjugate} if $\pi$ is a homeomorphism.
	From a dynamical point of view, conjugate $G$-systems can be considered equivalent.
	Let $\pi\colon (X,G)\lra (Y,G)$ be a factor map between two $G$-systems. Then 
	$$
	R_\pi=\{(x_1,x_2) \colon \pi(x_1)=\pi(x_2)\}
	$$
	is a closed $G$-invariant equivalence relation on $X$ and $Y=X/R_\pi$. 
	Moreover, there is a one-to-one correspondence between
	the collection of factors of $(X,G)$ and the collection of closed $G$-invariant equivalence relations on $X$ (see for example \cite[Appendix E.11-3]{Vries1993}).
	
	Let $\pi\colon (X,G)\lra (Y,G)$  and $\psi\colon (Z,G)\ra (Y,G)$ be two factor maps between $G$-systems. We say that $\psi$ is a \emph{factor of $\pi$} if there exists a factor map $\phi\colon (X,G)\ra (Z,G)$ such that $\pi=\phi\circ\psi$.
	\subsection{Sofic groups}
	\begin{defn}
		For $d \in \N$, we write $[d]$ for the set $\{1,\cdots,d\}$ and  \(\Sym(d)\) for the permutation group of $[d]$. A countable group $G$ is called \emph{sofic group} if there is a \emph{sofic approximation sequence} $\Sigma=\{\sigma_{i}:G \rightarrow \text{Sym}(d_{i})\}_{i=1}^{\infty}$ for $G$, namely the following three conditions are satisfied:
		\begin{itemize}
			\item [(1)]for any $s,t \in G$, one has 
			$$\lim_{i \rightarrow \infty}\frac{|\{v \in [d_{i}]:\sigma_{i}(s)\sigma_{i}(t)(v)=\sigma_{i}(st)(v)\}|}{d_{i}}=1;$$
			\item [(2)]for any distinct $s,t \in G$, one has $$\lim_{i \rightarrow \infty}\frac{|\{v \in [d_{i}]:\sigma_{i}(s)(v)=\sigma_{i}(t)(v)\}|}{d_{i}}=0;$$
			\item [(3)] $\lim_{i \rightarrow \infty}d_{i}=+\infty.$
		\end{itemize}
	\end{defn}
	For a map $\sigma$ from $G$ to  \(\Sym(d)\) for some $d \in \N$, we write $\sigma(s)(v)$ as $\sigma_{s}(v)$ or $sv$.

	Throughout this paper, we fix a countable sofic group $G$ and a sofic approximation sequence $\Sigma=\{\sigma_{i}:G \rightarrow \text{Sym}(d_{i})\}_{i=1}^{\infty}$ for $G$.  Denote by $\F(G)$  the collection of nonempty finite subsets of $G$.
	
	Let \(\rho\) be a continuous pseudometric on a compact  metric space \(X\). For a given \(d \in \mathbb{N}\), we define on the set of all maps from \([d]\) to \(X\) the pseudometrics
	\[
	\rho_2(\xi, \tau) = \left( \sum_{a \in [d]} \big( \rho(\xi(a), \tau(a)) \big)^2 \right)^{1/2},
	\]
	\[
	\rho_\infty(\xi, \tau) = \max_{a \in [d]} \rho(\xi(a), \tau(a)).
	\]
	
	\begin{defn}
		Let \(F \in \F(G)\) and \(\delta > 0\). Let \(\sigma\) be a map from $G$ to \(\Sym(d)\) for some \(d \in \N\). We define \(\Map(\rho, F, \delta, \sigma)\) to be the set of all maps \(\xi \colon [d] \to X\) such that 
		\[
		\rho_2(\xi \circ \sigma_s, \alpha_s \circ \xi) \leq \delta \quad \text{for all } s \in F.
		\]
		where $\alpha_{s}$ denotes the transformation $x \mapsto sx$ of $X$. We consider \(\Map(\rho, F, \delta, \sigma)\) to be a topological space with the topology inherited from \(X^d\).
	\end{defn}
	
	For a finite open cover $\alpha$ of \(X\) and $d\in\N$, we denote by $\alpha^d$ the finite open cover of \(X^{d}\) consisting of 
	\[
	U_1 \times U_2 \times \cdots \times U_d \quad \text{for } U_1, \dots, U_d \in \alpha.
	\]
	Note that \(\Map(\rho, F, \delta, \sigma)\) is a closed subset of \(X^{d}\). Consider the restriction 
	\[
	\alpha^d|_{\Map(\rho,F,\delta,\sigma)} = \alpha^d \cap \Map(\rho, F, \delta, \sigma)
	\] 
	of \(\alpha^d\) to \(\Map(\rho, F, \delta, \sigma)\). 
	
	We need the following lemma.
	\begin{lem}\label{li-2}\cite[Lemma 2.10]{LH}
		Let $\pi: (X,G)\to (Y,G)$ be a factor map between $G$-systems. Let $\rho_X$ and $\rho_Y$ be compatible metrics on $X$ and $Y$ respectively. Let $F$ be a nonempty finite subset of $G$ and $\delta > 0$. Then there exists $\delta' > 0$ such that for every map $\sigma$ from $G$ to $\operatorname{Sym}(d)$ for some $d \in \mathbb{N}$ and every $\xi \in \operatorname{Map}(\rho_X, F, \delta', \sigma)$, one has $\pi \circ \xi \in \operatorname{Map}(\rho_Y, F, \delta, \sigma)$.
	\end{lem}

	\section{Sofic conditional mean dimension of $(X,G)$ relative to $(Y,G)$}\label{sec3}
	
	In this section, we study the sofic conditional mean dimension of $(X,G)$ relative to $(Y,G)$. We begin by recalling the notion of (covering) dimension of compact metric spaces and establishing its fundamental properties, which will play a crucial role in our subsequent definitions and proofs.
	
	Let $X$ be a compact  metric space. For two finite open covers $\alpha$ and $\beta$ of $X$, the \emph{joining} $\alpha \vee \beta$ is defined as $\alpha \vee \beta=\{U \cap V: U \in \alpha, V \in \beta \}$. We say that $\beta$ \emph{refines} $\alpha$, defined $\beta \succ \alpha$, if every member of $\beta$ is contained in a member of $\alpha$. Denote
	\[
	\mathrm{ord}(\alpha) = \left( \max_{x \in X} \sum_{U \in \alpha} 1_U(x) \right) - 1\text{,~~}
	D(\alpha) = \min_{\beta \succ \alpha} \mathrm{ord}(\beta),
	\]
	where the minimum is taken over all finite open covers \(\beta\) of \(X\) such that \(\beta\) refines \(\alpha\). It is clear that
	$D(\alpha) \leq \mathrm{ord}(\alpha) < \infty$. 
	Then the (covering) dimension of $X$ is defined by $\sup_{\alpha}D(\alpha)$, where $\alpha$ ranges over all finite open covers of $X.$
	
	\begin{defn}\cite{LBB1}
		Let $\pi:X\to Y$ be a  continuous surjective  map and \(\alpha\)  a finite open cover of \(X\). Denote 
		$$	D(\alpha|\pi) = \min_{\{\pi^{-1}(y)\}_{y\in Y}\vee\beta \succ \alpha} \mathrm{ord}(\beta),$$
		where the minimum is taken over all finite open covers \(\beta\) of \(X\) such that $\{\pi^{-1}(y)\}_{y\in Y}\vee\beta$ refines \(\alpha\). It's obvious that $D(\alpha|\pi) \leq D(\alpha)$, if $Y$ is a singleton, then $	D(\alpha|\pi)=	D(\alpha)$.
	\end{defn} 
	
	\begin{prop}\label{jbwcjbc}
		Let $\pi:X\to Y$ be a  continuous surjective map. Suppose $\alpha$ and $\beta$ are two finite open covers of $X$ such that $\beta \succ \alpha$. Then one has $D(\beta|\pi) \geq D(\alpha|\pi)$.
	\end{prop}
	\begin{proof}
		Let $\gamma$ be a finite open cover of $X$ such that $\{\pi^{-1}(y)\}_{y\in Y}\vee \gamma  \succ \beta$ and $\mathrm{ord}(\gamma)=D(\beta|\pi)$, then we have $\{\pi^{-1}(y)\}_{y\in Y}\vee \gamma  \succ \beta \succ  \alpha$, which implies that $D(\beta|\pi)=\mathrm{ord}(\gamma) \geq D(\alpha|\pi)$.
	\end{proof}
	We need the following lemma, which is a conditional version of \cite[Proposition 2.4]{LW}.
	\begin{lem}\cite [Lemma 2.2]{LBB1}\label{polyhe}
		Let $\pi:X\to Y$ be a  continuous surjective  map and $\alpha$  a finite open cover of $X$. Then $D(\alpha|\pi )\leq k$ if and only if there exists a continuous map $f: X\to P$ for some polyhedron $P$ with $\dim (P)= k$, such that $\{f^{-1}(p)\cap \pi^{-1}(y)\}_{(p,y) \in P \times Y}$ refines $\alpha$.
	\end{lem}
	\begin{prop}\label{subbb}
		Let $\pi:X\to Y$ be a  continuous surjective map. Suppose $\alpha$ and $\beta$ are two finite open covers of $X$, then
		$$D(\alpha \vee \beta|\pi) \leq D(\alpha|\pi)+ D(\beta|\pi).$$
	\end{prop}
	\begin{proof}
		By   Lemma~\ref{polyhe},  there exists a continuous map $f: X \rightarrow P$ for some polyhedron $P$ with $\dim(P)=D(\alpha|\pi)$ such that
		$$\{f^{-1}(p)\cap \pi^{-1}(y)\}_{(p,y) \in P \times Y} \succ \alpha.$$
		Meanwhile, there exists a continuous map $g: X \rightarrow Q$ for some polyhedron $Q$ with $\dim (Q)=D(\beta|\pi)$ such that
		$$\{g^{-1}(q)\cap \pi^{-1}(y)\}_{(q,y) \in Q \times Y} \succ \beta.$$
		Define the map 
		$h:X \rightarrow P \times Q$ with $h: x \rightarrow (f(x),g(x))$. It's obvious that $h$ is continuous. Moreover, we have
		$$\{f^{-1}(p)\cap g^{-1}(q)\cap \pi^{-1}(y)\}_{(p,q,y) \in P \times Q \times Y} \succ \alpha \vee \beta.$$
		Using Lemma~\ref{polyhe}  again, we get
		\begin{align*}
			D(\alpha \vee \beta|\pi) \leq \dim(P \times Q)
			\leq \dim(P)+\dim(Q)=D(\alpha|\pi)+ D(\beta|\pi).
		\end{align*}
		The proof is complete.
	\end{proof}

	\begin{lem}\label{lemlk}
		Let $\pi_1:X_1\to Y$ and  $\pi_2:X_2\to Y$ be two  continuous surjective  maps. Suppose $f:X_1\to X_2$ is a continuous surjective map with $\pi_1=\pi_2\circ f$ and $\alpha$ is a finite open cover of $X_2$. Then $D(f^{-1}(\alpha)|\pi_1)\leq D(\alpha|\pi_2)$.
	\end{lem}	
	\begin{proof}
		Let $\beta$ be a finite open cover of $X_2$ such that $\{\pi_2^{-1}(y)\}_{y\in Y}\vee\beta\succ\alpha$ and $\ord(\beta)=D(\alpha|\pi_2)$. Then $f^{-1}(\beta)$ is a finite open cover of $X_1$, and $f^{-1}(\{\pi_2^{-1}(y)\}_{y\in Y}\vee\beta)=f^{-1}(\{\pi_2^{-1}(y)\}_{y\in Y})\vee f^{-1}(\beta)\succ f^{-1}(\alpha)$.
		For any $y\in Y$, $f^{-1}\circ\pi_2^{-1}(y)=\pi_1^{-1}(y)$, then $f^{-1}(\{\pi_2^{-1}(y)\}_{y\in Y})=\{\pi_1^{-1}(y)\}_{y\in Y}$.
		This implies that $\{\pi_1^{-1}(y)\}_{y\in Y}\vee f^{-1}(\beta)\succ f^{-1}(\alpha)$.
		Hence
		\begin{align*}
			D(f^{-1}(\alpha)|\pi_1)\leq\ord(f^{-1}(\beta))=\ord(\beta)=D(\alpha|\pi_2).
		\end{align*}
		This finishes the proof.
	\end{proof}
	
	\begin{prop}\label{dcjkdbcsd}
		For any $d \in \N$, let $\pi_{i}:X_{i} \rightarrow Y_{i}$ be a continuous surjective map and $\alpha_{i}$ a finite open cover of $X_{i}$, where $i \in \{1,2,\dots,d\}$. Then
		$$D(\prod_{i=1}^{d}\alpha_{i}|\prod^{d}_{i=1}\pi_{i}) \leq \sum_{i=1}^{d}D(\alpha_{i}|\pi_{i}).
		$$
		Moreover, if $X_{i}=X,Y_{i}=Y, \pi_{i}=\pi$ for all $i\in \{1,2,\dots,d\}$, then
		$$D(\alpha^d|\pi^{d})\leq d D(\alpha|\pi).$$
	\end{prop}	
	\begin{proof}
		Let $p_{i}:\prod_{i=1}^{d}X_{i} \rightarrow X_{i}$ be the projection on the $i$-th space, then 	$\prod_{i=1}^{d}\alpha_{i}=\bigvee_{i=1}^{d}p^{-1}_{i}(\alpha).$
		
		\begin{claim}\label{claim1}
			For any $i \in \{1,2,\dots,d \}$, we have
			$$D(p^{-1}_{i}(\alpha)|\prod^{d}_{i=1}\pi_{i})\leq D(p^{-1}_{i}(\alpha)|\pi_{i}\circ p_{i}).$$
		\end{claim}
		By Claim \ref{claim1}, we get
		Then
		\begin{align*}
			D(\prod_{i=1}^{d}\alpha_{i}|\prod^{d}_{i=1}\pi_{i}) &= D(\bigvee_{i=1}^{d}p^{-1}_{i}(\alpha)|\prod^{d}_{i=1}\pi_{i})\\
			&\leq \sum_{i=1}^{d}D(p^{-1}_{i}(\alpha)|\prod^{d}_{i=1}\pi_{i}) ~~(~\text{by~ Proposition ~\ref{subbb}} )\\
			&\leq \sum_{i=1}^{d}D(p^{-1}_{i}(\alpha)|\pi_{i}\circ p_{i})~~(~\text{by~ Claim ~\ref{claim1}} )\\
			&\leq \sum_{i=1}^{d}D(p^{-1}_{i}(\alpha)|\pi_{i})~~(~\text{by~ Lemma ~\ref{lemlk}} ).
		\end{align*}
		\emph{Proof of Claim \ref{claim1}}. Suppose $\beta$ is a finite open cover of $p^{-1}_{i}(\alpha)$, such that
		$\{p^{-1}_{i}\circ \pi_{i}^{-1}(y_i)\}_{y_i \in Y_{i}} \vee \beta \succ p^{-1}_{i}(\alpha),$
		and $D(p^{-1}_{i}(\alpha)|\pi_{i}\circ p_{i})=\ord(\beta)$.
		
		Note that 	$(\prod^{d}_{i=1}\pi_{i})^{-1}(y)
		\subset p^{-1}_{i}\circ \pi_{i}^{-1}(y_i),$
		for any $y=(y_i)_{i=1}^d \in \prod^{d}_{i=1}Y_{i}$,
		then
		$$\big\{(\prod^{d}_{i=1}\pi_{i})^{-1}(y)\big\}_{y=(y_i)_{i=1}^d \in \prod^{d}_{i=1}Y_{i}} \vee \beta \succ p^{-1}_{i}(\alpha),$$
		which implies that
		\begin{align*}
			D(p^{-1}_{i}(\alpha)|\prod^{d}_{i=1}\pi_{i}) \leq \ord(\beta)= D(p^{-1}_{i}(\alpha)|\pi_{i}\circ p_{i}).
		\end{align*}
		The proof of claim is complete.
	\end{proof}	
	\begin{defn}
		Let $\pi:X\to Y$ be a  continuous surjective  map, \(K\subset X\)  a  nonempty subset of \(X\) and \(\alpha\) a finite open cover of \(X\). Denote
		\[
		\mathrm{ord}(\alpha|_K) = \left( \max_{x \in K} \sum_{U \in \alpha} 1_U(x) \right) - 1
		\text{,~~}
		D(\alpha|_K\big|\pi) = \min_{\{\pi^{-1}(y)\}_{y\in Y}\vee\beta\succ \alpha} \mathrm{ord}(\beta|_K),
		\]
		where the minimum is taken over all finite open covers \(\beta\) of \(X\) such that $\{\pi^{-1}(y)\}_{y\in Y}\vee\beta$ refines \(\alpha\). In addition, we define
		$
		\mathrm{ord}(\alpha|_\emptyset) = -\infty
		$.
	\end{defn}

	\begin{lem}\label{subcover}
		Let $\pi:X\to Y$ be a  continuous surjective  map, \(K\subset X\)  a  nonempty closed subset and \(\alpha\) a finite open cover of \(X\). Then 
		
		$$D(\alpha|_K\big|\pi) = \min_{\{(\pi|_K)^{-1}(y)\}_{y\in \pi(K)}\vee\beta\succ \alpha|_K} \mathrm{ord}(\beta|_K),$$ 
		where the minimum is taken over all finite open covers \(\beta\) of $K$ such that $\{(\pi|_K)^{-1}(y)\}_{y\in \pi(K)}\vee\beta$ refines $\alpha|_K$.
	\end{lem}
	\begin{proof}
		Let $\beta$ be a finite open cover of $X$ such that $\{\pi^{-1}(y)\}_{y \in Y} \vee \beta \succ \alpha$ and  $D(\alpha|_K\big|\pi) =\mathrm{ord}(\beta|_K)$. Then
		$\{\pi^{-1}(y)\}_{y \in Y} \vee \beta|_K \succ \alpha|_K$, where $\beta|_K$ is the restriction of $\beta$ to $K$. Since $\{(\pi|_K)^{-1}(y)\}_{y \in \pi(K)}  \subset \{\pi^{-1}(y)\}_{y \in Y} $, then $\{(\pi|_K)^{-1}(y)\}_{y \in \pi(K)} \vee \beta|_K \succ \alpha|_K$. This implies that
		$$\min_{\{(\pi|_K)^{-1}(y)\}_{y\in \pi(K)}\vee\beta\succ \alpha|_K} \mathrm{ord}(\beta|_K) \leq \mathrm{ord}(\beta|_K)=D(\alpha|_K\big|\pi),$$ 
		where the minimum is taken over all finite open covers \(\beta\) of $K$ such that $\{(\pi|_K)^{-1}(y)\}_{y\in \pi(K)}\vee\beta$ refines $\alpha|_K$.
		
		On the other hand, suppose $\gamma$ is a finite open cover of $K$ such that $\{(\pi|_K)^{-1}(y)\}_{y\in \pi(K)}\vee\gamma \succ \alpha|_K$ such that 
		$$\min_{\{(\pi|_K)^{-1}(y)\}_{y\in \pi(K)}\vee\beta\succ \alpha|_K} \mathrm{ord}(\beta|_K) =\mathrm{ord}(\gamma|_K).$$ 
		There exists a finite open cover $\zeta$ of $X$ such that $\zeta|_K=\gamma$ and $\{\pi^{-1}(y)\}_{y \in Y} \vee \zeta \succ \alpha$, then
		$$D(\alpha|_K\big|\pi)  \leq \mathrm{ord}(\gamma|_K)=\mathrm{ord}(\zeta|_K)=\min_{\{(\pi|_K)^{-1}(y)\}_{y\in \pi(K)}\vee\beta\succ \alpha|_K} \mathrm{ord}(\beta|_K).
		$$
		This implies the desired result.
	\end{proof}
	
	Now we can give the definition of sofic conditional mean dimension  using finite open covers.
	
	\begin{defn}\label{defnnnnn1}
		Let $\pi: (X,G)\to (Y,G)$ be a factor map between $G$-systems, \(F \in \F(G)\) and \(\delta > 0\). For a finite open cover \(\alpha\) of \(X\) and a sofic approximation sequence for \(G\), \(\Sigma = \{\sigma_i: G \to \Sym(d_i)\}_{i=1}^{\infty}\), we define
		
		\begin{align*}
			D(\alpha|\pi, \rho, F, \delta, \sigma) :&= D\left(\alpha^{d}|_{\Map(\rho,F,\delta,\sigma)} \big| \pi^d|_{\Map(\rho,F,\delta,\sigma)}\right)\\&=\min_{\{(\pi^d|_{\Map(\rho,F,\delta,\sigma)})^{-1}(y)\}_{y\in \pi^d({\Map(\rho,F,\delta,\sigma)})}\bigvee\beta \succ \alpha^d|_{\Map(\rho,F,\delta,\sigma)}} \mathrm{ord}(\beta|_{\Map(\rho,F,\delta,\sigma)}),
		\end{align*}
		where $\pi^d : X^d\to Y^d$ sends $\xi$ to $\pi\circ\xi$ and the minimum is taken over all finite open covers \(\beta\) of $\Map(\rho,F,\delta,\sigma)$ such that $\{(\pi^d|_{\Map(\rho,F,\delta,\sigma)})^{-1}(y)\}_{y\in \pi^d({\Map(\rho,F,\delta,\sigma)})}\bigvee\beta$ refines \(\alpha^d|_{\Map(\rho,F,\delta,\sigma)}\). By Lemma~\ref{subcover}, one has
		\begin{align*}
			D(\alpha|\pi, \rho, F, \delta, \sigma) &= D\left(\alpha^{d}|_{\Map(\rho,F,\delta,\sigma)} \big| \pi^d\right)\\&=\min_{\{(\pi^d)^{-1}(y)\}_{y\in Y^d}\bigvee\beta \succ \alpha^d} \mathrm{ord}(\beta|_{\Map(\rho,F,\delta,\sigma)})
		\end{align*}
		where the minimum is taken over all finite open covers \(\beta\) of $X^d$ such that $\{(\pi^d)^{-1}(y)\}_{y\in Y^d}\bigvee\beta$ refines \(\alpha^d\).
		
		Define
		\[\mathrm{mdim}_{\Sigma}(\alpha|\pi, \rho, F, \delta) := \lim_{i \to \infty} \frac{D(\alpha|\pi, \rho, F, \delta, \sigma_i)}{d_i},
		\]
		\[
		\mathrm{mdim}_{\Sigma}(\alpha|\pi, \rho,F) := \inf_{\delta > 0} \mathrm{mdim}_{\Sigma}(\alpha|\pi, \rho, F, \delta)
		\]
		and
		\[
		\mathrm{mdim}_{\Sigma}(\alpha|\pi,\rho) := \inf_{F' \in \F(G)} \mathrm{mdim}_{\Sigma}(\alpha|\pi,\rho, F').
		\]
		Note that if \(\Map(\rho, F, \delta, s_i)\) is empty for all sufficiently large \(i\), then
		\[
		\mathrm{mdim}_{\Sigma}(\alpha|\pi, \rho, F, \delta) = -\infty.
		\]
		Define the \emph{sofic conditional mean dimension (with respect to \(\Sigma\)) of $(X,G)$ relevant to $\pi$} as
		\[
		\mathrm{mdim}_{\Sigma}(X|\pi,\rho) := \sup_{\alpha} \mathrm{mdim}_{\Sigma}(\alpha|\pi,\rho),
		\]
		where the supremum is taken over all finite open covers \(\alpha\) of \(X\). 
		
		Similar to the proof of \cite[Lemma 2.9]{LH}, we can prove that the quantities $\mathrm{mdim}_{\Sigma}(\alpha|\pi, \rho,F)$, $\mathrm{mdim}_{\Sigma}(\alpha|\pi, \rho)$ and $\mathrm{mdim}_{\Sigma}(X|\pi,\rho)$ do not depend on the choice of $\rho$, and we shall write them as $\mathrm{mdim}_{\Sigma}(\alpha|\pi, F)$, $\mathrm{mdim}_{\Sigma}(\alpha|\pi)$ and $\mathrm{mdim}_{\Sigma}(X|\pi)$ respectively. When $Y$ is a singleton, $\mathrm{mdim}_{\Sigma}(X|\pi)$ recovers as the \emph{sofic mean dimension} $\mathrm{mdim}_{\Sigma}(X)$ of $(X,G)$(see \cite[Definition 2.4]{LH}). To facilitate distinction, we will write  $D(\alpha|\pi, \rho, F, \delta, \sigma)$ and $\mathrm{mdim}_{\Sigma}(\alpha|\pi)$ replace by $D(\alpha, \rho, F, \delta, \sigma)$ and $\mathrm{mdim}_{\Sigma}(\alpha)$.
	\end{defn}
	
	We now state the elementary properties of sofic conditional mean dimension.
	
	\begin{prop}\label{lllplpl}
		Let $\pi: (X,G)\to (Y,G)$ be a factor map between $G$-systems, $\Sigma$ a sofic approximation sequence for $G$ and $\rho$ a compatible metric on $X$. Suppose $\alpha$  is a finite open cover of $X$, then $$\mathrm{mdim}_{\Sigma}(\alpha|\pi,\rho,F,\delta)\leq D(\alpha|\pi)$$
		and $$\mathrm{mdim}_{\Sigma}(\alpha|\pi)\leq D(\alpha|\pi).$$
	\end{prop}
	\begin{proof}
		Let $F\in \mathcal{F}(G)$, $\delta>0$, and $\sigma:G\to \text{Sym}(d)$ for some $d \in \N$. By Proposition \ref{dcjkdbcsd},
		we have  
		\begin{align*}
			D(\alpha|\pi,\rho,F,\delta,\sigma)&=D(\alpha^d|_{\text{Map}(\rho,F,\delta,\sigma)}|\pi^d)\\
			&=\min_{\{(\pi^d)^{-1}(y)\}_{y\in Y^d}\bigvee \beta\succ \alpha^d}\text{ord}(\beta|_{\text{Map}(\rho,F,\delta,\sigma)})\\
			&\leq \min_{\{(\pi^d)^{-1}(y)\}_{y\in Y^d}\bigvee
				\beta\succ \alpha^d}\text{ord}(\beta)\\
			&= D(\alpha^d|\pi^d)\\
			&\leq dD(\alpha|\pi).\end{align*}
		This implies that $$\text{mdim}_{\Sigma}(\alpha|\pi,\rho,F,\delta)=\varlimsup_{i\to \infty}\frac{D(\alpha|\pi,\rho,F,\delta,\sigma_i)}{d_i}\leq D(\alpha|\pi).$$
		Moreover, we get $\text{mdim}_{\Sigma}(\alpha|\pi)\leq D(\alpha|\pi).$
	\end{proof}
	
	\begin{prop}\label{impor}
		Let $\pi: (X,G)\to (Y,G)$ be a factor map between $G$-systems, $\Sigma$ a sofic approximation sequence for $G$ and $\alpha, \beta$ be two finite open covers of $X$. The following holds:
		\begin{itemize}
			\item [(1)]If $\alpha\succ \beta$, then $\mathrm{mdim}_{\Sigma}(\alpha|\pi)\geq \mathrm{mdim}_{\Sigma}(\beta|\pi).$
			\item[(2)] $\mathrm{mdim}_{\Sigma}(\alpha\vee \beta|\pi)\leq \mathrm{mdim}_{\Sigma}(\alpha|\pi)+\mathrm{mdim}_{\Sigma}(\beta|\pi).$
		\end{itemize}
	\end{prop}
	\begin{proof}
		Let $F\in \mathcal{F}(G)$, $\delta>0$, and $\sigma:G\to \text{Sym}(d)$ for some $d \in \N$.\\
		(1) Since $\alpha\succ \beta$, then $\alpha^d\succ \beta^d$, and in particular, $\alpha|_{\text{Map}(\rho,F,\delta,\sigma)}\succ \beta|_{\text{Map}(\rho,F,\delta,\sigma)}$. 
		By Proposition \ref{jbwcjbc}, we get
		\begin{align*}D(\alpha^d|_{\text{Map}(\rho,F,\delta,\sigma)}|\pi^d)&=
			D(\alpha|\pi,\rho,F,\delta,\sigma)\geq D(\beta^d|_{\text{Map}(\rho,F,\delta,\sigma)}|\pi^d)=D(\beta|\pi,\rho,F,\delta,\sigma).\end{align*}
		Therefore, $\text{mdim}_{\Sigma}(\alpha|\pi)\geq \text{mdim}_{\Sigma}(\beta|\pi).$\\
		(2)Since $(\alpha\bigvee \beta)^d=\alpha^d\bigvee \beta^d$, by Proposition \ref{subbb}, we get
		\begin{align*}
			D((\alpha\vee \beta)^d|_{\text{Map}(\rho,F,\delta,\sigma)}|\pi^d)&=D((\alpha\vee \beta)|\pi,\rho,F,\delta,\sigma)\\
			&\leq D(\alpha^d|_{\text{Map}(\rho,F,\delta,\sigma)}|\pi^d)+D(\beta^d|_{\text{Map}(\rho,F,\delta,\sigma)}|\pi^d)\\
			&=D(\alpha|\pi,\rho,F,\delta,\sigma)+D(\beta|\pi,\rho,F,\delta,\sigma).
		\end{align*}
		Thus $\text{mdim}_{\Sigma}(\alpha\vee \beta|\pi)\leq \text{mdim}_{\Sigma}(\alpha|\pi)+\text{mdim}_{\Sigma}(\beta|\pi).$
	\end{proof}
	
	In \cite{LBB2}, Liang gave another  definition  of sofic conditional mean dimension. Next we will prove the two definitions are equivalent (see Proposition~\ref{xoxoxoxo}).

	\begin{defn}\label{scdt2}\cite{LBB2}
		Let $\pi: (X,G)\to (Y,G)$ be a factor map between $G$-systems and $\rho_X$ a compatible metrics on $X$. For fixed $\varepsilon > 0$, $F \in \mathcal{F}(G)$, consider a map $\sigma: \Gamma \to \mathrm{Sym}(d)$ and a subset $K \subset X^{d}$. We say a continuous map $f: K \to P$ for a compact  metric space $P$ is an $(\varepsilon|Y, \rho_{X,\infty})$-embedding if $\rho_{X,\infty}(\xi, \xi') < \varepsilon$ holds for every $\xi, \xi' \in K$ satisfying $f(\xi) = f(\xi')$ and $\pi \circ \xi = \pi \circ \xi'$.
		
		Denote by $\mathrm{Wdim}_{\varepsilon}(K, \rho_{X,\infty}|Y)$ the minimal (covering) dimension $\dim (Q)$ of a compact  metric space $Q$ that admits an $(\varepsilon|Y, \rho_{X,\infty})$-embedding from $K$ to $Q$. If $Y$ is a singleton, we simplify $\mathrm{Wdim}_{\varepsilon}(K, \rho_{X,\infty}|Y)$ as $\mathrm{Wdim}_{\varepsilon}(K, \rho_{X,\infty})$.
		
		Let $\Sigma = \{ \sigma_i: G \to \mathrm{Sym}(d_i) \}_{i=1}^{\infty}$ be a sofic approximation of $G$. Define the \emph{conditional sofic mean dimension} of $(X,G)$ relative to $(Y,G)$ (with respect to \(\Sigma\)) as
		\[
		\mathrm{mdim}_{\Sigma}(\rho_X | Y) := \sup_{\varepsilon > 0} \inf_{F \in \mathcal{F}(G)} \inf_{\delta > 0} \limsup_{i \to \infty} \frac{1}{d_i} \mathrm{Wdim}_{\varepsilon}(\mathrm{Map}(\rho_X, F, \delta, \sigma_i), \rho_{X,\infty}|Y).
		\]
		When $\varepsilon$, $F$ and $\delta$ are fixed, we write $\mathrm{mdim}_{\Sigma,\varepsilon}(\rho_X | Y, F, \delta)$ for the above limit supremum.  If $\mathrm{Map}(\rho, F, \delta, \sigma_i)$ is empty for all sufficiently large $i$, the limit supremum is set to be $-\infty$.
		
		When $\varepsilon$ is fixed, we write $$\mathrm{mdim}_{\Sigma,\varepsilon}(\rho_X | Y)=\inf_{F \in \mathcal{F}(G)} \inf_{\delta > 0}\mathrm{mdim}_{\Sigma,\varepsilon}(\rho_X | Y, F, \delta).$$ 
		
	\end{defn}

	\begin{prop}\label{xoxoxoxo}
		Let $\pi: (X,G)\to (Y,G)$ be a factor map between $G$-systems, $\Sigma$ a sofic approximation sequence for $G$ and $\rho_X$ a compatible metric on $X$.  Then
		$$
		\mathrm{mdim}_{\Sigma}(X|\pi)=\mathrm{mdim}_{\Sigma}(\rho_X | Y).
		$$
	\end{prop}
	\begin{proof}
		We first prove $	\mathrm{mdim}_{\Sigma}(X|\pi)\leq\mathrm{mdim}_{\Sigma}(\rho_X | Y)$. Fix a finite open cover $\alpha$ of $X$. Take a Lebesgue number $\lambda$ of $\alpha$. It suffices to show $$D(\alpha|\pi, \rho, F, \delta, \sigma)\leq \mathrm{Wdim}_{\lambda}(\mathrm{Map}(\rho_X, F, \delta, \sigma), \rho_{X,\infty}|Y),$$ 
		for any $F \in \mathcal{F}(G), \delta>0$ and any map $\sigma$ from $G$ to $\Sym(d)$ for some $d \in \N$.
		This implies that $\mathrm{mdim}_{\Sigma}(\alpha|\pi)\leq \mathrm{mdim}_{\Sigma,\lambda}(\rho_X | Y) \leq \mathrm{mdim}_{\Sigma}(\rho_X | Y)$. Hence $\mathrm{mdim}_{\Sigma}(X|\pi)=\sup_{\alpha}\mathrm{mdim}_{\Sigma}(\alpha|\pi)\leq\mathrm{mdim}_{\Sigma}(\rho_X | Y)$.
		
		Let $f:\mathrm{Map}(\rho_X, F, \delta, \sigma)\ra P$ be a continuous map for a compact  metric space $P$ with $\dim(P)=\mathrm{Wdim}_{\lambda}(\mathrm{Map}(\rho_X, F, \delta, \sigma_i), \rho_{X,\infty}|Y)$, such that $\rho_{X,\infty}(\xi, \xi') < \lambda$ holds for every $\xi, \xi' \in \mathrm{Map}(\rho_X, F, \delta, \sigma)$ satisfying $f(\xi) = f(\xi')$ and $\pi \circ \xi = \pi \circ \xi'$. 
		
		Then for any $p=f(\xi) = f(\xi')\in P$ and $\tau=\pi \circ \xi = \pi \circ \xi'\in \pi^d (\mathrm{Map}(\rho_X, F, \delta, \sigma))$, $\diam(f^{-1}(p)\cap(\pi^{d}|_{\Map(\rho,F,\delta,\sigma)})^{-1}(\tau),\rho_{X,\infty})<\lambda$. By the choice of $\lambda$, one has
		$$\{f^{-1}(p)\cap(\pi^d|_{\Map(\rho,F,\delta,\sigma)})^{-1}(y)\}_{(p,\tau)\in P\times \pi^{d} (\mathrm{Map}(\rho_X, F, \delta, \sigma))}\succ\alpha^d|_{\mathrm{Map}(\rho_X, F, \delta, \sigma)}.$$
		
		By Lemma~\ref{polyhe}, one has
		\begin{align*}
			D(\alpha|\pi, \rho, F, \delta, \sigma)&\leq\dim(P)\\
			&= \mathrm{Wdim}_{\lambda}(\mathrm{Map}(\rho_X, F, \delta, \sigma), \rho_{X,\infty}|Y).
		\end{align*}
		
		Next, we show that
		$$D(\alpha|\pi, \rho, F, \delta, \sigma)\ge \mathrm{Wdim}_{\varepsilon}(\mathrm{Map}(\rho_X, F, \delta, \sigma), \rho_{X,\infty}|Y),$$
		for any $F \in \mathcal{F}(G), \delta>0$ and any map $\sigma$ from $G$ to $\Sym(d)$ for some $d \in \N$.
		This implies that $\mathrm{mdim}_{\Sigma}(\alpha|\pi)\geq\mathrm{mdim}_{\Sigma}(\rho_X | Y)$. Hence $\mathrm{mdim}_{\Sigma}(X|\pi)\geq \mathrm{mdim}_{\Sigma}(\alpha|\pi)\geq\mathrm{mdim}_{\Sigma}(\rho_X | Y)$.
		
		Fix $\ep>0$ and pick a finite open cover $\alpha$ of $X$ consisting of some open sets of the diameter less than $\ep$ under $\rho_X$.	
		By Lemma~\ref{polyhe}, there exists a continuous map 	
		$g:\mathrm{Map}(\rho_X, F, \delta, \sigma)\ra Q$ for a compact  metric space $P$ with $\dim (Q)=D(\alpha|\pi, \rho, F, \delta, \sigma)$, such that 
		$$\{f^{-1}(p)\cap(\pi^d|_{\Map(\rho,F,\delta,\sigma)})^{-1}(y)\}_{(q,\tau)\in Q\times \pi^{d} (\mathrm{Map}(\rho_X, F, \delta, \sigma))}\succ\alpha^d|_{\mathrm{Map}(\rho_X, F, \delta, \sigma)}.$$
		Then $\diam(f^{-1}(p)\cap(\pi^d|_{\Map(\rho,F,\delta,\sigma)})^{-1}(\tau),\rho_{X,\infty})<\ep$, for any $q\in Q$, $\tau \in \pi^d(\mathrm{Map}(\rho_X, F, \delta, \sigma))$. In other words, $f$ is an $(\ep|Y,\rho_{X,\infty})$-embedding. Thus
		\begin{align*}
			D(\alpha|\pi, \rho, F, \delta, \sigma)&=\dim (Q)\\
			&\ge \mathrm{Wdim}_{\varepsilon}(\mathrm{Map}(\rho_X, F, \delta, \sigma_i), \rho_{X,\infty}|Y).
		\end{align*}
		The proof is complete.
	\end{proof}
	
	The sofic conditional  mean dimension  exhibits strong connections with factor maps. We now investigate its properties relative to different types of factor maps.
	
	\begin{prop}\label{isomorp}
		Let $\pi: (X,G)\to (Y,G)$ be a factor map between $G$-systems, $\Sigma$ a sofic approximation sequence for $G$ and $\psi:(Z, G)\to(Y,G)$  a factor map of $\pi$. Let  $\phi:(X, G)\to(Z,G)$ be a factor map with $\pi=\psi\circ\phi$. If $\phi$ is an isomorphism, then $\mathrm{mdim}_{\Sigma} (X|\pi)= \mathrm{mdim}_{\Sigma} (Z|\psi)$.
	\end{prop}
	\begin{proof}
		Let $F \in \mathcal{F}(G), \delta>0$, by Lemma \ref{li-2}, there exists $\delta'>0$, for every map $\sigma$ from $G$ to $\Sym(d)$ for some $d \in \N$, one has
		$$\phi^d(\Map(\rho_X,F,\delta',\sigma) \subset \Map(\rho_Z,F,\delta,\sigma),$$
		where $\phi^d:X^d \ra Y^d$ sends $\xi$ to $\phi \circ \xi$.
		By continuity of $(\phi^d)^{-1}$, for any $\varepsilon>0$, there exists $\varepsilon'>0$, such that for any $\tau,\tau' \in \Map(\rho_Z,F,\delta,\sigma)$ satisfying $\rho_{Z,\infty}(\tau, \tau') < \varepsilon'$, one has
		$$\rho_{Z,\infty}((\phi^d)^{-1}(\tau), (\phi^d)^{-1}(\tau'))=\rho_{X,\infty}(\phi^{-1} \circ\tau, \phi^{-1} \circ\tau')<\varepsilon.$$
		We show that 
		$$\mathrm{Wdim}_{\varepsilon}(\mathrm{Map}(\rho_X, F, \delta', \sigma), \rho_{X,\infty}|Y)\leq\mathrm{Wdim}_{\varepsilon'}(\mathrm{Map}(\rho_Z, F, \delta, \sigma), \rho_{Z,\infty}|Y).$$
		This implies that $\mathrm{mdim}_{\Sigma}(\rho_X | Y) \leq \mathrm{mdim}_{\Sigma}(\rho_Z | Y)$. By Proposition \ref{xoxoxoxo}, we have that $\mathrm{mdim}_{\Sigma}(X|\pi, G)\leq \mathrm{mdim}_{\Sigma}(Z|\psi, G)$.
		
		By the definition, there exists an $(\ep'|Y,\rho_{Z,\infty})$-embedding $f:\mathrm{Map}(\rho_Z, F, \delta, \sigma)\ra P$ such that for any $\tau,\tau'\in\mathrm{Map}(\rho_Z, F, \delta, \sigma)$ satisfying $f(\tau)=f(\tau')$ and $\psi\circ \tau=\psi\circ \tau'$, we have $\rho_{Z,\infty}(\tau,\tau')<\ep'$. Assume that $\mathrm{Wdim}_{\varepsilon'}(\mathrm{Map}(\rho_Z, F, \delta, \sigma), \rho_{Z,\infty}|Y)=\dim (P)$. 
		
		Consider the continuous map $f \circ \phi^d: \Map(\rho_X,F,\delta',\sigma) \ra P$ with
		$$f \circ \phi^d(\xi)=f(\phi \circ\xi).$$
		\begin{claim}\label{shjxsq}
			$f \circ \phi^d$ is an $(\varepsilon|Y,\rho_{X,\infty})$-embedding.
		\end{claim}
		
		By Claim \ref{shjxsq}, we have
		\begin{align*}
			\mathrm{Wdim}_{\varepsilon}(\mathrm{Map}(\rho_X, F, \delta', \sigma), \rho_{X,\infty}|Y)&\leq \dim(P)\\
			&=\mathrm{Wdim}_{\varepsilon'}(\mathrm{Map}(\rho_Z, F, \delta, \sigma), \rho_{Z,\infty}|Y),
		\end{align*}
		which implies the desired result.
		
		\emph{Proof of Claim \ref{shjxsq}}. For every $\xi,\xi' \in \Map(\rho_X,F,\delta',\sigma)$ satisfying $f \circ \phi^d(\xi)=f \circ \phi^d(\xi')$ and $\pi \circ \xi=\pi \circ \xi'$, we have 
		$$f \circ \phi^d(\xi)=f \circ \phi^d(\xi')=f(\phi \circ \xi)=f(\phi \circ \xi'),$$
		and 
		$$\pi \circ \xi=\pi \circ \xi'=\psi (\phi \circ \xi)=\psi (\phi \circ \xi'),$$
		where $\phi \circ \xi, \phi \circ \xi' \in \Map(\rho_Z,F,\delta,\sigma)$. Thus
		$$\rho_{Z,\infty}(\phi \circ \xi, \phi \circ \xi') < \varepsilon'.$$
		This implies that 
		$$\rho_{X,\infty}(\phi^{-1}\circ \phi \circ \xi, \phi^{-1}\circ \phi \circ \xi')=\rho_{X,\infty}( \xi, \xi')<\varepsilon.$$
		Therefore, $f \circ \phi^d$ is an $(\varepsilon|Y,\rho_{X,\infty})$-embedding.
		
		On the other hand, since $\phi$ is an isomorphism, by Lemma \ref{li-2}, there exists $\delta''>0$, such that
		\begin{equation}\label{dwdcew}
			(\phi^{-1})^{d}(\Map(\rho_Z,F,\delta'',\sigma)= (\phi^d)^{-1}(\Map(\rho_Z,F,\delta'',\sigma) \subset \Map(\rho_X,F,\delta,\sigma),
		\end{equation}
		for any $F \in \mathcal{F}(G), \delta>0$, and any map $\sigma$ from $G$ to $\Sym(d)$ for some $d \in \N$. By continuity of $\phi^d$, for any $\varepsilon>0$, there exists $\varepsilon''>0$, such that for any $\xi,\xi' \in \Map(\rho_X,F,\delta,\sigma)$ satisfying $\rho_{X,\infty}(\xi, \xi') < \varepsilon''$, one has
		$$\rho_{Z,\infty}(\phi^d(\xi), \phi^d(\xi'))=\rho_{Z,\infty}(\phi \circ \xi, \phi \circ \xi')<\varepsilon.$$
		
		Now we show that 
		$$\mathrm{Wdim}_{\varepsilon}(\mathrm{Map}(\rho_Z, F, \delta'', \sigma), \rho_{Z,\infty}|Y)\leq\mathrm{Wdim}_{\varepsilon''}(\mathrm{Map}(\rho_X, F, \delta, \sigma), \rho_{X,\infty}|Y).$$
		This implies that $\mathrm{mdim}_{\Sigma}(\rho_Z | Y) \leq \mathrm{mdim}_{\Sigma}(\rho_X | Y)$. By Proposition \ref{xoxoxoxo}, we have that $\mathrm{mdim}_{\Sigma}(Z|\psi, G)\leq \mathrm{mdim}_{\Sigma}(X|\pi, G)$.
		
		There exists an $(\ep''|Y,\rho_{X,\infty})$-embedding $g:\mathrm{Map}(\rho_X, F, \delta, \sigma)\ra Q$ such that for any $\xi,\xi'\in\mathrm{Map}(\rho_X, F, \delta, \sigma)$ satisfying $g(\xi)=g(\xi')$ and $\pi\circ \xi=\pi\circ \xi'$, we have $\rho_{X,\infty}(\xi,\xi')<\ep''$. Assume that $\mathrm{Wdim}_{\varepsilon''}(\mathrm{Map}(\rho_X, F, \delta, \sigma), \rho_{X,\infty}|Y)=\dim (Q)$.
		
		Consider the continuous map $g \circ (\phi^d)^{-1}: \Map(\rho_Z,F,\delta'',\sigma) \ra Q$ with
		$$g \circ (\phi^d)^{-1}(\tau)=g(\phi^{-1} \circ \tau).$$
		\begin{claim}\label{shjxsqqq}
			$g \circ (\phi^d)^{-1}$ is an $(\varepsilon|Y,\rho_{Z,\infty})$-embedding.
		\end{claim}
		
		By Claim \ref{shjxsqqq}, we have
		\begin{align*}
			\mathrm{Wdim}_{\varepsilon}(\mathrm{Map}(\rho_Z, F, \delta'', \sigma), \rho_{Z,\infty}|Y)&\leq \dim (Q)\\
			&=\mathrm{Wdim}_{\varepsilon''}(\mathrm{Map}(\rho_X, F, \delta, \sigma), \rho_{X,\infty}|Y).
		\end{align*}
		which implies the desired result.
		
		\emph{Proof of Claim \ref{shjxsqqq}}. For every $\tau,\tau' \in \Map(\rho_Z,F,\delta'',\sigma)$, there exist $\xi,\xi' \in \Map(\rho_X,F,\delta,\sigma)$ such that $\tau=\phi\circ\xi, \tau'=\phi \circ \xi'$ (By (\ref{dwdcew})). If
		$\tau=\phi\circ\xi, \tau'=\phi \circ \xi'$ satisfying $g \circ (\phi^d)^{-1}(\tau)=g \circ (\phi^d)^{-1}(\tau')$ and $\psi \circ \tau=\psi \circ \tau'$, then
		$$g \circ (\phi^d)^{-1}(\tau)=g \circ (\phi^d)^{-1}(\tau')=g (\phi^{-1} \circ\tau)=g (\phi^{-1} \circ\tau')=g(\xi)=g(\xi'),$$
		and 
		$$\psi \circ \tau=\psi \circ \tau'=\psi \circ \phi \circ \xi=\psi \circ \phi \circ \xi'=\pi \circ \xi=\pi \circ \xi'.$$
		Thus
		$$\rho_{X,\infty}(\xi,\xi') < \varepsilon''.$$
		This implies that 
		$$\rho_{Z,\infty}( \phi \circ \xi,  \phi \circ \xi')=\rho_{Z,\infty}( \tau, \tau')<\varepsilon.$$
		Therefore, $g \circ (\phi^d)^{-1}$ is an $(\varepsilon|Y,\rho_{Z,\infty})$-embedding.
	\end{proof}

	\begin{prop}\label{isomorphism}
		Let $\pi: (X,G)\to (Y,G)$ be a factor map between $G$-systems, $\Sigma$ a sofic approximation sequence for $G$ and $\psi:(Z, G)\to(Y,G)$  a factor map of $\pi$. Let  $\phi:(X, G)\to(Z,G)$ be a factor map with $\pi=\psi\circ\phi$. Then $\mathrm{mdim}_{\Sigma} (X|\pi) \geq \mathrm{mdim}_{\Sigma} (X|\phi)$. Moreover, if $\psi$ is an isomorphism, then $\mathrm{mdim}_{\Sigma} (X|\pi)= \mathrm{mdim}_{\Sigma} (X|\phi)$.
	\end{prop}
	\begin{proof}
		We first show that $$\mathrm{Wdim}_{\varepsilon}(\mathrm{Map}(\rho_X, F, \delta, \sigma), \rho_{X,\infty}|Z)\leq\mathrm{Wdim}_{\varepsilon}(\mathrm{Map}(\rho_X, F, \delta, \sigma), \rho_{X,\infty}|Y),$$
		for any $F \in \mathcal{F}(G), \delta>0,\varepsilon>0$, and any map $\sigma$ from $G$ to $\Sym(d)$ for some $d \in \N$. 
		This implies that $\mathrm{mdim}_{\Sigma}(\rho_X | Z) \leq \mathrm{mdim}_{\Sigma}(\rho_X | Y)$. By Proposition \ref{xoxoxoxo}, we have that $\mathrm{mdim}_{\Sigma}(X|\phi, G)\leq \mathrm{mdim}_{\Sigma}(X|\pi, G)$.
		
		By the definition, there exists an $(\ep|Y,\rho_{X,\infty})$-embedding $f:\mathrm{Map}(\rho_X, F, \delta, \sigma)\ra P$ such that for any $\xi,\xi'\in\mathrm{Map}(\rho_X, F, \delta, \sigma)$ satisfying $f(\xi)=f(\xi')$ and $\pi\circ \xi=\pi\circ \xi'$, we have $\rho_{X,\infty}(\xi,\xi')<\ep$. Assume that $\mathrm{Wdim}_{\varepsilon}(\mathrm{Map}(\rho_X, F, \delta, \sigma), \rho_{X,\infty}|Y)=\dim (P)$. Then we show that $f$ is an $(\ep|Z,\rho_{X,\infty})$-embedding.
		For any  $\xi,\xi'\in\mathrm{Map}(\rho_X, F, \delta, \sigma)$ satisfying $f(\xi)=f(\xi')$ and $\phi\circ \xi=\phi\circ \xi'$, one has $\psi\circ\phi\circ \xi=\psi\circ\phi\circ \xi'$, i.e. $\pi\circ \xi=\pi\circ \xi'$. Thus, $\rho_{X,\infty}(\xi,\xi')<\ep$, which implies that $f$ is an $(\ep|Z,\rho_{X,\infty})$-embedding.
		Hence,
		\begin{align*}
			\mathrm{Wdim}_{\varepsilon}(\mathrm{Map}(\rho_X, F, \delta, \sigma), \rho_{X,\infty}|Z)&\leq\dim (P)\\
			&=\mathrm{Wdim}_{\varepsilon}(\mathrm{Map}(\rho_X, F, \delta, \sigma), \rho_{X,\infty}|Y),
		\end{align*}
		which implies the desired result.
		
		Now assume that $\psi$ is an isomorphism. Then we show that $$\mathrm{Wdim}_{\varepsilon}(\mathrm{Map}(\rho_X, F, \delta, \sigma), \rho_{X,\infty}|Z)\ge\mathrm{Wdim}_{\varepsilon}(\mathrm{Map}(\rho_X, F, \delta, \sigma), \rho_{X,\infty}|Y),$$
		for any $F \in \mathcal{F}(G), \delta>0,\varepsilon>0$, and any map $\sigma$ from $G$ to $\Sym(d)$ for some $d \in \N$. 
		This implies that $\mathrm{mdim}_{\Sigma}(\rho_X | Z) \geq \mathrm{mdim}_{\Sigma}(\rho_X | Y)$. By Proposition \ref{xoxoxoxo}, we have that $\mathrm{mdim}_{\Sigma}(X|\pi)\geq \mathrm{mdim}_{\Sigma}(X|\phi)$.
		
		By the definition, there exists an $(\ep|Z,\rho_{X,\infty})$-embedding $g:\mathrm{Map}(\rho_X, F, \delta, \sigma)\ra Q$ such that for any $\xi,\xi'\in\mathrm{Map}(\rho_X, F, \delta, \sigma)$ satisfying $g(\xi)=g(\xi')$ and $\phi\circ \xi=\phi\circ \xi'$, we have $\rho_{X,\infty}(\xi,\xi')<\ep$. Assume that $\mathrm{Wdim}_{\varepsilon}(\mathrm{Map}(\rho_X, F, \delta, \sigma), \rho_{X,\infty}|Z)=\dim (Q)$. Then we show that $g$ is an $(\ep|Y,\rho_{X,\infty})$-embedding.
		
		Since $\pi=\psi\circ\phi$ and $\psi$ is an isomorphism, for any  $\xi,\xi'\in\mathrm{Map}(\rho_X, F, \delta, \sigma)$ satisfying $f(\xi)=f(\xi')$ and $\pi\circ \xi=\pi\circ \xi'$, one has $\phi\circ \xi=\phi\circ \xi'$. Thus, $\rho_{X,\infty}(\xi,\xi')<\ep$, which implies that $g$ is an $(\ep|Y,\rho_{X,\infty})$-embedding.
		Hence,
		\begin{align*}
			\mathrm{Wdim}_{\varepsilon}(\mathrm{Map}(\rho_X, F, \delta, \sigma), \rho_{X,\infty}|Z)&=\dim (Q)\\
			&\geq\mathrm{Wdim}_{\varepsilon}(\mathrm{Map}(\rho_X, F, \delta, \sigma), \rho_{X,\infty}|Y).
		\end{align*}
		This implies the desired result.
		
		By above discuss, we get
		$$\mathrm{mdim}_{\Sigma} (X|\pi) =\mathrm{mdim}_{\Sigma} (X|\phi).$$
		The proof is complete.
	\end{proof}

	\begin{prop}\label{subset}
		Let $\pi: (X,G)\to (Y,G)$ be a factor map between $G$-systems, $\Sigma$ a sofic approximation sequence for $G$ and \(K\) be a closed \(G\)-invariant subset of \(X\), $Y':=\pi(K)$ and $\pi_2:=\pi_1|_K:(K,G)\to(Y',G)$. Then $\mathrm{mdim}_{\Sigma} (K|\pi_2) \leq \mathrm{mdim}_{\Sigma} (X|\pi_1)$.
	\end{prop}
	\begin{proof}
		Suppose $\rho_X$ is a compatible metric defined on $X$.  This metric restricts to a compatible metric $\rho_K$ on the subset $K$. Let $\alpha$ be a finite open cover of \(K\). Then there exists a finite cover $\beta$ of $X$ such that $\alpha$ is exactly the covering obtained by intersecting $\beta$ with $K$. It is clear that $\mathrm{Map}(\rho_K, F, \delta, \sigma)\subset\mathrm{Map}(\rho_X, F, \delta, \sigma)$, for any $F \in \mathcal{F}(G)$, $\d>0$ and any map $\sigma$ from $G$ to $\Sym(d)$ for some $d\in\N$. Furthermore, the restriction of $\beta^{d}|_{\Map(\rho_X,F,\delta,\sigma)}$ on $\Map(\rho_K,F,\delta,\sigma)$ is exactly $\alpha^{d}|_{\Map(\rho_K,F,\delta,\sigma)}$. 
		
		Next we prove that
		\begin{equation}\label{jjhwdqoqq}
			D(\alpha^{d}|_{\Map(\rho_K,F,\delta,\sigma)}\big| \pi_2^d)\leq D(\beta^{d}|_{\Map(\rho_X,F,\delta,\sigma)}\big| \pi_1^d).
		\end{equation}
		Let $\gamma$ be a finite open cover of $\Map(\rho_X,F,\delta,\sigma)$, such that $$\{(\pi_1^d|_{\Map(\rho_X,F,\delta,\sigma)})^{-1}(y)\}_{y\in \pi_1^d({\Map(\rho_X,F,\delta,\sigma)})}\bigvee\gamma \succ \beta^d|_{\Map(\rho_X,F,\delta,\sigma)},$$
		and $D(\beta^d|_{\Map(\rho_X,F,\delta,\sigma)}|\pi_1^d)=\ord(\gamma|_{\Map(\rho_X,F,\delta,\sigma)})$.
		
		Then
		$$\{(\pi_1^d|_{\Map(\rho_X,F,\delta,\sigma)})^{-1}(y)\}_{y\in \pi_1^d({\Map(\rho_X,F,\delta,\sigma)})}\bigvee\gamma\big|_{\Map(\rho_K,F,\delta,\sigma)} \succ \alpha^d|_{\Map(\rho_K,F,\delta,\sigma)}.$$
		It's obvious that $\gamma|_{\Map(\rho_K,F,\delta,\sigma)}$ is a finite open cover of $\Map(\rho_K,F,\delta,\sigma)$ and 
		$$\{(\pi_2^d|_{\Map(\rho_K,F,\delta,\sigma)})^{-1}(y)\}_{y\in \pi_2^d({\Map(\rho_K,F,\delta,\sigma)})}\bigvee\gamma|_{\Map(\rho_K,F,\delta,\sigma)} \succ \alpha^d|_{\Map(\rho_K,F,\delta,\sigma)}.$$
		Hence,
		\begin{align*}
			D(\alpha^{d}|_{\Map(\rho_K,F,\delta,\sigma)}\big| \pi_2^d)&\leq\ord(\gamma|_{\Map(\rho_K,F,\delta,\sigma)})\\
			&\leq\ord(\gamma|_{\Map(\rho_X,F,\delta,\sigma)})= D(\beta^{d}|_{\Map(\rho,F,\delta,\sigma)}\big| \pi_1^d).
		\end{align*}
		Thus, by (\ref{jjhwdqoqq}), one has
		$$D(\alpha|\pi_2,\rho_{K},F,\delta,\sigma) \leq D(\alpha|\pi_1,\rho_{X},F,\delta,\sigma).$$
		This implies that
		$$\mathrm{mdim}_{\Sigma} (\alpha|\pi_2)\leq\mathrm{mdim}_{\Sigma} (\beta|\pi_1)\leq\mathrm{mdim}_{\Sigma} (X|\pi_1).$$	
		Since $\alpha$ is an arbitrary finite cover open cover of $K$, we have	
		$$\mathrm{mdim}_{\Sigma}(K|\pi_2)\leq\mathrm{mdim}_{\Sigma}(X|\pi_1).$$
		The proof is complete.
	\end{proof}
	
	Next we will discuss other fundamental properties of the sofic conditional mean dimension.
	
	\begin{prop}
		Let $\pi_n: (X_n,G)\to (Y_n,G)$ be a factor map between $G$-systems for each \(1 \leq n < R\), where \(R \in \mathbb{N} \cup \{\infty\}\) and $\Sigma$ a sofic approximation sequence for $G$. Consider the factor map $\pi:=\prod_{1 \leq n < R}\pi_n$ between the product systems $(X := \prod_{1 \leq n < R}X_n,G)$ and $(Y := \prod_{1 \leq n < R}Y_n,G)$. Then
		\begin{equation*}
			\mathrm{mdim}_{\Sigma} (X|\pi) \leq \sum_{1 \leq n < R} \mathrm{mdim}_{\Sigma} (X_n|\pi_n).
		\end{equation*}
	\end{prop}
	
	\begin{proof}
		For any $n\in \N$, let $\rho^{(n)}_X$ be compatible metrics on $X_n$ and denote by $p_n$ the projection of $X$ onto $X_n$ and $p'_n$ the projection of $Y$ onto $Y_n$.
		Let $\rho_X$ be compatible metrics on $X$  and $\alpha$  a finite open cover of $X$. Then there exist an $N \in \mathbb{N}$ with $N < R$ and a finite open covers $\beta_n$ of $X_n$ for all $1 \leq n \leq N$ such that
		$$\beta:=\bigvee_{n=1}^{N}p^{-1}_{n}(\beta_{n}) \succ\alpha.$$
		
		Let $F \in \mathcal{F}(G),\delta>0$, by Lemma \ref{li-2}, we can find a $\delta'>0$, such that for any map $\sigma$ from $G$ to $\Sym(d)$ for some $d \in \N$, one has
		\begin{equation}\label{qqwsswqwww}
			p^{d}_{n}(\Map(\rho_X,F,\delta',\sigma)) \subset \Map(\rho_{X}^{(n)},F,\delta,\sigma),
		\end{equation}
		
		for all $1 \leq n \leq N$. It follows that we have a continuous map $\Phi_n :=p^{d}_{n}|_{\Map(\rho_X,F,\delta',\sigma)}:\Map(\rho_X,F,\delta',\sigma) \rightarrow p^{d}_{n}(\Map(\rho_X,F,\delta',\sigma))$ sending $\xi$ to $p_{n} \circ \xi$ for each $1 \leq n \leq N$. Note that
		$$\beta^{d}|_{\Map(\rho_X,F,\delta',\sigma)}=\bigvee_{n=1}^{N}\Phi_n^{-1}(\beta_{n}^{d}|_{p^{d}_{n}(\Map(\rho_X,F,\delta',\sigma))}).
		$$
		Hence
		\begin{align*}
			D(\alpha|\pi,\rho_X,F,\delta',\sigma) &\leq D(\beta|\pi,\rho_{X},F,\delta',\sigma)\\
			&=D(\bigvee_{n=1}^{N}\Phi_n^{-1}(\beta_{n}^{d}|_{p^{d}_{n}(\Map(\rho_X,F,\delta',\sigma))})|\pi^d)\\
			&\leq \sum^{N}_{n=1}D(\Phi_n^{-1}(\beta_{n}^{d}|_{p^{d}_{n}(\Map(\rho_X,F,\delta',\sigma))})|\pi^d) (~\text{by~ Proposition ~\ref{subbb}}).
		\end{align*}
		
		\begin{claim}\label{nannanana}
			For every $1\leq n \leq N$, we have
			$$D(\Phi_n^{-1}(\beta_{n}^{d}|_{p^{d}_{n}(\Map(\rho_X,F,\delta',\sigma))})|\pi^d) \leq D((\beta_{n}^{d}|_{p^{d}_{n}(\Map(\rho_X,F,\delta',\sigma))})|\pi_n^d).
			$$
		\end{claim}
		
		By Claim \ref{nannanana}, we get
		\begin{align*}
			D(\alpha|\pi,\rho_X,F,\delta',\sigma) &\leq D(\beta|\pi,\rho_{X},F,\delta',\sigma)\\
			&\leq \sum^{N}_{n=1}D(\Phi_n^{-1}(\beta_{n}^{d}|_{p^{d}_{n}(\Map(\rho_X,F,\delta',\sigma))}|\pi^{d})\\
			&\leq \sum^{N}_{n=1}D((\beta_{n}^{d}|_{p^{d}_{n}(\Map(\rho_X,F,\delta',\sigma))}|\pi^d_n)(~\text{by~ Claim ~\ref{nannanana}})\\
			&\leq \sum^{N}_{n=1}D((\beta_{n}^{d}|_{(\Map(\rho_X^{(n)},F,\delta,\sigma))}|\pi^d_n)(~\text{by~ (\ref{qqwsswqwww}))}\\
			&=\sum^{N}_{n=1}D(\beta_{n}|\pi_{n},\rho^{(n)}_X,F,\delta,\sigma).
		\end{align*}
		Then 
		$$\mathrm{mdim}_{\Sigma}(\alpha|\pi,\rho_{X},F,\delta')\leq \sum_{n=1}^{N}\mathrm{mdim}_{\Sigma}(\beta_{n}|\pi_n,\rho^{(n)}_{X},F,\delta),$$
		which implies that
		$$\mathrm{mdim}_{\Sigma}(\alpha|\pi)\leq \sum^{N}_{n=1}\mathrm{mdim}_{\Sigma}(\beta_n|\pi_n)\leq \sum^{N}_{n=1}\mathrm{mdim}_{\Sigma}(X_n|\pi_n).$$ 
		Since $\alpha$ is arbitrary, we have
		$$\mathrm{mdim}_{\Sigma}(X|\pi)\leq \sum^{N}_{n=1}\mathrm{mdim}_{\Sigma}(X_n|\pi_n) \leq \sum_{1\leq n <R}\mathrm{mdim}_{\Sigma}(X_n|\pi_n),$$
		which implies the desired result.
		
		\emph{Proof of Claim \ref{nannanana}}. Since for any $1\leq n \leq N$, we know that
		$$\Phi_n:\Map(\rho_{X},F,\delta',\sigma) \ra p^{d}_{n}(\Map(\rho_{X}^{(n)},F,\delta',\sigma))$$
		is a continuous surjective map, and $$\Psi_n:=\pi^{d}_{n}|_{p^{d}_{n}(\Map(\rho_{X},F,\delta',\sigma))}:p^{d}_{n}(\Map(\rho_{X},F,\delta',\sigma)) \ra \pi^d_n \circ p^d_n(\Map(\rho_{X},F,\delta',\sigma))$$ is a continuous surjective map. By Lemma \ref{lemlk}, we get
		\begin{align*}
			D(\Phi_n^{-1}(\beta_{n}^{d}|_{p^{d}_{n}(\Map(\rho_X,F,\delta',\sigma))})|\Psi_n \circ \Phi_n) 
			&\leq D(\beta_{n}^{d}|_{p^{d}_{n}(\Map(\rho_X,F,\delta',\sigma))}|\Psi_n))\\
			&=D(\beta_{n}^{d}|_{p^{d}_{n}(\Map(\rho_X,F,\delta',\sigma))}|\pi_n^d).
		\end{align*}
		
		Note that
		\begin{align*}
			D(\Phi_n^{-1}(\beta_{n}^{d}|_{p^{d}_{n}(\Map(\rho_X,F,\delta',\sigma))})|\pi^d)&= D(\Phi_n^{-1}(\beta_{n}^{d})|_{\Map(\rho_X,F,\delta',\sigma)})|\pi^d)\\
			&= D(\Phi_n^{-1}(\beta_{n}^{d})|_{\Map(\rho_X,F,\delta',\sigma)})|\pi^d|_{\Map(\rho_X,F,\delta',\sigma)}))
		\end{align*}
		and
		$$
		D(\Phi_n^{-1}(\beta_{n}^{d}|_{p^{d}_{n}(\Map(\rho_X,F,\delta',\sigma))})|\Psi_n \circ \Phi_n)=D(\Phi_n^{-1}(\beta_{n}^{d})|_{\Map(\rho_X,F,\delta',\sigma)})|\Psi_n \circ \Phi_n). $$
		Thus we only need to show that
		$$D(\Phi_n^{-1}(\beta_{n}^{d})|_{\Map(\rho_X,F,\delta',\sigma)})|\pi^d|_{\Map(\rho_X,F,\delta',\sigma)}))\leq D(\Phi_n^{-1}(\beta_{n}^{d})|_{\Map(\rho_X,F,\delta',\sigma)})|\Psi_n \circ \Phi_n).$$

		Let $\Xi:=\pi^d |_{\Map(\rho_X,F,\delta',\sigma}),\Omega_n:=(p'_n)^{d} |_{\pi^{d}({\Map(\rho_X,F,\delta',\sigma}))}$, then $\Omega_n \circ \Xi =\Psi_n \circ \Phi_n$. Hence
		$$ (\Psi_n \circ \Phi_n)^{-1}(\tau^{(n)})=(\Omega_n \circ \Xi)^{-1}(\tau^{(n)}) \supseteq\Xi^{-1}(\tau),$$
		for every $\tau=(\tau^{(1)},\cdots,\tau^{(n)},\cdots) \in \pi^d(\Map(\rho_{X},F,\delta',\sigma)).$
		
		Suppose $\gamma$ is a finite open cover of $\Map(\rho_X,F,\delta',\sigma)$, such that
		$$\{(\Psi_n \circ \Phi_n)^{-1}(\tau^{(n)})\}_{\tau^{(n)} \in \pi_n^d \circ p_n^d((\Map(\rho_{X},F,\delta',\sigma))} \bigvee \gamma \succ \Phi_n^{-1}(\beta_{n}^{d})|_{\Map(\rho_{X},F,\delta',\sigma)}$$
		and $D(\Phi_n^{-1}(\beta_{n}^{d})|_{\Map(\rho_X,F,\delta',\sigma)})|\Psi_n \circ \Phi_n)=\ord(\gamma|\Map(\rho_{X},F,\delta',\sigma))$. Then
		$$\{\Xi^{-1}(\tau)\}_{\tau=(\tau^{(1)},\cdots,\tau^{(n)},\cdots) \in \pi^d(\Map(\rho_{X},F,\delta',\sigma))} \bigvee \gamma \succ \Phi_n^{-1}(\beta_{n}^{d})|_{\Map(\rho_{X},F,\delta',\sigma)}.
		$$
		This implies that 
		\begin{align*}
			D(\Phi_n^{-1}(\beta_{n}^{d}|_{p^{d}_{n}(\Map(\rho_X,F,\delta',\sigma))})|\pi^d) &= D(\Phi_n^{-1}(\beta_{n}^{d})|_{\Map(\rho_X,F,\delta',\sigma)})|\Xi)\\
			&\leq \ord(\gamma|\Map(\rho_{X},F,\delta',\sigma))\\
			&=D(\Phi_n^{-1}(\beta_{n}^{d})|_{\Map(\rho_X,F,\delta',\sigma)})|\Psi_n \circ \Phi_n).
		\end{align*}
		which completes the proof of Claim \ref{nannanana}.
	\end{proof}
	
	Factor maps with zero sofic conditional  mean dimension constitute a significant proportion, making their study particularly important. We now introduce a key concept about zero sofic conditional  mean dimension which named zero sofic conditional  mean dimension extension, and investigate its fundamental properties. 
	
	\begin{defn}
		Let $\pi: (X,G)\to (Y,G)$ be a factor map between $G$-systems and $\Sigma$ a sofic approximation sequence for $G$. We say that $\pi$ is a \emph{zero sofic conditional  mean dimension extension (with respect to \(\Sigma\))} if $\mathrm{mdim}_{\Sigma} (X|\pi)=0$.
	\end{defn}
	
	Recall that if \(\{\psi_\lambda : (X_\lambda,G) \to (Z,G)\}_{\lambda \in \Lambda}\) is a set of extensions of $(Z,G)$, then \emph{the product relative to $(Z,G)$} of this set is the extension $\psi: (X,G) \to (Z,G)\}$, where $(X,G)$ is the subsystem of \((\prod_{\lambda \in \Lambda} X_\lambda, G)\) on the closed and invariant set
	\[
	X := \left\{ (x_\lambda)_{\lambda} \in \prod_{\lambda \in \Lambda} X_\lambda : \forall \kappa,\lambda \in \Lambda, \psi_\kappa(x_\kappa) = \psi_\lambda(x_\lambda) \right\}
	\]
	\[
	= \left\{ (x_\lambda)_{\lambda} \in \prod_{\lambda \in \Lambda} X_\lambda : \exists z \in Z \ \forall \lambda \in \Lambda, \psi_\lambda(x_\lambda) = z \right\},
	\]
	and \(\psi : X \to Z\) is (unambiguously) defined as follows: if \(\mathbf{x} = (x_\lambda)_{\lambda} \in X
	\) then
	\[
	\psi(\mathbf{x}) := \psi_\kappa(x_\kappa) \quad \text{for some (every) } \kappa \in \Lambda.
	\]
	(See \cite[Appendix E.14-2]{Vries1993} for details.)
	\begin{cor}\label{prod}
		Let $\Sigma$ be a sofic approximation sequence for $G$ and \(\{\psi_n : (X_n,G) \to (Z,G)\}_{n \in \N}\) a set of zero sofic conditional mean dimension extensions of $(Z,G)$. Let $\psi: (X,G)\to (Z,G)$ be the product relative to $(Z,G)$ of this set. Then $\mathrm{mdim}_{\Sigma} (X|\psi) \leq 0$. 
	\end{cor}
	\begin{proof}
		By the definition, $$X= \left\{ (x_n)_{n\in\N} \in \prod_{n \in \N} X_n : \exists z \in Z \ \forall n \in \N, \psi_n(x_n) = z \right\}\subset\prod_{n\in\N}X_n,$$
		and $\psi(X)=\Delta_{\prod_{n\in\N}Z}:=\{(z)_{n\in\N}\in\prod_{n\in\N}Z:z\in Z\}$. By Proposition~\ref{subset} and Proposition~\ref{prod}, one has $\mathrm{mdim}_{\Sigma} (X|\prod_{n \in \N} \psi_n|_X) \leq 0$. It is clear that $(\Delta_{\prod_{n\in\N}Z},G)$ is conjugate to $(Z,G)$. Then by Proposition \ref{isomorphism}. $\mathrm{mdim}_{\Sigma} (X|\psi)=\mathrm{mdim}_{\Sigma} (X|\prod_{n \in \N} \psi_n|_X) \leq 0$.
	\end{proof}
	
	We need the following key lemma.
	
	\begin{lem}\label{factor}
		Let $\pi: (X,G)\to (Y,G)$ be a factor map between $G$-systems 
		with \(\mathrm{mdim}_{\Sigma} (X|\pi) \geq 0\) and $\Sigma$ a sofic approximation sequence for $G$. Then for any factor $\psi:(Z,G)\to(Y,G)$ of $\pi$, one has $\mathrm{mdim}_{\Sigma} (Z|\psi) \geq 0$.
	\end{lem}
	\begin{proof}
		Since \(\mathrm{mdim}_{\Sigma} (X|\pi) \geq 0\), for any finite subset $F$ of $G$, any $\d > 0$, and any $N\in \N$, there is some $i\ge N$ such that
		$\mathrm{Map}(\rho_X, F, \delta, \sigma_i)$ is nonempty. For any factor $\psi:(Z,G)\to(Y,G)$ of $\pi$, there exists a factor map $\phi:(X,G)\to(Z,G)$ with $\pi=\psi\circ\phi$. Then by Lemma~\ref{li-2}, there exists $\delta' > 0$ such that for every map $\sigma$ from $G$ to $\operatorname{Sym}(d)$ for some $d \in \mathbb{N}$ and every $\xi \in \operatorname{Map}(\rho_X, F, \delta', \sigma_i)$, one has $\phi \circ \xi \in \operatorname{Map}(\rho_Z, F, \delta, \sigma_i)$. Then $\mathrm{mdim}_{\Sigma} (Z|\psi) \geq 0.$
	\end{proof}

	Now we can present the first main result of this paper.

	\begin{thm}\label{main1}
		Let $\pi: (X,G)\to (Y,G)$ be a factor map between $G$-systems 
		with \(\mathrm{mdim}_{\Sigma} (X|\pi) \geq 0\) and $\Sigma$ a sofic approximation sequence for $G$. Then there is a maximal zero sofic conditional  mean dimension factor $\psi:(X_{\pi}^{\Sigma}, G)\to(Y,G)$ of $\pi$. That is, $\psi:(X_{\pi}^{\Sigma}, G)\ra (Y,G)$ is a zero sofic conditional  mean dimension extension
		and a factor of $\pi$;
		and every zero sofic conditional mean dimension  extension $\psi':(Z,G)\ra (Y,G)$ of $\pi$ is a factor of $\psi$.
	\end{thm}
	
	\begin{proof}
		For each factor $\psi_W:(W,G)\ra (Y,G)$ of $\pi$, denote by $R_W$ the closed subset
		\[
		R_W = \left\{ (x, y) \in R_\pi : \psi_W(x) = \psi_W(y) \right\} \subset R_\pi.
		\]
		Denote by \(R\) the 
		set \(\bigcap_{W} R_W\) where the intersection ranges over zero sofic conditional  mean dimension factor $\psi_W$ of $\pi$ (i.e., \(\mathrm{mdim}_{\Sigma}(W|\psi_W) = 0
		\)). 
		
		Since $R_\pi$ is a compact  metric space, it has a countable base. Then for any subset $K$ of $R_\pi$ with the topology inherited from $R_\pi$, $K$ is a Lindelöf space. That is, every open cover of $K$ has a countable subcover (see for example \cite[Page 49]{JLK}). 
		
		Taking \(K = X^
		2 \setminus R\) and considering the open cover of \( R_\pi \setminus R\) consisting of \( R_\pi \setminus R_W\) for all zero sofic conditional mean dimension factors $\psi_W$ of  $\pi$ , we find zero sofic conditional mean dimension factors \(W_1, W_2, \dots\) of  $\pi$
		such that	$\bigcap_{n=
			1
		}^{\infty} R_{W_n} = R$.
		Since $\psi_{W_n}:(W_n,G)\ra (Y,G)$ is a factor map of $\pi$, there exists a factor $\phi_{W_n}:(X,G)\ra (W_n,G)$ such that $\pi=\phi_{W_n}\circ\psi_{W_n}$, for any $n\in\N$.

		Consider the 
		map \(\phi : X \to \prod _{n=1}^{\infty}W_n\) sending \(x\) to \(\left( \phi_{W_n}(x) \right)_{n=1
		}^{\infty}\). Let $X_{\pi}^{\Sigma} := \phi(X)$. Then
		$(X_{\pi}^{\Sigma},G)$ is a $G$-system and $R_\phi=R\subset R_\pi$. Thus there exists a factor map $\psi:(X_{\pi}^{\Sigma},G)\ra(Y,G)$ with $\pi=\phi\circ\psi$.
		Now we prove that $\psi$ is the maximal zero sofic conditional mean dimension factor of $\pi$.
		
		Let $\Psi:(Z,G)\ra (Y,G)$ be the product relative to $(Y,G)$ of $\{\psi_{W_n}:(W_n,G)\ra(Y,G) \}_{n\in\N}$. By the definition of $X_{\pi}^{\Sigma}$, $Z$ and $\psi$, $X_{\pi}^{\Sigma}\subset Z$ and $\Psi|_{X_{\pi}^{\Sigma}}=\psi$. Then by Proposition~\ref{subset} and Corollary~\ref{prod}, $\mathrm{mdim}_{\Sigma}(X_{\pi}^{\Sigma}|\psi) \leq 0$. As $\mathrm{mdim}_{\Sigma} (X|\pi) \geq 0$, combining with Lemma~\ref{factor}, we have $\mathrm{mdim}_{\Sigma}(X_{\pi}^{\Sigma}|\psi) = 0$ .
		
		Since 
		$R_W \supset R = R_\phi$
		for every factor zero sofic conditional mean dimension factor $\psi_W$ of $\pi$, $\psi$ is the maximal zero sofic conditional mean dimension factor of $\pi$. The proof is complete.
	\end{proof}
	\begin{rem}
		In the proof of Theorem~\ref{main1}, let $S_\pi^{\Sigma}=R_\phi$. Then $X_{\pi}^{\Sigma}=X/S_\pi^{\Sigma}$. $S_\pi^{\Sigma}$ is said to be the \emph{relative zero sofic conditional mean dimension relation}. It is clear that
		\begin{enumerate}
			\item [(1)] $S_\pi^{\Sigma}=\Delta_X$ if and only if $\pi$ is a zero sofic conditional mean dimension extension (with respect to \(\Sigma\)).
			\item [(2)] $S_\pi^{\Sigma}=R_\pi$ if and only if $\pi$ contains no proper factor.
		\end{enumerate}
	\end{rem}

	\section{Sofic conditional mean dimension tuples, sofic \texorpdfstring{$X|\pi$-}-CPCMD and sofic \texorpdfstring{$X|\pi$-}-UPCMD}\label{sec4}
	
	In this section, we study the local properties of sofic conditional mean dimension. Firstly, we introduce the notions of sofic conditional mean dimension tuples. Then we define completely positive sofic conditional mean dimension and uniform positive sofic conditional mean dimension, respectively, and study the properties of them.
	
	\subsection{Sofic conditional mean dimension tuples}\label{sbsec1}

	\begin{defn}
		Let $(x_i)_{i=1}^n \in X^n \setminus \Delta_n(X)$ and $\alpha$ a finite open cover of $X$. We say $\alpha$ is an \emph{admissible open cover} with respect to $(x_i)_{i=1}^n$ if for any $U \in \alpha$, $\{x_1, x_2, \ldots, x_n\} \not\subset \overline{U}$. A finite open cover $\alpha= \{U_1, U_2,\dots,U_n\}$ of $X$ is called \textit{non-dense} if  $\ov{U_i}\neq X$, for $i=1, 2,\dots,n$. 
	\end{defn}

	\begin{defn}
		Let $\pi \colon (X, G) \to (Y, G)$ be a nontrivial factor map between $G$-systems and $\Sigma$ a sofic approximation sequence for $G$. A tuple $(x_i)_{i=1}^n\in X^n$ is said to be a \emph{sofic conditional mean dimension tuple relevant to $\Sigma$ and $\pi$} if for every admissible open cover $\alpha$ with respect to $(x_i)_{i=1}^n$, we have $\text{mdim}_{\Sigma}(\alpha|\pi)>0$.
	\end{defn}
	
	For $n \geq 2$, denote by $\mathrm{D}_{n}^{\mathrm{md}}(X|\pi, G,\Sigma)$ the set of all sofic conditional mean dimension $n$-tuples  relevant to $\Sigma$ and $\pi$, and write it as $\mathrm{D}_{n}^{\mathrm{md}}(X, G,\Sigma)$ (the set of all sofic  mean dimension $n$-tuples), when $(Y,G)$ is trivial (see for example \cite{GRFYG} for $n=2$). 
	
	The following result is to make full use of the property of admissible cover.

	\begin{prop}\label{admis}
		Let $\pi:(X,G)\to (Y,G)$ be a factor map between $G$-systems, $n \geq 2$ and $\Sigma$ a sofic approximation sequence for $G$. Then 
		$(x_i)_{i=1}^n \in X^n \setminus \Delta_n(X)$ is a sofic conditional mean dimension tuple relevant to $\Sigma$ and $\pi$ if and only if for all open covers with the form $\alpha= (U_1, U_2,\dots,U_n)$, where $U_i^{c}$ is a neighborhood of $x_i$ such that $U_i^{c}\cap U_j^{c}=\emptyset$ when $x_i\neq x_j$, $1\leq i<j\leq n$, we have $\mathrm{mdim}_{\Sigma}(\alpha|\pi)>0$.
	\end{prop}
	
	\begin{proof}
		Assume that $\beta= \{V_1, V_2,\dots,V_l\}$ is an admissible open finite cover of $X$ with respect to $(x_i)_{i=1}^n$. Without loss of generality, let $\{x_1,x_2,\dots,x_k\}$ be the set of different points in the  set $\{x_1,x_2,\dots,x_n\}$. Then there exists a small enough $\ep>0$ such that for every $i\in\{1,2,\dots,l\}$, one can find $j_i\in\{1,2,\dots,k\}$ satisfying $B_\ep(x_{j_i})\subset \ov{V_{i}}^c$ and  $\overline{B_\ep(x_i)}\cap\overline{B_\ep(x_j)}=\emptyset$, for any  $i\neq j$, $i,j\in\{1,\dots,k\}$. Let $U_i=\overline{B_\ep(x_i)}^c$, $i=1,2,\dots,k$.
		Since $\bigcap_{i=1}^k\overline{B_{\varepsilon}(x_i)}=\emptyset,$ we know that $\bigcup_{i=1}^k\overline{B_{\varepsilon}(x_i)}^c=\bigcup_{i=1}^kU_i=X.$ Therefore,
		$$\text{mdim}_{\Sigma}((U_1,\dots,U_k)|\pi)>0.$$
		Since $(V_1,\dots,V_l)\succ(U_1,\dots,U_k)$, then 
		$$0<\text{mdim}_{\Sigma}((U_1,\dots,U_k)|\pi)\leq \text{mdim}_{\Sigma}((V_1,\dots,V_l)|\pi),$$
		which implies that $(x_i)_{i=1}^n $ is a sofic conditional mean dimension tuple relevant to $\Sigma$ and $\pi$.
		
		On the other hand, assume that $(x_i)_{i=1}^n \in X^n \setminus \Delta_n(X)$ is a sofic conditional mean dimension tuple relevant to $\Sigma$ and $\pi$. Let $\mathcal{W}=(W_1,\dots,W_n)$ be a finite open cover of $X$ satisfying $W_i^c$ is a neighborhood of $x_i$ such that $W_i^c\cap W_j^c=\emptyset,$ where $x_i\neq x_j,1\leq i<j\leq n.$ It's obvious that $\mathcal{W}=(W_1,\dots,W_n)$ is an admissible open cover of $X$ with respect to $(x_i)_{i=1}^n$ (If it's not true, we can find $j\in \{1,\dots,n\}$ such that $\{x_1,\dots,x_n\}\subseteq \overline{W_j}$. Since $x_j\in W_j^c$, then $x_j\in \overline{W_j}\cap W_j^c$, a contradiction). Thus $\text{mdim}_{\Sigma}(\mathcal{W}|\pi)>0.$
	\end{proof}
	Next we will discuss the properties of sofic conditional mean dimension tuples.
	\begin{prop}\label{subset Rpi}
		Let $\pi:(X,G)\to (Y,G)$ be a factor map between $G$-systems, $n \geq 2$ and $\Sigma$ a sofic approximation sequence for $G$. Then 
		$$\mathrm{D}_{n}^{\mathrm{md}}(X|\pi, G,\Sigma)\subset R_\pi^n\setminus \Delta_n(X).$$
	\end{prop}
	
	\begin{proof}
		It suffices to show that for any $(x_i)_{i=1}^n\in X^n\setminus R_\pi^n$, there exists an admissible open cover $\alpha$ with respect to $(x_i)_{i=1}^n$ such that $\mdim_{\Sigma}(\alpha|\pi)\leq 0$.
		
		There exists an admissible open cover $\alpha=(\pi^{-1}(U_1),\dots,\pi^{-1}(U_n))$ with respect to $(x_i)_{i=1}^n$, where $U_1,\dots,U_n$ are some open sets of $Y$. Then we show that $D(\alpha|\pi)=0$.
		
		Since 	$$	D(\alpha|\pi) = \min_{\{\pi^{-1}(y)\}_{y\in Y}\vee\beta \succ \alpha} \mathrm{ord}(\beta),$$
		and $\{\pi^{-1}(y)\}_{y\in Y}\succ(\pi^{-1}(U_1),\dots,\pi^{-1}(U_n))$,  we can choose $\beta=\{X,\emptyset\}$ such that $\{\pi^{-1}(y)\}_{y\in Y}\vee\beta \succ \alpha$, and then $D(\alpha|\pi)=\mathrm{ord}(\beta)\leq 0$. Therefore, by Proposition \ref{lllplpl}, we have
		$\mdim_{\Sigma}(\alpha|\pi)\leq D(\alpha|\pi) \leq0$. The proof is complete.
	\end{proof}

	\begin{prop}\label{propo322}
		Let $\pi:(X,G)\to (Y,G)$ be a factor map between $G$-systems, $n \geq 2$ and $\Sigma$ a sofic approximation sequence for $G$. If $\mathrm{mdim}_{\Sigma}(X|\pi)>0$, then there exists an admissible open cover $\alpha=(U_1,U_2,\dots,U_n)$ of $X$ with respect to some $(x_i)_{i=1}^n\in X^n$ such that $\mathrm{mdim}_{\Sigma}(\alpha|\pi)>0.$
	\end{prop}
	\begin{proof}
		Since $\mathrm{mdim}_{\Sigma}(X|\pi)>0$, one can take a finite open cover $\alpha=(U_1,U_2,\dots,U_k)$ of $X$ such that $\text{mdim}_{\Sigma}(\alpha|\pi)>0.$ Without loss of generality, we can assume that for every $i\in \{1,2,\dots,k\}$, $U_i$ is not contained in the union of the other $U_j$'s. That is, for every $i\in \{1,2,\dots,k\}$
		\begin{equation}\label{equa111}
			U_i\setminus \bigcup_{j\in \{1,2,\dots,k\}\setminus \{i\}}U_j\neq \emptyset.
		\end{equation}

		For each $i\in \{1,\dots,k\}$ and each $x\in U_i$, choose an open ball $B^i(x)$ such that $\overline{B^i(x)}\subset U_i$. Let $(W_1,\dots,W_p)$ be a finite subcover of the collection $\{B^i(x)\}_{i\in \{1,\dots,k\}, x\in U_i}$. 
		
		By (\ref{equa111}), for each $i\in \{1,\dots,k\}$, there exists some $j\in \{1,\dots,p\}$ with $\overline{W_j}\subset U_i$. For each $i$, define
		\[
		Z_i = \{\overline{W_j}^c : U_i^c \subseteq \overline{W_j}^c\},
		\]
		which is nonempty. Moreover, when $U_i^c\neq \emptyset$, we have
		\[
		U_i^c \subset V_i := \bigcap_{W\in Z_i} W \neq \emptyset.
		\]
		
		Now observe that $\beta_i = (U_i,V_i)$ forms an open cover of $X$ for each $i\in \{1,\dots,k\}$. Note that every $\overline{W_j}^c$ appears in some $Z_i$ (specifically, if $\overline{W_j}\subset U_i$, then $\overline{W_j}^c\in Z_i$). This implies
		\[
		\bigcup_{i=1}^k \bigcup_{W\in Z_i} W^c = X,
		\]
		and consequently,
		\[
		\bigcap_{i=1}^k \bigcap_{W\in Z_i}W = \bigcap_{i=1}^k V_i = \emptyset.
		\]
		
		Consider now the join $\bigvee_{i=1}^k \beta_i$. Any nonempty member of this join must be contained in some $U_i$, since the intersection $\bigcap_{i=1}^k V_i$ is empty. This shows that $\bigvee_{i=1}^k \beta_i$ refines $\alpha$.	
		
		Then, by Proposition \ref{impor}, we get
		$$0< \text{mdim}_{\Sigma}(\alpha|\pi)\leq \text{mdim}_{\Sigma}(\bigvee_{i=1}^k\beta_i|\pi)\leq \sum_{i=1}^k\text{mdim}_{\Sigma}(\beta_i|\pi).$$
		Hence there exists $i\in \{1,\dots,k\}$ such that $\text{mdim}_{\Sigma}(\beta_i|\pi)>0.$ If $U_i$ and $V_i$ are not dense, then $\beta_i$ is a non-dense open cover, there exists $(x_1,x_2)\in X\times X$ such that $\beta_i$ distinguishes $(x_1,x_2)$. Assume that $U_i$ is dense, there exists $x_1\in U_i$ such that $\overline{B_{\varepsilon}(x_1)}\subset V_i$. It holds that $\mathcal{V}=(U_i\setminus \overline{B_{\varepsilon}(x_1)},V_i)$ is also an open cover of $X$. Since $\mathcal{V}\succ \beta,$ then $\text{mdim}_{\Sigma}(\mathcal{V}|\pi)>0.$ If $V_n$ is not dense, then $\mathcal{V}$ is a non-dense open cover, we can find a $x_2\in X$ such that $\mathcal{V}$ distinguishes $(x_1,x_2)$. If $V_n$ is dense, then we can find $x_2\in V_i$ such that $\overline{B_{\varepsilon'}(x_2)}\subset U_i\setminus \overline{B_{\varepsilon}(x_1)}$
		It holds that $\mathcal{V}':=(U_i\setminus \overline{B_{\varepsilon}(x_1)},V_i\setminus \overline{B_{\varepsilon'}(x_2)})$ is a non-dense open cover. Moreover, since $\mathcal{V}'\succ \mathcal{V} $, we have $\text{mdim}_{\Sigma}(\mathcal{V}'|\pi)>0$, where $\mathcal{V}'$ distinguishes $(x_1,x_2)$.
		
		By following the above steps, we always find a non-dense open cover $\mathcal{M}=(U,V)$ distinguishes $(x_1,x_2)\in X\times X$, and $\text{mdim}_{\Sigma}(\mathcal{M}|\pi)>0$.
		
		Choose $x_3,x_4,\dots,x_n\in U$. Consider the open balls $\{B_{\varepsilon^{(i)}}(x_i)\}_{i=3}^{n}$, it's easy to know that $\mathcal{N}=(U,V, B_{\varepsilon^{(2)}}(x_2),\dots,B_{\varepsilon^{(n)}}(x_n))$ is an open cover of $X$. Moreover, $\{x_1,x_2,\dots,x_n\}\not\subset \overline{U}$, $\{x_1,x_2,\dots,x_n\}\not\subset \overline{V}$ and $\{x_1,x_2,\dots,x_n\}\not\subset B_{\varepsilon^{(i)}}(x_i)$ for every $i\in \{3,4,\dots,n\}$, which implies that $\mathcal{N}=(U,V, B_{\varepsilon^{(2)}}(x_2),\dots,B_{\varepsilon^{(n)}}(x_n))$ is an admissible open cover of $X$ with $(x_1,x_2,\dots,x_n)\in X^{n}$. Since $\mathcal{N}\succ \mathcal{M},$ then $\text{mdim}_{\Sigma}(\mathcal{N}|\pi)>0.$
	\end{proof}
	
	\begin{prop}\label{ppppro324}
		Let $\pi:(X,G)\to (Y,G)$ be a factor map between $G$-systems, $n \geq 2$ and $\Sigma$ a sofic approximation sequence for $G$. For any finite open cover $\alpha=(U_1,U_2,\dots,U_n)$ of $X$ satisfying $\mathrm{mdim}_{\Sigma}(\alpha|\pi)>0$, there exist $x_i\in U_i^c, i=1,\dots,n$ such that $(x_i)_{i=1}^n\in \mathrm{D}_{n}^{\mathrm{md}}(X|\pi,G,\Sigma).$ 
	\end{prop}
	\begin{proof}
		We first show that there exists an open cover $\beta_1=(B_1^{1},B_2^{1},\dots,B_n^{1})$ of $X$ satisfying $\text{mdim}_{\Sigma}(\beta_1|\pi)>0$, and $B_i^{1}\supseteq U_i$ with $\text{diam}((B_i^{1})^c)\leq \frac12\text{diam}((U_i)^c)$ for $i=1,2,\dots,n.$ By induction, for every $i\in \{1,2,\dots,n\}$, one gets a strictly decreasing sequence of nonempty closed sets $(B_i^{m})^c$ converging to $x_i$ for all $i=1,2,\dots,n$, and  $(x_i)_{i=1}^{n}\in \mathrm{D}_{n}^{\mathrm{md}}(X|\pi,G,\Sigma).$
		
		If $U_i^c$ is a singleton, let $B_1^{1}=U_1.$ If $U_i^c$ is not a singleton, there exist at least two distinct points $y,y'\in U_1^{c}$ and $\rho(y,y')>0$. Fix $\varepsilon_1$ such that $0<\varepsilon_1<\frac12\rho(y,y')$, and constructs a cover of $U_1^c$ by open balls with radius $\varepsilon_1$ centered in $U_1^c$, call it $\mathcal{V}.$ Since $U_1^c$ is a compact set, there exists a finite subcover $(W_1,\dots,W_k)$ of $\mathcal{V}.$ It's obvious that $k\geq 2.$ Let $F_i=U_1^c\cap \overline{W_i}$, by the choice of $\mathcal{V}$, we know that $F_i$ is a proper subset of $U_1^c$. 
		
		Since $\bigvee_{i=1}^k(F_i^c,U_2,\dots,U_n)\succ \alpha$, then by Proposition \ref{impor}, we have $$0<\text{mdim}_{\Sigma}(\alpha|\pi)\leq \text{mdim}_{\Sigma}(\bigvee_{i=1}^k(F_i^c,U_2,\dots,U_n)|\pi)\leq \sum_{i=1}^k\text{mdim}_{\Sigma}((F_i^c,U_2,\dots,U_n)|\pi).$$
		Hence, there exists $i\in \{1,2,\dots,k\}$ such that $\text{mdim}_{\Sigma}((F_i^c,U_2,\dots,U_n)|\pi)>0.$ Put $B_1^{1}=F_i^c$, do the same for $U_2$, thus obtaining $B_2^{1}$, in sequence, construct $B_3^{1},\dots,B_n^{1}$, such that the open cover $\beta_1=(B_1^{1},\dots,B_n^{1})$ of $X$ satisfying $\text{mdim}_{\Sigma}(\beta_1|\pi)>0$ and $B_i^{1}\supseteq U_i$ with $\text{diam}((B_i^{1})^c)\leq \frac12\text{diam}((U_i)^c), i=1,2,\dots,n.$
		
		Iterate infinitely times, we get a strictly decreasing sequence of nonempty sets $(B_i^{(m)})^c$ converging to $x_i, i=1,2,\dots,n$ and  $$\text{mdim}_{\Sigma}(\beta_m|\pi)>0,$$
		where $\beta_m=(B_1^{m},B_2^{m},\dots,B_n^{m}).$
		Since $\alpha=(U_1,U_2,\dots,U_n)$ is a cover of $X$,  $\bigcap_{i=1}^nU_i^c=\emptyset$. As $x_i\in U_i^c$, for any $i$, then $(x_i)_{i=1}^n\notin\Delta_n(X).$
		
		Next we show that $(x_i)^{n}_{i=1}\in \mathrm{D}_{n}^{\mathrm{md}}(X|\pi,G,\Sigma).$  Given an admissible open cover $\mathcal{V}=(V_1,\dots,V_l)$ of $X$ with respect to $(x_i)_{i=1}^n.$ There exists $\varepsilon>0$ such that for every $i\in \{1,2,\dots,l\}$, we can find $j_i\in \{1,2,\dots,n\}$ satisfies $B_{\varepsilon}(x_{j_i})\subset V_i^c$. Choose a sufficiently large $m\in \N$ such that $(B_i^{m})^c\subset B_{\varepsilon}(x_{i})$ for every $i\in \{1,2\dots,n\}$. Then for every $i\in \{1,\dots,n\}$, we have $V_i\subset (B_{\varepsilon}(x_{j_i}))^c\subset B_{j_i}^{m}$. Hence $\mathcal{V}\succ \beta_m$, by Proposition \ref{impor}, which implies that
		$$0<\text{mdim}_{\Sigma}(\beta_m|\pi)\leq \text{mdim}_{\Sigma}(\mathcal{V}|\pi).$$
		This proof is complete.
	\end{proof}
	
	Then we can give a equivalent characterization of factor maps with positive sofic conditional mean dimension via sofic conditional mean dimension tuples.

	\begin{thm}\label{main3}
		Let $\pi:(X,G)\to (Y,G)$ be a factor map between $G$-systems and $\Sigma$ a sofic approximation sequence for $G$. Then $\mathrm{mdim}_{\Sigma}(X|\pi)>0$ if and only if $D_{n}^{\mathrm{md}}(X|\pi, G,\Sigma)\neq \emptyset$ for some $n\ge 2$.
	\end{thm}
	\begin{proof}
		Assuming that $\mathrm{D}_{n}^{\mathrm{md}}(X|\pi, G,\Sigma)\neq \emptyset$, for some $n\ge 2$. Let $\varepsilon>0$ such that $\bigcap_{i=1}^n\overline{B_{\v}(x_i)}=\emptyset.$ Define $U_i=X\setminus \overline{B_{\v}(x_i)}$, then $\bigcup_{i=1}^nU_i=X.$ Note that $B_{\varepsilon}(x_i)\subset U_i^c$ for every $i=1,2,\dots,n.$ Then $(U_1,\dots,U_n)$ is an admissible open cover of $X$ with respect to $(x_i)_{i=1}^n,$ which implies $\text{mdim}_{\Sigma}((U_1,\dots,U_n)|\pi)>0$. Hence $$\text{mdim}_{\Sigma}(X|\pi)\geq \text{mdim}_{\Sigma}((U_1,\dots,U_n)|\pi)>0.$$
		
		Now suppose that $\text{mdim}_{\Sigma}(X|\pi)>0$. Combining Propositions \ref{propo322}  and \ref{ppppro324}, we get the desired result.
	\end{proof}

	\begin{prop}\label{closed}
		Let $\pi:(X,G)\to (Y,G)$ be a factor map between $G$-systems, $n \geq 2$ and $\Sigma$ a sofic approximation sequence for $G$. 
		Then $\mathrm{D}_{n}^{\mathrm{md}}(X|\pi,G,\Sigma)\cup \Delta_n(X)$ is a closed set and $\overline{\mathrm{D}_{n}^{\mathrm{md}}(X|\pi,G,\Sigma)}\subset \mathrm{D}_{n}^{\mathrm{md}}(X|\pi,G,\Sigma)\cup \Delta_n(X).$
	\end{prop}
	\begin{proof}
		For any $(x_i)_{i=1}^n\in \overline{\mathrm{D}_{n}^{\mathrm{md}}(X|\pi, G,\Sigma)}\setminus \Delta_n(X),$ there exist $\{(x_i^{(m)})_{i=1}^n:m\in \N\}\subset \mathrm{D}_{n}^{\mathrm{md}}(X|\pi, G,\Sigma)$ such that $(x_i^{(m)})_{i=1}^n\to (x_i)_{i=1}^n$, when $m\to\infty$. Let $\alpha$ be an admissible open cover of $X$ with respect to $(x_i)_{i=1}^n$. By definition, for sufficiently large $m\in \N$, $\alpha$ is   an admissible open cover of $X$ with respect to $(x_i^{(m)})_{i=1}^n$, which implies $\text{mdim}_{\Sigma}(\alpha|\pi)>0$. Thus $(x_i)_{i=1}^n\in \mathrm{D}_{n}^{\mathrm{md}}(X|\pi, G,\Sigma)$. This implies that $\overline{\mathrm{D}_{n}^{\mathrm{md}}(X|\pi, G,\Sigma)}\subset \mathrm{D}_{n}^{\mathrm{md}}(X|\pi, G,\Sigma)\cup \Delta_n(X).$ Moreover, $\mathrm{D}_{n}^{\mathrm{md}}(X|\pi, G,\Sigma)\cup \Delta_n(X)$ is a closed set of $X^{n}.$
		%
		%
	\end{proof}

	\begin{lem}\label{lllllem30}
		Let $\pi:(X,G)\to (Y,G)$ be a factor map between $G$-systems and $\Sigma$ a sofic approximation sequence for $G$. Let $\psi:(Z, G)\to(Y,G)$  a factor map of $\pi$ and $\phi:(X, G)\to(Z,G)$ a factor map with $\pi=\psi\circ\phi$. Then for any finite open cover $\alpha$ of $Z$, $\mathrm{mdim}_{\Sigma}(\phi^{-1}(\alpha)|\pi)\leq \mathrm{mdim}_{\Sigma}(\alpha|\psi).$
	\end{lem}
	\begin{proof}
		For any finite open cover $\mathcal{U}$ of $Z^d$ and $\xi\in X^d,$ we have 
		$$\sum_{U\in \alpha} 1_{U}(\phi^d (\xi))=\sum_{V\in (\phi^d)^{-1}(\mathcal{U})}1_{V}(\xi).$$
		Note that if $\{(\psi^d)^{-1}(y)\}_{y\in Y^d}\vee\beta \succ\alpha^d$, then 
		\begin{align*}
			(\phi^d)^{-1}(\{(\psi^d)^{-1}(y)\}_{y\in Y^d} \bigvee\beta)&=\{((\psi\circ\phi)^d)^{-1}(y)\}_{y\in Y^d} \bigvee (\phi^d)^{-1}(\beta)\\
			&=\{(\pi^d)^{-1}(y)\}_{y\in Y^d}\bigvee (\phi^d)^{-1}(\beta)\\
			&\succ (\phi^d)^{-1}(\alpha^d).
		\end{align*}
		Let $F \in \mathcal{F}(G), \delta>0$ and $\sigma:G \rightarrow \Sym(d)$ for some $d \in \N$. By Lemma \ref{li-2}, there exists $\delta'>0$ such that 
		\begin{align*}
			&\min_{\{(\pi^d)^{-1}(y)\}_{y\in Y^d}\vee\beta\succ (\phi^{-1}(\alpha))^d}\left(\max_{\xi\in \text{Map}(\rho_X,F,\delta',\sigma)}\sum_{U\in \beta}1_U(\xi)\right)\\
			&\leq \min_{\{(\psi^d)^{-1}(y)\}_{y\in Y^d}\vee\beta \succ\alpha^d}\left(\max_{\xi\in \text{Map}(\rho_X,F,\delta',\sigma)}\sum_{U\in \beta}1_U(\phi\circ\xi)\right)\\
			&\leq \min_{\{(\psi^d)^{-1}(y)\}_{y\in Y^d}\vee\beta \succ\alpha^d}\left(\max_{\xi'\in \text{Map}(\rho_Z,F,\delta,\sigma)}\sum_{U\in \beta}1_U(\xi')\right).
		\end{align*}
		Then $$D((\phi^{-1}(\alpha))^d|_{\text{Map}(\rho_X,F,\delta',\sigma)}|\pi^d)\leq D(\alpha^d|_{\text{Map}(\rho_Z,F,\delta,\sigma)}|\psi^d),$$
		which implies that $\text{mdim}_{\Sigma}(\phi^{-1}(\alpha)|\pi)\leq \text{mdim}_{\Sigma}(\alpha|\psi).$
	\end{proof}

	\begin{prop}\label{lift}
		Let $\pi:(X,G)\to (Y,G)$ be a factor map between $G$-systems, $n \geq 2$ and $\Sigma$ a sofic approximation sequence for $G$. Let $\psi:(Z, G)\to(Y,G)$ be a factor map of $\pi$ and $\phi:(X, G)\to(Z,G)$ a factor map with $\pi=\psi\circ\phi$. If $(x_i)_{i=1}^{n}\in D_{n}^{\mathrm{md}}(X|\pi,G,\Sigma)$ with $(\phi(x_i))_{i=1}^{n}\not\in\Delta_n(Z)$, then $(\phi(x_i))_{i=1}^{n}\in \mathrm{D}_{n}^{\mathrm{md}}(Z|\psi,G,\Sigma).$
	\end{prop}
	\begin{proof}
		Let $(x_i)_{i=1}^{n}\in \mathrm{D}_{n}^{\mathrm{md}}(X|\pi,G,\Sigma)$ with $(\phi(x_i))_{i=1}^{n}\not\in\Delta_n(Z)$. Suppose that $\alpha=(U_1,\dots,U_n)$ is an admissible open cover of $Z$ with respect to $(\phi(x_i))_{i=1}^{n},$ where $U_i^c$ is a closed neighborhood of $\phi(x_i)$ such that $U_i^c\cap U_j^c=\emptyset$ when $\phi(x_i)\neq \phi(x_j)$, $1\leq i,j\leq n$. 
		By the continuity of $\phi,$ we have $(\phi^{-1}(U_i))^c=\phi^{-1}(U_i^c)$ is a closed neighborhood of $x_i$ for every $i=1,2\dots,n.$ Moreover,
		$$(\phi^{-1}(U_i))^c\cap (\phi^{-1}(U_j))^c=\phi^{-1}(U_i^c)\cap \phi^{-1}(U_j^c)=\phi^{-1}(U_i^c\cap U_j^c)=\emptyset,$$
		where $x_i\neq x_j,1\leq i,j\leq n.$ Note that $\phi^{-1}(\alpha)$ is a finite open cover of $X.$ Thus $\phi^{-1}(\alpha)=(\phi^{-1}(U_1),\phi^{-1}(U_2),\dots,\phi^{-1}(U_n))$ is an admissible open cover of $X$ with respect to $(x_i)_{i=1}^n$. By Lemma \ref{lllllem30}, we have
		$$0<\text{mdim}_{\Sigma}(\phi^{-1}(\alpha)|\pi)\leq \text{mdim}_{\Sigma}(\alpha|\psi).$$
		Combining with Proposition \ref{admis},  we get  $(\phi(x_1),\dots,\phi(x_n))\in \mathrm{D}_{n}^{\mathrm{md}}(Z|\psi,G,\Sigma).$
	\end{proof}
	
	\begin{prop}
		Let $\pi_1:(X,G)\to (Y,G)$ be a factor map between $G$-systems, $n \geq 2$ and $\Sigma$ a sofic approximation sequence for $G$. Let $Z\subset X$ be a closed $G$-invariant subset. If $(x_i)_{i=1}^{n}\in \mathrm{D}_{n}^{\mathrm{md}}(Z|\pi_2, G,\Sigma),$  where $\pi_2=\pi_1|_Z:(Z,G)\to (\pi_1(Z),G)$, then $(x_i)_{i=1}^{n}\in \mathrm{D}_{n}^{\mathrm{md}}(X|\pi, G,\Sigma)$. 
	\end{prop}
	\begin{proof}
		Let $\alpha$ be an admissible open cover of $X$ with respect to $(x_i)_{i=1}^{n}\in Z^n\setminus \Delta_n(Z)$. Clearly, $\alpha|_Z$, its restriction to $Z$, is an admissible open cover of $Z$ with respect to $(x_i)_{i=1}^{n}$. According to the proof of Proposition \ref{subset}, we have $\text{mdim}_{\Sigma}(\alpha|\pi_1)\geq \text{mdim}_{\Sigma}(\alpha|_Z|\pi_2)$. Combining with $\text{mdim}_{\Sigma}(\alpha|_Z|\pi_2)>0$, we get the desired result.
	\end{proof}

	\subsection{Completely/Uniform positive sofic conditional mean dimension}\label{sbsec2}

	\begin{defn}
		Let $\pi:(X,G)\to (Y,G)$ be a nontrivial factor map between $G$-systems  and $\Sigma$ a sofic approximation sequence for $G$. We say $\pi$ has \emph{completely positive sofic conditional mean dimension} (sofic $X|\pi$-CPCMD) if for any proper factor $\psi:(Z,G)\to(Y,G)$ of $\pi$, one has $\mdim_{\Sigma}(Z|\psi)>0$. 
	\end{defn}
	
	The following result is an observation.
	\begin{prop}\label{univ1}
		Let $\pi:(X,G)\to (Y,G)$ be a factor map between $G$-systems with \(|X| \geq 2\) and \(\mathrm{mdim}_{\Sigma} (X|\pi) \geq 0\). Let $\Sigma$ a sofic approximation sequence for $G$ and $S_\pi^{\Sigma}$ the relative zero sofic conditional  mean dimension relation. Then $\pi$ has sofic $X|\pi$-CPCMD if and only if $S_\pi^{\Sigma}=R_\pi$.
	\end{prop}

	\begin{defn}
		Let $\pi \colon (X, G) \to (Y, G)$ be a factor map between $G$-systems. A finite open cover $\alpha$ of $X$ is called \textit{non-dense-on-fibre}, if there is $y \in Y$ such that $\pi^{-1}(y)$ is not contained in any element of $\overline{\alpha}$ (which consists of the closures of elements of $\alpha$ in $X$). Clearly, if a finite open cover $\alpha = \{U_1, U_2\}$ is non-dense-on-fibre, then $\pi(U_1) \cap \pi(U_2) \neq \emptyset$.
	\end{defn}
	\begin{defn}
		Let $\pi \colon (X, G) \to (Y, G)$ be a nontrivial factor map between $G$-systems and $\Sigma$ a sofic approximation sequence for $G$.
		\begin{enumerate}
			\item Let $n \geq 2$. Say $\pi$ has \textit{sofic $X|\pi$-UPCMD of order $n$}, if any non-dense-on-fibre open cover $\alpha$ of $X$ by $n$-sets has positive sofic conditional mean dimension, i.e., $\text{mdim}_{\Sigma}(\alpha|\pi)>0$. (When $n = 2$, we say simply $\pi$ or $(X, G)$ has \textit{sofic $X|\pi$-UPCMD})
			\item Say $(X, G)$ or $\pi$ has \textit{sofic $X|\pi$-UPCMD of all orders}, if any non-dense-on-fibre open cover $\alpha$ of $X$ by finitely many sets has positive sofic conditional mean dimension, i.e., it has sofic $X|\pi$-UPCMD of order $n$ for any $n \geq 2$.
		\end{enumerate}
	\end{defn}
	\begin{rem}\label{upmd-rmk}
		\
		\begin{enumerate}
			\item[(1)]By Proposition~\ref{subset Rpi},
			it is easy to show that, for all $n \geq 2$, $\pi$ has sofic $X|\pi$-UPCMD of order $n$ iff	
			$\mathrm{D}_{n}^{\mathrm{md}}(X|\pi, G,\Sigma) = R^{n}_\pi \setminus \Delta_n(X)$.

			\item[(2)] Properties of sofic $X|\pi$-CPCMD and sofic $X|\pi$-UPCMD are both stable under factor maps. Precisely, let $\pi \colon (X, G) \to (Y, G)$ be a given factor map between $G$-systems, if $\pi$ has property P (P denotes one of the properties: sofic $X|\pi$-UPCMD and sofic $X|\pi$-CPCMD) and $(Z, G)$ is a proper factor of $(X, G)$ with respect to $(Y, G)$; then $(Z, G)$ also has property P with respect to $(Y, G)$. Note that sofic $X|\pi$-UPCMD (resp. sofic $X|\pi$-CPCMD) recovers sofic UPMD (resp. sofic CPMD) introduced in \cite{GRFYG}, when $(Y,G)$ is trivial.
		\end{enumerate}
	\end{rem}
	
	\begin{thm}\label{univ}
		Let $\pi:(X,G)\to (Y,G)$ be a factor map between $G$-systems with \(\mathrm{mdim}_{\Sigma} (X|\pi) \geq 0\), $\Sigma$ a sofic approximation sequence for $G$ and $S_\pi^{\Sigma}$ the relative zero sofic conditional mean dimension relation. The following holds:
		\begin{enumerate}
			\item The smallest closed \(G\)-invariant equivalence relation that contains $\mathrm{D}_{2}^{\mathrm{md}}(X|\pi, G,\Sigma)$ is contained in $S_\pi^{\Sigma}$.
			\item $\pi$ has sofic $X|\pi$-CPCMD if the smallest closed \(G\)-invariant equivalence relation that contains $\mathrm{D}_{2}^{\mathrm{md}}(X|\pi, G,\Sigma)$ is $R_\pi$.
		\end{enumerate}
	\end{thm}
	
	\begin{proof}
		We begin by proving statement (1). Let $Q$ be the smallest closed $G$-invariant equivalence relation containing $\mathrm{D}_{2}^{\mathrm{md}}(X|\pi, G, \Sigma)$, and let $\phi: (X, G) \to (X^{\Sigma}_\pi, G)$ be the factor map where $X^{\Sigma}_\pi = X / S_\pi^{\Sigma}$. To establish (1), it suffices to show that $\mathrm{D}_{2}^{\mathrm{md}}(X|\pi, G, \Sigma) \subset S_\pi^{\Sigma}$, which would imply $Q \subset S_\pi^{\Sigma}$.
		
		Suppose for contradiction that there exists $(x, y) \in \mathrm{D}_{2}^{\mathrm{md}}(X|\pi, G, \Sigma)$ with $(x, y) \notin S_\pi^{\Sigma}$. By the definition of $S_\pi^{\Sigma}$, there exist factor maps $\theta_1: (X, G) \to (Z, G)$ and $\theta_2: (Z, G) \to (Y, G)$ such that:
		\begin{itemize}
			\item $\pi = \theta_2 \circ \theta_1$,
			\item $\mathrm{mdim}_{\Sigma}(Z|\theta_2) \leq 0$,
			\item $\theta_1(x) \neq \theta_1(y)$.
		\end{itemize}
		
		By Proposition~\ref{lift}, $(\theta_1(x), \theta_1(y)) \in \mathrm{D}_{2}^{\mathrm{md}}(Z|\theta_2, G, \Sigma)$. However, Theorem~\ref{main3} implies $\mathrm{mdim}_{\Sigma}(Z|\theta_2) > 0$, leading to a contradiction. Thus, $Q \subset S_\pi^{\Sigma}$, proving (1).
		
		Statement (2) follows directly from (1). If $Q = R_\pi$, then (1) implies $S_\pi^{\Sigma} = R_\pi$. By the definition of $S_\pi^{\Sigma}$, this means $\pi$ admits no non-trivial intermediate factors with non-positive sofic conditional mean dimension.
	\end{proof}
	
	\begin{cor}
		Let $\pi \colon (X, G) \to (Y, G)$ be a nontrivial factor map between $G$-systems with \(|X| \geq 2\) and $\Sigma$ a sofic approximation sequence for $G$. If $\pi$ has sofic $X|\pi$-UPCMD, then $\pi$ has sofic $X|\pi$-CPCMD.
	\end{cor}
	\begin{proof}
		Assume the factor map $\pi \colon (X, G) \to (Y, G)$ has sofic $X|\pi$-UPCMD.
		By Remark~\ref{upmd-rmk}, this implies the relation 
		$ \mathrm{D}_{2}^{\mathrm{md}}(X|\pi, G, \Sigma) = R_\pi \setminus \Delta_2(X).$
		
		Since $|X| \geq 2$, $\pi$ admits non-trivial non-dense-on-fibre open covers. The sofic $X|\pi$-UPCMD assumption consequently yields 
		$\mathrm{mdim}_{\Sigma}(X|\pi) > 0$.
		
		This allows us to apply Theorem~\ref{univ}-(1), which establishes the inclusion
		$ R_\pi \setminus \Delta_2(X) \subset S_{\pi}^{\Sigma}.$
		This inclusion immediately forces the equality $S_{\pi}^{\Sigma} = R_\pi$.
		Then the conclusion now follows by invoking Proposition~\ref{univ1}, which completes the proof.
	\end{proof}

	\section{Sofic mean dimension of $(Y,G)$ relative to $(X,G)$ }\label{sec5}
	Similar to Section~\ref{sec3}, in this section,  we further study the sofic mean dimension of $(Y,G)$ relative to $(X,G)$ defined by Li and Liang in \cite{LB}. 
	
	Let $\pi:(X,G)\to (Y,G)$ be a factor map between $G$-systems and $\rho$ be a continuous pseudometric on $X$. $\rho$ is called \emph{dynamically generating} if $\sup_{g\in G}\rho(x,y)>0$ for any distinct $x,y\in X$. Let $F \in \mathcal{F}(G)$ and $\delta>0$. Let $\sigma$ be a map from $G$ to $\Sym(d)$ for some $d \in \N$. Consider the closed subset $\pi^{d}(\Map(\rho,F,\delta,\sigma))$ of $Y^d$. For a finite open cover $\alpha$ of $X$, denote $D(\alpha^d|_{\pi^d(\Map(\rho,F,\delta,\sigma))})$ by $D(\pi,\alpha,\rho,F,\delta,\sigma)$. 
	\begin{defn}
		Let $\pi:(X,G)\to (Y,G)$ be a factor map between $G$-systems and $\rho$ a dynamically generating continuous pseudometric on $X$. Let $F \in \mathcal{F}(G)$ and $\delta > 0$. For a finite open cover $\alpha$ of $Y$ we define
		\[
		\mdim_{\Sigma}(\pi, \alpha, \rho, F, \delta) = \lim_{i \to \infty} \frac{D(\pi, \alpha, \rho, F, \delta, \sigma_i)}{d_i},
		\]
		and
		\[
		\mdim_{\Sigma}(\pi, \alpha, \rho) = \inf_{F \in \mathcal{F}(\Gamma)} \inf_{\delta > 0} \mdim_{\Sigma}(\pi, \alpha, \rho, F, \delta).
		\]
		If $\Map(\rho, F, \delta, \sigma_i) = \emptyset$ for all sufficiently large $i \in \mathbb{N}$, we set $\mdim_{\Sigma}(\pi, \alpha, \rho, F, \delta) = -\infty$.
		
		We define the \emph{sofic mean topological dimension of $(Y,G)$ relative to the extension $(X,G)$} as
		\[
		\mdim_{\Sigma}(Y|X, \rho,G) := \sup_{\alpha} \mdim_{\Sigma}(\pi, \alpha, \rho)
		\]
		for $\alpha$ ranging over all finite open covers of $Y$. The quantities $\mdim_{\Sigma}(\pi, \alpha, \rho)$ and $\mdim_{\Sigma}(Y|X, \rho)$ do not depend on the choice of $\rho$ (see for example \cite[Lemma 8.2]{LB}), and we shall write them as $\mdim_{\Sigma}(\pi,\alpha)$ and $\mdim_{\Sigma}(Y|X)$ respectively.  When $\pi$ is a homeomorphism,  $\mathrm{mdim}_{\Sigma}(Y|X)$ recovers as the \emph{sofic mean dimension} $\mathrm{mdim}_{\Sigma}(Y)$ of $(Y,G)$.
	\end{defn}
	
	Meanwhile, Liang \cite{LBB2} gave an equivalent definition of sofic mean topological dimension of $(Y,G)$ relative to the extension $(X,G)$.
	\begin{defn}
		Let $\pi:(X,G)\to (Y,G)$ be a factor map between $G$-systems and $\Sigma$ a sofic approximation sequence for $G$. For every $\varepsilon>0, F \in \mathcal{F}(G)$, we set
		$$
		\operatorname{mdim}_{\Sigma,\varepsilon}\left(\rho_Y | \rho_X, F, \delta\right):=\limsup _{i \rightarrow \infty} \frac{1}{d_i} \operatorname{Wdim}_{\varepsilon}\left(\pi^{d_i}\left(\operatorname{Map}\left(\rho_X, F, \delta, \sigma_i\right)\right), \rho_{Y, \infty}\right)
		$$
		If $\operatorname{Map}\left(\rho, F, \delta, \sigma_i\right)$ is empty for all sufficiently large $i$, we set $\operatorname{mdim}_{\Sigma,\varepsilon}\left(\rho_Y | \rho_X, F, \delta\right)=$ $-\infty$. Set
		$$
		\operatorname{mdim}_{\Sigma,\varepsilon}\left(\rho_Y | \rho_X\right):=\inf _{F \in \mathcal{F}(G)} \inf _{\delta>0} \operatorname{mdim}_{\Sigma,\varepsilon}\left(\rho_Y | \rho_X, F, \delta\right)
		$$
		and
		$$
		\operatorname{mdim}_{\Sigma}\left(\rho_Y | \rho_X\right):=\sup _{\varepsilon>0} \operatorname{mdim}_{\Sigma,\varepsilon}\left(\rho_Y | \rho_X\right).
		$$
		
		When $\rho_X$ and $\rho_Y$ are compatible metrics the \emph{sofic mean dimension $\operatorname{mdim}_{\Sigma}(Y |X)$ of $(Y,G)$ relative to $(X,G)$} is defined as the value $\operatorname{mdim}_{\Sigma}\left(\rho_Y|\rho_X\right)$. It's known that $\operatorname{mdim}_{\Sigma}\left(\rho_Y|\rho_X\right)$ is independent of the choice of compatible metrics $\rho_X$ and $\rho_Y$. 
	\end{defn}
	
	Now we will show the fundamental properties of relative sofic mean dimension.
	
	\begin{prop}
		Let $\pi:(X,G)\to (Y,G)$ be a factor map between $G$-systems, $\Sigma$ a sofic approximation sequence for $G$ and $\alpha$ a finite open cover of $Y$, then $$\mathrm{mdim}_{\Sigma}(\pi,\alpha,\rho,F,\delta)\leq D(\alpha).$$
		Moreover $\mathrm{mdim}_{\Sigma}(\pi,\alpha)\leq D(\alpha).$
	\end{prop}
	\begin{proof}
		Let $F\in \mathcal{F}(G)$, $\delta>0,d\in \N$ and $\sigma:G\to \text{Sym}(d)$.	Since $\pi^d(\text{Map}(\rho,F,\delta,\sigma) \subset Y^d$, then 
		\begin{align*}
			D(\alpha^d|_{\pi^d(\text{Map}(\rho,F,\delta,\sigma))})
			&=D(\pi,\alpha,\rho,F,\delta,\sigma)\\
			&\leq D(\alpha^d)\\
			&\leq dD(\alpha).\end{align*}
		This implies that $$\text{mdim}_{\Sigma}(\pi,\alpha,\rho,F,\delta)=\varlimsup_{i\to \infty}\frac{D(\pi,\alpha,\rho,F,\delta,\sigma_i)}{d_i}\leq D(\alpha).$$
		Moreover,$$\text{mdim}_{\Sigma}(\pi,\alpha)=\inf_{F'\in \mathcal{F}(G)}\inf_{\delta>0}\text{mdim}_{\Sigma}(\pi,\alpha,\rho,F',\delta)\leq D(\alpha).$$
		This proof is complete.
	\end{proof}
	
	\begin{prop}\label{jajajajaj}
		Let $\pi:(X,G)\to (Y,G)$ be a factor map between $G$-systems, $\Sigma$ a sofic approximation sequence for $G$ and $\alpha, \beta$ be finite open covers of $Y$. The following holds:
		\begin{itemize}
			\item [(1)]If $\alpha\succ \beta$, then $\mathrm{mdim}_{\Sigma}(\pi,\alpha)\geq \mathrm{mdim}_{\Sigma}(\pi,\beta).$
			\item[(2)] $\mathrm{mdim}_{\Sigma}(\pi,\alpha\bigvee \beta)\leq \mathrm{mdim}_{\Sigma}(\pi,\alpha)+\mathrm{mdim}_{\Sigma}(\pi,\beta).$
		\end{itemize}
		\begin{proof}
			Let $F\in \mathcal{F}(G)$, $\delta>0$, and $\sigma:G\to \text{Sym}(d)$ for some $d \in \N$.\\
			(1) Since $\alpha\succ \beta$, then $\alpha^d|_{\pi^d(\text{Map}(\rho,F,\delta,\sigma))}\succ \beta^d|_{\pi^d(\text{Map}(\rho,F,\delta,\sigma))}$. By Proposition \ref{jbwcjbc} (when $Y$ is a singleton), we have \begin{align*}
				D(\pi,\alpha,\rho,F,\delta,\sigma)&=D(\alpha^d|_{\pi^d(\text{Map}(\rho,F,\delta,\sigma))})\\
				&\geq D(\beta^d|_{\pi^d(\text{Map}(\rho,F,\delta,\sigma))})\\
				&=D(\pi,\beta,\rho,F,\delta,\sigma),\end{align*}
			which implies that $\text{mdim}_{\Sigma}(\pi,\alpha)\geq \text{mdim}_{\Sigma}(\pi,\beta).$\\
			(2)Since $(\alpha\vee \beta)^d=\alpha^d\vee \beta^d$, by Proposition \ref{subbb} (when $Y$ is a singleton), we get
			\begin{align*}
				D((\alpha\vee \beta)^d|_{\pi^d(\text{Map}(\rho,F,\delta,\sigma))})&=D((\alpha^d\vee \beta^d)|_{\pi^d(\text{Map}(\rho,F,\delta,\sigma))})\\
				&=D(\pi,\alpha\vee \beta,\rho,F,\delta,\sigma)\\
				&\leq D(\alpha^d|_{\pi^d(\text{Map}(\rho,F,\delta,\sigma))})+D(\beta^d|_{\pi^d(\text{Map}(\rho,F,\delta,\sigma))})\\
				&=D(\pi,\alpha,\rho,F,\delta,\sigma)+D(\pi,\beta,\rho,F,\delta,\sigma).
			\end{align*}
			Thus $\text{mdim}_{\Sigma}(\pi,\alpha\bigvee \beta)\leq \text{mdim}_{\Sigma}(\pi,\alpha)+\text{mdim}_{\Sigma}(\pi,\beta).$
		\end{proof}
	\end{prop}	
	
	Relative sofic mean dimension also maintains an intrinsic connection with factor maps, so we will also discuss its properties relative to different types of factor maps.
	
	\begin{prop}\label{isomorphism2}
		Let $\pi: (X,G)\to (Y,G)$ be a factor map between $G$-systems, $\Sigma$ a sofic approximation sequence for $G$ and $\psi:(Z, G)\to(Y,G)$  a factor map of $\pi$. If $\psi$ is an isomorphism, then $\mathrm{mdim}_{\Sigma} (Y|X)= \mathrm{mdim}_{\Sigma} (Z|X)$.
	\end{prop}
	\begin{proof}
		
		Let $F\in \mathcal{F}(G)$, $\delta>0$, and $\sigma:G\to \text{Sym}(d)$ for some $d \in \N$. For any $\ep>0$, there exists $\ep'>0$ such that for any $z,z'\in Z^d$ with $\rho_{Z,\infty}(z,z')<\ep'$, then one has $\rho_{Y,\infty}(\psi^d(z),\psi^d(z'))<\ep$. Let  $\phi:(X, G)\to(Z,G)$ be a factor map with $\pi=\psi\circ\phi$. Then for any $\tau,\tau'\in \phi^d(\Map(\rho_X,F,\d,\sigma))$ with $\rho_{Z,\infty}(\tau,\tau')<\ep'$, we have $\rho_{Y,\infty}(\psi^d (\tau) ,\psi^d (\tau'))=\rho_{Y,\infty}(\psi \circ \tau,\psi\circ \tau')<\ep$.  Note that $\psi^d:\phi^d(\Map(\rho_X,F,\d,\sigma))\to\pi^d(\Map(\rho_X,F,\d,\sigma))$ sending $\phi\circ \xi$ to $\pi\circ \xi$, for any $\xi \in \Map(\rho_X,F,\d,\sigma)$. 
		
		First we show that $$\mathrm{Wdim}_{\varepsilon}(\pi^d(\mathrm{Map}(\rho_X, F, \delta, \sigma), \rho_{Y,\infty}))\leq\mathrm{Wdim}_{\varepsilon'}(\phi^d(\mathrm{Map}(\rho_X, F, \delta, \sigma), \rho_{Z,\infty})),$$
		for any nonempty finite subset $F$ of $G$, $\d>0$ and any map $\sigma$ from $G$ to $\Sym(d)$ for some $d\in\N$.
		
		By the definition, there exists an $(\ep',\rho_{Z,\infty})$-embedding $f:\phi^d(\mathrm{Map}(\rho_X, F, \delta, \sigma))\ra P$ such that for any $\xi,\xi'\in\mathrm{Map}(\rho_X, F, \delta, \sigma)$ satisfies $f(\phi\circ \xi)=f(\phi\circ \xi')$ , we have $\rho_{Z,\infty}(\phi\circ \xi,\phi\circ \xi')<\ep'$. Assume that $\mathrm{Wdim}_{\varepsilon'}(\phi^d\mathrm{Map}(\rho_X, F, \delta, \sigma), \rho_{Z,\infty})=\dim(P).$
		
		Since $\psi$ is an isomorphism and $\pi=\psi\circ\phi$, one can well define a continuous map $h:\pi^d(\Map(\rho_X,F,\d,\sigma))\to\phi^d(\Map(\rho_X,F,\d,\sigma)))$ sending $\pi\circ \xi$ to $\phi\circ \xi$, for any $\xi\in \Map(\rho_X,F,\d,\sigma)$.
		We claim that $f\circ h:\pi^d(\mathrm{Map}(\rho_X, F, \delta, \sigma))\ra P $ is an $(\ep,\rho_{Y,\infty})$-embedding.
		
		Indeed, for any  $\xi,\xi'\in\mathrm{Map}(\rho_X, F, \delta, \sigma)$ satisfying $f\circ h(\pi\circ \xi)=f\circ h(\pi\circ \xi')$  i.e. $f (\phi\circ \xi)=f(\phi\circ \xi')$, which implies  that $\rho_{Z,\infty}(\phi\circ \xi,\phi\circ \xi')<\ep'$. Then  $\rho_{Y,\infty}(\psi\circ\phi\circ \xi,\psi\circ\phi\circ \xi')=\rho_{Y,\infty}(\pi\circ \xi,\pi\circ \xi')<\ep$ which implies that $f\circ h$ is an $(\ep,\rho_{Y,\infty})$-embedding.
		Hence,
		\begin{align*}
			\mathrm{Wdim}_{\varepsilon}(\pi^d(\mathrm{Map}(\rho_X, F, \delta, \sigma)), \rho_{Y,\infty})&\leq\dim(P)\\
			&=\mathrm{Wdim}_{\varepsilon'}(\phi^d(\mathrm{Map}(\rho_X, F, \delta, \sigma)), \rho_{Z,\infty}).
		\end{align*}
		Then $\mathrm{mdim}_{\Sigma} (Y|X)\le \mathrm{mdim}_{\Sigma} (Z|X)$.

		By the continuity of $\psi^{-1}$, there exists $\varepsilon''>0$ such that for any $\pi\circ \xi,\pi\circ \xi'\in\pi^d(\mathrm{Map}(\rho_X, F, \delta, \sigma))$ 
		with $\rho_{Y,\infty}(\pi\circ \xi,\pi\circ \xi')<\varepsilon''$, one has $\rho_{Z,\infty}(\phi\circ \xi,\phi\circ \xi')<\ep$. Next we will show that 
		$$\mathrm{Wdim}_{\varepsilon}(\phi^d(\mathrm{Map}(\rho_X, F, \delta, \sigma)), \rho_{Z,\infty})\leq \mathrm{Wdim}_{\varepsilon''}(\pi^d(\mathrm{Map}(\rho_X, F, \delta, \sigma)), \rho_{Y,\infty}).$$
		
		By the definition, there exists an $(\ep'',\rho_{Y,\infty})$-embedding $g:\pi^d(\mathrm{Map}(\rho_X, F, \delta, \sigma_i))\ra Q$ such that for any $\xi,\xi'\in\mathrm{Map}(\rho_X, F, \delta, \sigma)$ satisfying $g(\pi\circ \xi)=g(\pi\circ \xi')$, we have $\rho_{Y,\infty}(\pi\circ \xi,\pi\circ \xi')<\ep''$. Assume that $\mathrm{Wdim}_{\varepsilon''}(\pi^d(\mathrm{Map}(\rho_X, F, \delta, \sigma)), \rho_{Y,\infty})=\dim(Q)$. Then we show that $g\circ\psi^d:\phi^d(\mathrm{Map}(\rho_X, F, \delta, \sigma))\ra Q$ is an $(\ep,\rho_{Z,\infty})$-embedding.
		
		In fact, for any  $\xi,\xi'\in\mathrm{Map}(\rho_X, F, \delta, \sigma)$ satisfying $g\circ \psi^d(\phi\circ \xi)=g\circ \psi^d(\phi\circ \xi')$  i.e. $g (\pi\circ \xi)=g(\pi\circ \xi')$, which implies that $\rho_{Y,\infty}(\pi\circ \xi,\pi\circ \xi')=\rho_{Y,\infty}(\psi\circ\phi\circ \xi,\psi\circ\phi\circ \xi')<\ep''$. Then $\rho_{Z,\infty}(\phi\circ \xi,\phi\circ \xi')<\ep$.
		Thus, $g\circ \psi^d$ is an $(\ep,\rho_{Z,\infty})$-embedding. Hence,
		\begin{align*}
			\mathrm{Wdim}_{\varepsilon}(\phi^d(\mathrm{Map}(\rho_X, F, \delta, \sigma)), \rho_{Z,\infty}) &\leq \dim(Q)\\
			&=\mathrm{Wdim}_{\varepsilon''}(\pi^d(\mathrm{Map}(\rho_X, F, \delta, \sigma)), \rho_{Y,\infty}).
		\end{align*}
		Then $\mathrm{mdim}_{\Sigma} (Z|X)\leq \mathrm{mdim}_{\Sigma} (Y|X)$.
	\end{proof}

	\begin{prop}\label{isomorp2w}
		Let $\pi: (X,G)\to (Y,G)$ be a factor map between $G$-systems, $\Sigma$ a sofic approximation sequence for $G$ and $\psi:(Z, G)\to(Y,G)$  a factor map of $\pi$. Let  $\phi:(X, G)\to(Z,G)$ be a factor map with $\pi=\psi\circ\phi$. Then $\mathrm{mdim}_{\Sigma} (Y|X)\leq  \mathrm{mdim}_{\Sigma} (Y|Z)$. If $\phi$ is an isomorphism, then $\mathrm{mdim}_{\Sigma} (Y|X)= \mathrm{mdim}_{\Sigma} (Y|Z)$.
	\end{prop}
	\begin{proof}
		Let $F \in \mathcal{F}(G), \delta>0, \varepsilon>0$, by Lemma \ref{li-2}, there exists $\delta'>0$, for every map $\sigma$ from $G$ to $\Sym(d)$ for some $d \in \N$, one has
		$$\phi^d(\Map(\rho_X,F,\delta',\sigma) \subset \Map(\rho_Z,F,\delta,\sigma),$$
		where $\phi^d:X^d \ra Y^d$ sends $\xi$ to $\phi \circ \xi$. Thus
		$$\pi^d(\Map(\rho_X,F,\delta',\sigma)= \psi^d \circ\phi^d(\Map(\rho_X,F,\delta',\sigma)\subset \psi^d(\Map(\rho_Z,F,\delta,\sigma)).
		$$
		This implies that
		$$\mathrm{Wdim}_{\varepsilon}(\pi^d(\mathrm{Map}(\rho_X, F, \delta', \sigma), \rho_{Y,\infty}))\leq\mathrm{Wdim}_{\varepsilon}(\psi^d(\mathrm{Map}(\rho_Z, F, \delta, \sigma), \rho_{Y,\infty})).$$
		Hence
		$$\operatorname{mdim}_{\Sigma,\varepsilon}\left(\rho_Y | \rho_X\right) \leq\operatorname{mdim}_{\Sigma,\varepsilon}\left(\rho_Y | \rho_Z\right),$$
		which implies that $\mathrm{mdim}_{\Sigma} (Y|X)\leq  \mathrm{mdim}_{\Sigma} (Y|Z).$
		
		On the other hand, as $\phi$ is an isomorphism, by Lemma \ref{li-2}, there exists $\delta''>0$, such that
		$$(\phi^d)^{-1}(\Map(\rho_Z,F,\delta'',\sigma) \subset \Map(\rho_X,F,\delta,\sigma),$$
		then
		$$\Map(\rho_Z,F,\delta'',\sigma) \subset \phi^d(\Map(\rho_X,F,\delta,\sigma)).$$
		Thus
		$$\psi^d(\Map(\rho_Z,F,\delta'',\sigma)) \subset \psi^d \circ \phi^d(\Map(\rho_X,F,\delta,\sigma))=\pi^d(\Map(\rho_X,F,\delta,\sigma)).$$
		This implies that
		$$\mathrm{Wdim}_{\varepsilon}(\psi^d(\mathrm{Map}(\rho_Z, F, \delta'', \sigma), \rho_{Y,\infty}))\leq\mathrm{Wdim}_{\varepsilon}(\pi^d(\mathrm{Map}(\rho_X, F, \delta, \sigma), \rho_{Y,\infty})).$$
		Hence
		$$\operatorname{mdim}_{\Sigma,\varepsilon}\left(\rho_Y | \rho_Z\right) \leq\operatorname{mdim}_{\Sigma,\varepsilon}\left(\rho_Y | \rho_X\right),$$
		which implies that $\mathrm{mdim}_{\Sigma} (Y|Z)\leq  \mathrm{mdim}_{\Sigma} (Y|X).$
	\end{proof}
	\begin{prop}\label{subset2}
		Let $\pi_1: (X,G)\to (Y,G)$ be a factor map between $G$-systems, $\Sigma$ a sofic approximation sequence for $G$ and $Z$ be a closed \(G\)-invariant subset of $Y$, $X':=\pi_1^{-1}(Z)$ and $\pi_2:=\pi_1|_{X'}:(X',G)\to(Z,G)$. Then $\mathrm{mdim}_{\Sigma} (Z|X') \leq \mathrm{mdim}_{\Sigma} (Y|X)$.
	\end{prop}
	\begin{proof}

		Let $\rho_X$ and $\rho_Y$ be compatible metrics on $X$ and $Y$ respectively. Then $\rho_X$ restricts to a compatible metric $\rho_{X'}$ on $X'$ and $\rho_Y$ restricts to a compatible metric $\rho_Z$ on $Z$. Let \(\alpha\) be a finite open cover of \(Z\). Then there exists a finite cover $\beta$ of $Y$ such that $\alpha$ is exactly the covering obtained by intersecting $\beta$ with $Z$. Note that $\mathrm{Map}(\rho_{X'}, F, \delta, \sigma)\subset\mathrm{Map}(\rho_X, F, \delta, \sigma)$, for any nonempty finite subset $F$ of $G$, $\d>0$ and any map $\sigma$ from $G$ to $\Sym(d)$ for some $d\in\N$. Then $\pi_2^d(\mathrm{Map}(\rho_{X'}, F, \delta, \sigma))\subset\pi_1^d(\mathrm{Map}(\rho_X, F, \delta, \sigma))$, for any nonempty finite subset $F$ of $G$, $\d>0$ and any map $\sigma$ from $G$ to $\Sym(d)$ for some $d\in\N$. Furthermore, the restriction of $\beta^{d}|_{\pi_1^d(\mathrm{Map}(\rho_{X}, F, \delta, \sigma))}$ on $\pi_2^d(\mathrm{Map}(\rho_{X'}, F, \delta, \sigma))$ is exactly $\alpha^{d}|_{\pi_2^d(\mathrm{Map}(\rho_{X'}, F, \delta, \sigma))}$. Thus
		$D(\pi_2,\alpha,\rho_{X'},F,\d,\sigma)\leq D(\pi_1,\beta,\rho_{X},F,\d,\sigma)$.
		It follows that,
		$$\mdim_{\Sigma}(\pi_2,\alpha,\rho_{X'})\leq\mdim_{\Sigma}(\pi_1,\beta,\rho_{X})\leq\mdim_{\Sigma}(Y|X).$$
		Since $\alpha$ is an arbitrary finite open cover of $Z$, one has $$\mdim_{\Sigma}(Z|X')=\sup_{\alpha}\mdim_{\Sigma}(\pi_2,\alpha,\rho_{X'})\leq\mdim_{\Sigma}(Y|X).$$
		This proof is finished.
	\end{proof}
	
	\begin{prop}
		Let $\pi_n: (X_n,G)\to (Y_n,G)$ be a factor map between $G$-systems for each \(1 \leq n < R\), where \(R \in \mathbb{N} \cup \{\infty\}\) and $\Sigma$ a sofic approximation sequence for $G$. Consider the factor map $\pi:=\prod_{1 \leq n < R}\pi_n$ between the product systems $(X := \prod_{1 \leq n < R}X_n,G)$ and $(Y := \prod_{1 \leq n < R}Y_n,G)$. Then
		\[
		\mathrm{mdim}_{\Sigma} (Y|X) \leq \sum_{1 \leq n < R} \mathrm{mdim}_{\Sigma} (Y_n|X_n).
		\]
	\end{prop}
	\begin{proof}
		For any $n\in \N$, let $\rho^{(n)}_X$ be compatible metrics on $X_n$ and denote by $p_n$ the projection of $X$ onto $X_n$ and $p'_n$ the projection of $Y$ onto $Y_n$.
		Let $\rho_X$ be compatible metrics on $X$  and $\alpha$  a finite open cover of $Y$. Then there exists an $N \in \mathbb{N}$ with $N < R$ and a finite open covers $\beta_n$ of $Y_n$ for all $1 \leq n \leq N$ such that
		\[
		\beta := \bigvee_{n=1}^N {(p_n')}^{-1}(\beta_n) \succ \alpha.
		\]
		
		\begin{claim}\label{jcwjkcqqqqq}
			For any nonempty finite subset $F$ of $G$ and $\d>0$, we can find $\d'>0$ such that for any map $\sigma$ from $G$ to $\Sym(d)$ for some $d\in\N$, we have
			$${(p_n')}^d \circ\pi^d(\mathrm{Map}(\rho_X, F, \delta', \sigma))\subset\pi_n^d (\mathrm{Map}(\rho_X^{(n)}, F, \delta, \sigma)).$$
		\end{claim} 
		\begin{proof}[Proof of claim \ref{jcwjkcqqqqq}]
			Let $F$ be a nonempty finite subset of $G$ and $\d>0$. We first show that
			$${(p_n')}^d\circ\pi^d(\mathrm{Map}(\rho_X, F, \delta, \sigma))=\pi_n^d\circ p_n^d (\mathrm{Map}(\rho_X, F, \delta, \sigma)),$$
			for any nonempty finite subset $F$ of $G$, $\d>0$ and any map $\sigma$ from $G$ to $\Sym(d)$ for some $d\in\N$.
			
			For any $\xi\in \mathrm{Map}(\rho_X, F, \delta, \sigma)$, we know that $p'_n\circ\pi\circ\xi\in {(p_n')}^d\circ\pi^d(\mathrm{Map}(\rho_X, F, \delta, \sigma))$, which implies $\pi_n\circ p_n\circ\xi=p'_n\circ \pi\circ\xi\in {(p_n')}^d\circ\pi^d(\mathrm{Map}(\rho_X, F, \delta, \sigma))$. That is $\pi_n^d\circ p_n^d (\mathrm{Map}(\rho_X, F, \delta, \sigma))\subset {(p_n')}^d\circ\pi^d(\mathrm{Map}(\rho_X, F, \delta, \sigma))$. By symmetry, we have $\pi_n^d\circ p_n^d (\mathrm{Map}(\rho_X, F, \delta, \sigma))\supset {(p_n')}^d\circ\pi^d(\mathrm{Map}(\rho_X, F, \delta, \sigma))$. This implies the desired result.

			By Lemma~\ref{li-2},  there exists $\d'>0$ such that
			$ p_n^d (\mathrm{Map}(\rho_X, F, \delta, \sigma))\subset \mathrm{Map}(\rho_X^{(n)}, F, \delta, \sigma)$, which implies $$\pi_n^d\circ p_n^d (\mathrm{Map}(\rho_X, F, \delta', \sigma))\subset \pi_n^d (\mathrm{Map}(\rho_X^{(n)}, F, \delta, \sigma)).$$

			This proves the claim.
		\end{proof}	
		By claim \ref{jcwjkcqqqqq}, we can find a continuous map
		\[
		\Phi_n : \pi^d(\Map(\rho, F, \delta', \sigma)) \to \pi_n^d(\Map(\rho^{(n)}, F, \delta, \sigma))
		\]
		sending $\pi\circ\xi$ to $p'_n \circ\pi \circ \xi$ for each $1 \leq n \leq N$.
		Note that
		\[
		\beta^d|_{\pi^d(\Map(\rho,F,\delta',\sigma))} = \bigvee_{n=1}^N \Phi_n^{-1} \left( \beta_n^d|_{\pi_n^d(\Map(\rho^{(n)},F,\delta,\sigma))} \right).
		\]
		Thus
		\begin{align*}
			D(\pi,\alpha, \rho_X, F, \delta', \sigma) \leq D(\pi,\beta, \rho_X, F, \delta', \sigma)& \leq \sum_{n=1}^N D( \Phi_n^{-1} ( \beta_n^d|_{\pi_n^d(\Map(\rho^{(n)},F,\delta,\sigma)} )))\\&\leq \sum_{n=1}^N D(\pi,\beta_n, \rho_X^{(n)}, F, \delta, \sigma),
		\end{align*}
		and hence
		$$\mdim_{\Sigma} (\pi,\alpha,\rho_X)\leq\mdim_{\Sigma} (\pi,\alpha,\rho_X, F,\d')\leq \sum_{n=1}^N\mdim_{\Sigma} (\pi,\alpha,\rho_X^{(n)}, F,\d).$$
		Since $F$ and $\d$ are arbitrary, we have
		$$\mdim_{\Sigma} (\pi,\alpha,\rho_X)\leq \sum_{n=1}^N\mdim_{\Sigma} (\pi,\alpha,\rho_X^{(n)}, F,\d)\leq \sum_{n=1}^N \mdim_{\Sigma}(Y_n|X_n)\leq \sum_{1 \leq n < R} \mdim_{\Sigma}(Y_n|X_n)$$	
		Then
		\[
		\mathrm{mdim}_{\Sigma} (Y|X) \leq \sum_{1 \leq n < R} \mathrm{mdim}_{\Sigma} (Y_n|X_n).
		\]	
		The proof is complete.
	\end{proof}
	
	\begin{lem}\label{factor2}
		Let $\pi: (X,G)\to (Y,G)$ be a factor map between $G$-systems 
		with \(\mathrm{mdim}_{\Sigma} (Y|X) \geq 0\) and $\Sigma$ a sofic approximation sequence for $G$. Then for any factor $\psi:(Z,G)\to(Y,G)$ of $\pi$, one has $\mathrm{mdim}_{\Sigma} (Y|Z) \geq 0$.
	\end{lem}
	\begin{proof}
		Since \(\mathrm{mdim}_{\Sigma} (Y|X) \geq 0\), for any finite subset $F$ of $G$, any $\d > 0$, and any $N\in \N$, there is some $i\ge N$ such that
		$\mathrm{Map}(\rho_X, F, \delta, \sigma_i)$ is nonempty. For any factor $\psi:(Z,G)\to(Y,G)$ of $\pi$, there exists a factor map $\phi:(X,G)\to(Z,G)$ with $\pi=\psi\circ\phi$. Then by Lemma~\ref{li-2}, there exists $\delta' > 0$ such that for every map $\sigma$ from $G$ to $\operatorname{Sym}(d)$ for some $d \in \mathbb{N}$ and every $\xi \in \operatorname{Map}(\rho_X, F, \delta', \sigma_i)$, one has $\phi \circ \xi \in \operatorname{Map}(\rho_Z, F, \delta, \sigma_i)$. Then $\mathrm{mdim}_{\Sigma} (Y|Z) \geq 0.$
	\end{proof}

	For relative sofic mean dimension, one can also define the notion of relative zero sofic mean dimension extension.	
	
	\begin{defn}
		Let $\pi: (X,G)\to (Y,G)$ be a factor map between $G$-systems and $\Sigma$ a sofic approximation sequence for $G$. We say that $\pi$ is a \emph{relative zero sofic mean dimension extension (with respect to \(\Sigma\) and $X$)} if $\mathrm{mdim}_{\Sigma} (Y|X)=0$.
	\end{defn}

	\begin{cor}\label{prod2}
		Let $\Sigma$ be a sofic approximation sequence for $G$ and \(A=\{\psi_n : (X_n,G) \to (Z,G)\}_{n \in \N}\), where $\psi_n$ is a relative zero sofic  mean dimension extensions of $(Z,G)$ with respect to $X_n$ for any $n\in\N$. Let $\psi: (X,G)\to (Z,G)$ be the product relative to $(Z,G)$ of this set. Then $\mathrm{mdim}_{\Sigma} (Z|X) \leq 0$. 
	\end{cor}
	\begin{proof}
		By the definition, $$X= \left\{ (x_n)_{n\in\N} \in \prod_{n \in \N} X_n : \exists z \in Z \ \forall n \in \N, \psi_n(x_n) = z \right\}\subset\prod_{n\in\N}X_n,$$
		and $\psi(X)=\Delta_{\prod_{n\in\N}Z}:=\{(z)_{n\in\N}\in\prod_{n\in\N}Z:z\in Z\}\subset\prod_{n\in\N}Y_n$. By Proposition~\ref{subset2} and Proposition~\ref{prod2}, one has $\mathrm{mdim}_{\Sigma} (\Delta_{\prod_{n\in\N}Z}|X) \leq 0$. It is clear that $(\Delta_{\prod_{n\in\N}Z},G)$ is conjugate to $(Z,G)$. Then by Proposition \ref{isomorphism2}, we have $\mathrm{mdim}_{\Sigma} (Z|X)=\mathrm{mdim}_{\Sigma} (\Delta_{\prod_{n\in\N}Z}|X) \leq 0$.
	\end{proof}	
	Following the idea in Theorem~\ref{main1}, we also have maximal relative zero sofic mean dimension factor by Proposition~\ref{subset2}, Lemma \ref{factor2}, and Corollary~\ref{prod2}.
	\begin{thm}\label{main2}
		Let $\pi: (X,G)\to (Y,G)$ be a factor map between $G$-systems 
		with \(\mathrm{mdim}_{\Sigma} (Y|X) \geq 0\) and $\Sigma$ a sofic approximation sequence for $G$. Then there is a  maximal relative zero sofic mean dimension factor $\psi:(X_{Y|X}^{\Sigma}, G)\to(Y,G)$ of $\pi$ (with respect to \(\Sigma\) and X). That is, $\psi:(X_{Y|X}^{\Sigma}, G)\ra (Y,G)$ is a relative zero  sofic mean dimension extension (with respect to \(\Sigma\) and $X_{Y|X}^{\Sigma}$)
		and a factor of $\pi$;
		and every  relative zero sofic mean dimension   extension $\psi':(Z,G)\ra (Y,G)$ of $\pi$ (with respect to \(\Sigma\) and $Z$) is a factor of $\psi$.
	\end{thm}

	\section{Relative mean dimension tuples, sofic \texorpdfstring{$Y| X$-}-CPMD and sofic \texorpdfstring{$Y| X$-}-UPMD}\label{sec6}
	
	Like  sofic conditional mean dimension tuples, completely positive sofic conditional mean dimension and uniform positive sofic conditional mean dimension, we study relative sofic mean dimension tuples, completely positive relative sofic mean dimension and uniform positive relative sofic mean dimension, respectively, in this section.
	
	\subsection{Relative sofic mean dimension tuples}

	\begin{defn}
		Let $\pi \colon (X, G) \to (Y, G)$ be a nontrivial factor map between $G$-systems and $\Sigma$ a sofic approximation sequence for $G$. A tuple $(y_i)_{i=1}^n\in Y^n$ is said to be a \emph{relative sofic mean dimension tuple relevant to $\Sigma$ and $X$} if for every admissible open cover $\alpha$ with respect to $(y_i)_{i=1}^n$, we have $\text{mdim}_{\Sigma}(\pi,\alpha)>0$. For $n \geq 2$, denote by $\mathrm{D}_{n}^{\mathrm{md}}(Y|X, G,\Sigma)$ the set of all sofic conditional mean dimension $n$-tuples  relevant to $\Sigma$ and $X$.
	\end{defn}
	\begin{rem}
		When $\pi$ is a homeomorphism, the set of all sofic conditional mean dimension $n$-tuples $\mathrm{D}_{n}^{\mathrm{md}}(Y|X, G,\Sigma)$ and the set of all sofic  mean dimension $n$-tuples $\mathrm{D}_{n}^{\mathrm{md}}(Y, G,\Sigma)$ coincide.
	\end{rem}
	
	Similar to the sofic conditional mean dimension tuple,	relative sofic mean dimension tuple also has the following result by the idea in the proof of Proposition \ref{jajajajaj}.
	\begin{prop}\label{admis-Y}
		Let $\pi:(X,G)\to (Y,G)$ be a factor map between $G$-systems, $n \geq 2$ and $\Sigma$ a sofic approximation sequence for $G$. Then 
		$(y_i)_{i=1}^n \in Y^n\in Y^n\setminus\Delta_n(Y)$ is a relative sofic mean dimension tuple relevant to $\Sigma$ and $X$ if and only if for all open covers to be of the form $\alpha=(U_1,\dots,U_n)$, where $U_i^c$ is a neighborhood of $x_i$ such that $U_i^c\cap U_j^c=\emptyset,$ when $y_i\neq y_j,1\leq i<j\leq n,$ we have $\mathrm{mdim}_{\Sigma}(\pi,\alpha)>0.$
	\end{prop}

	\begin{prop}\label{Y-admi}
		Let $\pi:(X,G)\to (Y,G)$ be a factor map between $G$-systems, $n \geq 2$ and $\Sigma$ a sofic approximation sequence for $G$. If $\mathrm{mdim}_{\Sigma}(Y|X,G)>0$, then there exists an admissible open cover $\beta=(V_1,V_2,\dots,V_n)$ of $Y$ with respect to some $(y_i)_{i=1}^n\in Y^n$ such that $\mathrm{mdim}_{\Sigma}(\pi,\beta)>0.$
	\end{prop}

	\begin{proof}
		Combining with Proposition \ref{jajajajaj}, using the same method as Proposition \ref{propo322}, we can find an admissible open cover $\beta=(V_1,V_2,\dots,V_n)$ with respect to $(y_i)_{i=1}^{n}\in Y^n$ and $\text{mdim}_{\Sigma}(\pi,\beta)>0.$
	\end{proof}
	\begin{prop}\label{proooo325}
		Let $\pi:(X,G)\to (Y,G)$ be a factor map between $G$-systems, $n \geq 2$ and $\Sigma$ a sofic approximation sequence for $G$. For any finite open cover $\alpha=(U_1,U_2,\dots,U_n)$ of $Y$ satisfying $\mathrm{mdim}_{\Sigma}(\pi,\alpha)>0$, there exist $y_i\in U_i^c, i=1,\dots,n$ such that $(y_i)_{i=1}^{n}\in \mathrm{D}_{n}^{\mathrm{md}}(Y|X, G,\Sigma).$   
	\end{prop}
	\begin{proof}
		Combining with Proposition \ref{jajajajaj}, using the same method as Proposition \ref{ppppro324}, we can find $y_i\in U_i^c, i=1,\dots,n$ such that $(y_i)_{i=1}^{n}\in \mathrm{D}_{n}^{\mathrm{md}}(Y|X, G,\Sigma).$
	\end{proof}

	\begin{thm}\label{main4}
		Let $\pi:(X,G)\to (Y,G)$ be a factor map between $G$-systems, $n \geq 2$ and $\Sigma$ a sofic approximation sequence for $G$. Then $\mathrm{mdim}_{\Sigma}(Y|X)>0$ if and only if $\mathrm{D}_{n}^{\mathrm{md}}(Y|X, G,\Sigma)\neq \emptyset,$ for some $n\ge 2$.
	\end{thm}
	\begin{proof}
		Assume $\mathrm{D}_{n}^{\mathrm{md}}(Y|X, G,\Sigma)\neq \emptyset$, for some $n\ge 2$. Let $(y_i)^{n}_{i=1}\in \mathrm{D}_{n}^{\mathrm{md}}(Y|X, G,\Sigma)$. Let $\varepsilon>0$ such that $\bigcap_{i=1}^n\overline{B_{\v}(y_i)}=\emptyset.$ Define $U_i=Y\setminus \overline{B_{\v}(y_i)}$, then $\bigcup_{i=1}^nU_i=Y.$ Note that $B_{\v}(y_i)\subset U_i^c$ for every $i=1,2,\dots,n$, and $U_i\cap U_j\neq \emptyset$ for $i\neq j.$ Then $(U_1,\dots,U_n)$ is an admissible open cover of $Y$ with respect to $(y_i)_{i=1}^n$, which implies $\text{mdim}_{\Sigma}(\pi,(U_1,\dots,U_n))>0$. Hence 
		$$\text{mdim}_{\Sigma}(Y|X)\geq \text{mdim}_{\Sigma}(\pi,(U_1,\dots,U_n))>0.$$
		
		Now suppose that $\text{mdim}_{\Sigma}(Y|X)>0$. Combining Propositions \ref{Y-admi} and \ref{proooo325}, we get the desired result.
	\end{proof}
	
	Similar to the proof of Proposition~\ref{closed}, we also have the following result.

	\begin{prop}\label{closed2}
		Let $\pi:(X,G)\to (Y,G)$ be a factor map between $G$-systems, $n \geq 2$ and $\Sigma$ a sofic approximation sequence for $G$. 
		Then 
		$\mathrm{D}_{n}^{\mathrm{md}}(Y|X,G,\Sigma)\cup \Delta_n(Y)$ is a closed set and $\overline{\mathrm{D}_{n}^{\mathrm{md}}(Y|X,G,\Sigma)}\subset \mathrm{D}_{n}^{\mathrm{md}}(Y|X,G,\Sigma)\cup \Delta_n(Y).$
	\end{prop}

	\begin{lem}\label{lemmmm31}
		Let $\pi:(X,G)\to (Y,G)$ be a factor map between $G$-systems and $\Sigma$ a sofic approximation sequence for $G$. Let $\psi:(Z, G)\to(Y,G)$  a factor map of $\pi$ and $\phi:(X, G)\to(Z,G)$ a factor map with $\pi=\psi\circ\phi$. Then for any  finite open cover $\alpha$ of $Y$, one has $\mathrm{mdim}_{\Sigma}(\phi,\psi^{-1}(\alpha))\leq \mathrm{mdim}_{\Sigma}(\pi,\alpha).$
	\end{lem}
	\begin{proof}
		For every open cover $\alpha$ of $Y^d$ and $z \in Z^d$, it holds 
		$$\sum_{V\in \alpha}1_{V}(\psi^d(z))=\sum_{V\in (\psi^d)^{-1}(\alpha)}1_{V}(z).$$
		Let $F \in \mathcal{F}(G), \delta>0$ and $\sigma:G \rightarrow \Sym(d)$ for some $d \in \N$. Note that \begin{align*}
			\psi^d\circ \phi^d(\text{Map}(\rho_X,F,\delta,\sigma))&=(\psi\circ \phi)^d(\text{Map}(\rho_X,F,\delta,\sigma))\\
			&=\pi^d(\text{Map}(\rho_X,F,\delta,\sigma)).
		\end{align*}
		Assume that $\xi \in \Map(\rho_X,F,\delta,\sigma)$ 
		Then 
		\begin{align*}
			&\min_{\beta\succ(\psi^{-1}(\alpha))^{d}}\left(\max_{\phi\circ \xi\in \phi^d( \text{Map}(\rho_X,F,\delta,\sigma))}\sum_{U\in \beta}1_U(\phi\circ \xi)\right)\\
			&\leq \min_{\beta\succ\alpha^{d}}\left(\max_{\phi\circ \xi\in \phi^d( \text{Map}(\rho_X,F,\delta,\sigma))}\sum_{U\in \beta}1_U(\psi\circ \phi\circ \xi)\right)\\
			&= \min_{\beta\succ\alpha^{d}}\left(\max_{\psi\circ \phi (\xi)=\pi\circ \xi\in \pi^d( \text{Map}(\rho_X,F,\delta,\sigma))}\sum_{U\in \beta}1_U(\pi\circ \xi)\right).
		\end{align*}
		Thus $$D(\phi,\psi^{-1}(\alpha),\rho_X,F,\delta,\sigma)\leq D(\pi,\alpha,\rho_X,F,\delta,\sigma).$$
		This implies that $$\text{mdim}_{\Sigma}(\phi,\psi^{-1}(\alpha))\leq \text{mdim}_{\Sigma}(\pi,\alpha).$$
		The proof is complete.
	\end{proof}

	\begin{prop}\label{liftY}
		Let $\pi:(X,G)\to (Y,G)$ be a factor map between $G$-systems, $n \geq 2$ and $\Sigma$ a sofic approximation sequence for $G$. Let $\psi:(Z, G)\to(Y,G)$ be a factor map of $\pi$ and $\phi:(X, G)\to(Z,G)$ a factor map with $\pi=\psi\circ\phi$. If $(z_i)_{i=1}^n\in \mathrm{D}_{n}^{\mathrm{md}}(Z|X,G,\Sigma)$ with $(\psi(z_i))_{i=1}^n\not\in\Delta_n(Y)$, then $(\psi(z_i))_{i=1}^n\in \mathrm{D}_{n}^{\mathrm{md}}(Y|X,G,\Sigma).$
	\end{prop}

	\begin{proof}
		Let $(z_i)_{i=1}^n\in \mathrm{D}_{n}^{\mathrm{md}}(Z|X,G,\Sigma)$ with $(\psi(z_i))_{i=1}^n\not\in\Delta_n(Y)$. Suppose that $\alpha=(U_1,\dots,U_n)$ is an admissible open cover of $Y$ with respect to $(\psi(z_i))_{i=1}^n,$ where $U_i^c$ is a closed neighborhood of $\psi(z_i)$ such that $U_i^c\cap U_j^c=\emptyset$ when $\psi(z_i)\neq \psi(z_j)$, $1\leq i,j\leq n.$ By the continuity of $\psi,$ we have $(\psi^{-1}(U_i))^c=\psi^{-1}(U_i^c)$ is a closed neighborhood of $z_i$ for every $i=1,2\dots,n.$ Moreover,
		$$(\psi^{-1}(U_i))^c\cap (\psi^{-1}(U_j))^c=\psi^{-1}(U_i^c)\cap \psi^{-1}(U_j^c)=\psi^{-1}(U_i^c\cap U_j^c)=\emptyset,$$
		when $x_i\neq x_j,1\leq i,j\leq n.$ Note that $\psi^{-1}(\alpha)$ is a finite open cover of $Z.$ Thus $\psi^{-1}(\alpha)=(\psi^{-1}(U_1),\psi^{-1}(U_2),\dots,\psi^{-1}(U_n))$ is an admissible open cover of $Z$ with respect to $(z_i)_{i=1}^n$. By Lemma \ref{lemmmm31},
		$$0<\text{mdim}_{\Sigma}(\phi,\psi^{-1}(\alpha))\leq \text{mdim}_{\Sigma}(\pi,\alpha).$$
		Combining with Proposition \ref{admis-Y}, we have $(\psi(z_i))_{i=1}^n\in \mathrm{D}_{n}^{\mathrm{md}}(Y|X,G,\Sigma).$
	\end{proof}	
	\begin{prop}
		Let $\pi_1:(X,G)\to (Y,G)$ be a factor map between $G$-systems, $n \geq 2$ and $\Sigma$ a sofic approximation sequence for $G$. Let $Z\subset Y$ be a closed $G$-invariant subset and $\pi_2=\pi_1|_Z:(\pi_2^{-1}(Z),G)\to (Z,G)$ a factor map between $G$-systems.  If $(y_i)_{i=1}^n\in \mathrm{D}_{n}^{\mathrm{md}}(Z|\pi_2^{-1}(Z), G,\Sigma),$ then $(y_i)_{i=1}^n\in \mathrm{D}_{n}^{\mathrm{md}}(Y|X, G,\Sigma)$.
	\end{prop}

	\begin{proof}
		Let $\alpha$ be an admissible open cover of $y$ with respect to $(y_i)_{i=1}^n\in (\pi^{-1}(Z))^n\setminus \Delta_n(\pi^{-1}(Z))$. Clearly, $\alpha|_Z$, its restriction to $Z$, is an admissible open cover of $Z$ with respect to $(y_i)_{i=1}^n$. Since $\text{mdim}_{\Sigma}(\pi,\alpha)\geq \text{mdim}_{\Sigma}(\pi,\alpha|_Z).$ Combining with $\text{mdim}_{\Sigma}(\pi,\alpha|_Z)>0$, we get the desire result.
	\end{proof}

	\subsection{Relative completely/uniform positive sofic mean dimension} 
	
	This subsection will start with introducing the notion of relative completely positive sofic mean dimension and relative uniform positive sofic mean dimension, and study the properties of them.
	
	\begin{defn}
		Let $\pi:(X,G)\to (Y,G)$ be a nontrivial factor map between $G$-systems and $\Sigma$ a sofic approximation sequence for $G$. We say $\pi$ has \emph{relative completely positive sofic mean dimension relevant to $X$} (sofic $Y|X$-CPMD) if for any proper factor $(Z,G)$ of $(Y,G)$, one has $\mdim_{\Sigma}(Z|X)>0$. 
	\end{defn}

	\begin{defn}
		Let $\pi:(X,G)\to (Y,G)$ be a nontrivial factor map between $G$-systems and $\Sigma$ a sofic approximation sequence for $G$. We say $\pi$ has \emph{relative uniform positive sofic mean dimension relevant to $X$} (sofic $Y|X$-UPMD) if for any non-dense open cover $\alpha$, one has $\mdim_{\Sigma}(\pi,\alpha)>0$. 
	\end{defn}

	\begin{prop}\label{dyx}
		Let $\pi:(X,G)\to (Y,G)$ be a factor map between $G$-systems and $\Sigma$ a sofic approximation sequence for $G$. Then  $\pi$ has sofic $Y|X$-UPMD if and only if $\mathrm{D}_{n}^{\mathrm{md}}(Y|X,G)=Y^n\setminus \Delta_n(Y)$, for some $n\ge 2$.
	\end{prop}
	\begin{proof}
		Assume that  $\pi$ has sofic $Y|X$-UPMD. Let $(y_i)_{i=1}^n\in Y^n\setminus \Delta_n(Y)$, for every admissable open cover $\alpha$ with respect to $(y_i)_{i=1}^n$, it holds $\text{mdim}_{\Sigma}(\pi,\alpha)>0$. This implies that $(y_i)_{i=1}^n\in \mathrm{D}_{n}^{\mathrm{md}}(Y|X, G,\Sigma)$. Hence $\mathrm{D}_{n}^{\mathrm{md}}(Y|X, G,\Sigma)=Y^n\setminus \Delta_n(Y).$

		On the other hand, assume that $\mathrm{D}_{n}^{\mathrm{md}}(Y|X, G,\Sigma)=Y^n\setminus \Delta_n(Y),$ for some $n\ge 2$. Since every admissible open cover with respect to some  $(y_i)_{i=1}^n\in Y^n\setminus \Delta_n(Y)$, then for all non-dense open covers $\alpha$, it holds $\text{mdim}_{\Sigma}(\pi,\alpha)>0$, which implies the desired result.
	\end{proof}

	\begin{thm}\label{univ-Y}
		Let $\pi:(X,G)\to (Y,G)$ be a factor map between $G$-systems with \(\mathrm{mdim}_{\Sigma} (Y|X) \geq 0\) and $\Sigma$ a sofic approximation sequence for $G$. If the smallest closed \(G\)-invariant equivalence relation that contains $\mathrm{D}_{2}^{\mathrm{md}}(Y|X, G,\Sigma)$ is $Y\times Y$, then $\pi$ has sofic $Y|X$-CPMD. 
		
	\end{thm}
	\begin{proof}	
		Assume that $\pi$ does not have sofic $Y|X$-CPMD. Then there exists a proper factor $(Z,G)$ of $(Y,G)$ with $\mdim_{\Sigma}(Z|X)\leq0$ and $\phi \colon (Y, G) \to (Z, G)$ is the nontrivial factor map. By Theorem~\ref{main4}, we have $\mathrm{D}_{2}^{\mathrm{md}}(Z|X, G,\Sigma)=\emptyset$.  By Proposition~\ref{liftY}, we have $\phi\times\phi(\mathrm{D}_{2}^{\mathrm{md}}(Y|X, G,\Sigma))\subset \Delta(Z)$. That is $R_\phi\supset \mathrm{D}_{2}^{\mathrm{md}}(Y|X, G,\Sigma)$. Then $R_\phi=Y\times Y$. A contradiction. Thus $\pi$ has sofic $Y|X$-CPMD.
	\end{proof}
	
	\begin{cor}
		Let $\pi \colon (X, G) \to (Y, G)$ be a nontrivial factor map between $G$-systems  and $\Sigma$ a sofic approximation sequence for $G$. If $\pi$ has sofic $Y|X$-UPMD, then $\pi$ has sofic $Y|X$-CPMD.
	\end{cor}
	\begin{proof}
		Assume the factor map $\pi \colon (X, G) \to (Y, G)$ has sofic $Y|X$-UPMD.
		By Proposition \ref{dyx}, 	$\mathrm{D}_{2}^{\mathrm{md}}(Y|X, G,\Sigma) = Y^2\setminus \Delta_2(Y)$.
		Then by Theorem~\ref{univ-Y}, $\pi$ has sofic $Y|X$-CPMD.

	\end{proof}
	
	Now we will investigate the relationship between sofic conditional mean dimension tuple and relative sofic mean dimension tuple. To start with, we need the following lemma.	
	\begin{lem}\label{lllllem344}
		Let $\pi:(X,G)\to (Y,G)$ be a factor map between $G$-systems , $\Sigma$ a sofic approximation sequence for $G$. Then for any  finite open cover $\alpha$ of $Y$, one has $$\mathrm{mdim}_{\Sigma}(\pi^{-1}(\alpha)|\pi)\leq \mathrm{mdim}_{\Sigma}(\pi^{-1}(\alpha))\leq \mathrm{mdim}_{\Sigma}(\pi,\alpha)\leq \mathrm{mdim}_{\Sigma}(\alpha).$$
	\end{lem}
	\begin{proof}
		Let $F\in \mathcal{F}(G)$, $\delta>0$, and $\sigma:G\to \text{Sym}(d)$ for some $d \in \N$, we have 
		\begin{align*}
			&\min_{\{(\pi^d)^{-1}(y)\}_{y\in Y^d}\vee \beta\succ (\pi^{-1}(\alpha))^d} \left(\max_{\xi\in \text{Map}(\rho_X,F,\delta,\sigma)}\sum_{U\in \beta}1_U(\xi)\right)\\
			&\leq\min_{\beta\succ(\pi^{-1}(\alpha))^d}\left(\max_{\xi\in \text{Map}(\rho_X,F,\delta,\sigma)}\sum_{U\in \beta}1_U(\xi)\right)\\
			&\leq\min_{\beta\succ\alpha^d}\left(\max_{\xi\in \text{Map}(\rho_X,F,\delta,\sigma)}\sum_{U\in \beta}1_U(\pi\circ \xi)\right)\\
			&=\min_{\beta\succ\alpha^d}\left(\max_{\pi\circ \xi\in \pi^d(\text{Map}(\rho_X,F,\delta,\sigma))}\sum_{U\in \beta}1_U(\pi\circ \xi)\right)\\
			&=\min_{\beta\succ\alpha^d}\left(\max_{\xi'\in \pi^d(\text{Map}(\pi,\rho_X,F,\delta,\sigma))}\sum_{U\in \beta}1_U(\xi')\right).
		\end{align*}
		Then
		\begin{align*}
			D((\pi^{-1}(\alpha))^d|_{\text{Map}(\rho_X,F,\delta,\sigma)}|\pi^d)&=D(\pi^{-1}(\alpha)|\pi,\rho_X,F,\delta,\sigma)\\
			&\leq D(\alpha,\rho_X,F,\delta,\sigma)\\
			&\leq D(\alpha^d|_{\pi^d(\text{Map}(\rho_X,F,\delta,\sigma)})\\
			&=D(\pi,\alpha,\rho_X,F,\delta,\sigma),
		\end{align*}
		which implies that $$\text{mdim}_{\Sigma}(\pi^{-1}(\alpha)|\pi)\leq \text{mdim}_{\Sigma}(\pi^{-1}(\alpha))\leq \text{mdim}_{\Sigma}(\pi,\alpha).$$
		On the other hand, by Lemma \ref{li-2}, there exists  $\delta'>0$, such that
		$$\pi^d(\Map(\rho_X,F,\delta',\sigma)) \subset \Map(\rho_Y,F,\delta,\sigma).$$
		This implies that
		$$D(\pi,\alpha,\rho_X,F,\delta',\sigma)\leq D(\alpha,\rho_Y,F,\delta,\sigma).$$
		Hence, we get $ \text{mdim}_{\Sigma}(\pi,\alpha) \leq \mathrm{mdim}_{\Sigma}(\alpha).$
	\end{proof}

	\begin{thm}
		Let $\pi:(X,G)\to (Y,G)$ be a factor map between $G$-systems, $n \geq 2$ and $\Sigma$ a sofic approximation sequence for $G$. If $(x_i)_{i=1}^{n}\in \mathrm{D}^{\mathrm{md}}_{n}(X|\pi,G,\Sigma)$ with $(\pi(x_i))_{i=1}^{n}\not\in\Delta_n(Y)$, then $(\pi(x_i))_{i=1}^{n}\in \mathrm{D}^{\mathrm{md}}_{n}(Y|X,G,\Sigma)$. Moreover,
		\begin{align*}
			\pi^n( \mathrm{D}^{\mathrm{md}}_{n}(X|\pi,G,\Sigma))\subset \pi^n( \mathrm{D}^{\mathrm{md}}_{n}(X,G,\Sigma))\subset \mathrm{D}^{\mathrm{md}}_{n}(Y|X,G,\Sigma)\subset \mathrm{D}^{\mathrm{md}}_{n}(Y,G,\Sigma).
		\end{align*}
	\end{thm}
	\begin{proof}
		By Lemma \ref{lllllem344}, we only need to show that
		$$\pi^n( \mathrm{D}^{\mathrm{md}}_{n}(X,G,\Sigma))\subset \mathrm{D}^{\mathrm{md}}_{n}(Y|X,G,\Sigma).$$
		
		Let $(x_i)_{i=1}^{n}\in \mathrm{D}^{\mathrm{md}}_{n}(X,G,\Sigma)$  with $(\pi(x_i))_{i=1}^{n}\not\in\Delta_n(Y)$. Suppose that $\alpha$ is an admissible open cover of $Y$ with respect to $(\pi(x_i))_{i=1}^{n}$, then $\pi^{-1}(\alpha)$ is  an admissible open cover of $X$ with respect to $(x_i)_{i=1}^{n}$. By Lemma \ref{lllllem344},
		$$0<\text{mdim}_{\Sigma}(\pi^{-1}(\alpha))\leq \text{mdim}_{\Sigma}(\pi,\alpha),$$
		which implies that $(\pi(x_i))_{i=1}^{n}\in \mathrm{D}^{\mathrm{md}}_{n}(Y|X,G,\Sigma).$
	\end{proof}
	
	\section{Sofic conditional  mean dimension sets and relative sofic mean dimension sets }\label{sec7}
	
	In this section, we will introduce and study the notions of sofic conditional  mean dimension sets and relative sofic mean dimension sets. We begin by presenting their formal definitions.
	
	\begin{defn}
		Let $\pi:(X,G)\to (Y,G)$ be a factor map  between $G$-systems and $\Sigma$ a sofic approximation sequence for $G$. We say that a subset $K$ of $X$ is a \emph{sofic conditional  mean dimension set (of $X$ relevant to $\pi$)} if $\#(K)\geq 2$ and if for any admissible finite open cover $\alpha$ with respect to $K$, one has $\text{mdim}(\alpha|\pi)>0.$ 
	\end{defn}
	
	\begin{defn}
		Let $\pi:(X,G)\to (Y,G)$ be a factor map between $G$-systems and $\Sigma$ a sofic approximation sequence for $G$. We say that a subset $K$ of $Y$ is a \emph{relative sofic mean dimension set (of $Y$ relevant to $X$)} if $\#(K)\geq 2$ and if for any admissible finite open cover $\alpha$ with respect to $K$, one has $\text{mdim}(\pi,\alpha)>0.$ 
	\end{defn}
	
	The following conclusion shows the relationship between sofic conditional  mean dimension sets and sofic conditional mean dimension tuples, and the relationship between relative sofic mean dimension sets and relative sofic mean dimension tuples.

	\begin{thm}\label{ttttthm39}
		Let $\pi:(X,G)\to (Y,G)$ be a factor map between $G$-systems and $\Sigma$ a sofic approximation sequence for $G$. The following holds:
		\begin{enumerate}
			\item[(1)] A subset $K$ of $X$ is a sofic conditional mean dimension set (of $X$ relative to $\pi$) if and only if $\#(K)\geq 2$ and if for any distinct $n$ points $\{x_1,x_2,\dots,x_n\}\subset K$, the tuple $(x_i)_{i=1}^{n}$ is a sofic conditional mean dimension tuple (of $X$ relevant to $\pi$).
			\item[(2)] A subset $K$ of $Y$ is a relative sofic mean dimension set (of $Y$ relative to $X$) if and only if $\#(K)\geq 2$ and if each $n\geq 2$ and for any distinct $n$ points $\{y_1,y_2,\dots,y_n\}\subset K$, the tuple $(y_i)_{i=1}^{n}$ is a relative sofic mean dimension tuple (of $Y$ relevant to $X$).
		\end{enumerate}
		
	\end{thm}
	\begin{proof}
		We only prove (1), the proof of (2) is similar.
		
		If $K\subset X$ is a sofic conditional mean dimension set, then for each $n\geq 2$ with respect to $\{x_1,x_2,\dots,x_n\}\subset K$, and any admissible finite open cover $\alpha$ of $X$ with respect to $(x_i)_{i=1}^{n}$ is also an admissible finite open cover of $X$ with respect to $K$. Thus $\text{mdim}_{\Sigma}(\alpha|\pi)>0$ and the tuple $(x_i)_{i=1}^{n}$  is a sofic conditional  mean dimension tuple.
		
		On the other hand, for each $n\geq 2,$ let $\alpha=\{U_1,U_2,\dots,U_n\}$ be any admissible finite open cover of $X$ with respect to $K$. Then there exists $x_i\in K$ such that $x_i\not\in \overline{U_i}, 1\leq i\leq n.$ Notice that $\alpha$ is an  admissible finite open cover of $X$ with respect to $(x_i)_{i=1}^{n}$, we have $\text{mdim}_{\Sigma}(\alpha|\pi)>0$.
	\end{proof}
	
	\begin{thm}
		Let $\pi:(X,G)\to (Y,G)$ be a factor map between $G$-systems and $\Sigma$ a sofic approximation sequence for $G$. The following holds:
		\begin{enumerate}
			\item [(1)] If $K$ is a sofic conditional mean dimension set, then $\overline{K}$ is a sofic conditional mean dimension set.
			\item [(2)] If $K$ is a relative sofic mean dimension set, then $\overline{K}$ is a relative sofic mean dimension set.
		\end{enumerate}
		
	\end{thm}
	\begin{proof}
		(1). Let $n\geq 2$ and $\{x_1,x_2,\dots,x_n\}\subset \overline{K}$ with $x_i\neq x_j$ if $x\neq j.$ Assume that for each $1\leq j\leq n$, $x_j^k\to x_j$ when $k\to \infty$ and $x_j^k\in K.$ Then by Theorem \ref{ttttthm39}, $(x_i^k)_{i=1}^n\in D^{\mathrm{md}}_{n}(X|\pi,G,\Sigma)$ when $k$ is sufficiently large. By Proposition \ref{closed}, $D^{\mathrm{md}}_{n}(X|\pi,G,\Sigma)\cup \Delta_n(X)$ is a closed set, then we have $(x_i)_{i=1}^{n}\in D^{\mathrm{md}}_{n}(X|\pi,G,\Sigma).$ Therefore, $\overline{K}$ is a sofic conditional mean dimension set.
		
		(2). Similarly, combining with Proposition \ref{closed2}, we can prove (2).
	\end{proof}
	
	\begin{thm}\label{ttttthmm40}
		Let $\pi:(X,G)\to (Y,G)$ be a factor map between $G$-systems and $\Sigma$ a sofic approximation sequence for $G$. Let $\psi:(Z, G)\to(Y,G)$  a factor map of $\pi$ and $\phi:(X, G)\to(Z,G)$ a factor map with $\pi=\psi\circ\phi$. Then the image of the sofic conditional mean dimension set of $X$ relevant to $\pi$ under $\phi$ is a sofic conditional mean dimension set of $Z$ relative to $\psi$ if it does not collapse to a single point.
	\end{thm}
	\begin{proof}
		Let $K$ be a sofic conditional mean dimension set of $X$ relevant to $\pi$, by Theorem \ref{ttttthm39}, for any distinct $n$ points $\{x_1,x_2,\dots,x_n\}\subset K$, we have $(x_i)_{i=1}^{n}\in D^{\mathrm{md}}_{n}(X|\pi,G,\Sigma)$. By Proposition \ref{lift}, if $(\phi(x_i))_{i=1}^{n}\not\in\Delta_n(Z)$, then $(\phi(x_i))_{i=1}^{n}\in D^{\mathrm{md}}_{n}(Z|\psi,G,\Sigma),$ which implies $\phi(K)$ is a sofic conditional mean dimension set of $Z$ relative to $\psi$. 
	\end{proof}
	
	\begin{thm}\label{ttttthmm41}
		Let $\pi:(X,G)\to (Y,G)$ be a factor map between $G$-systems and $\Sigma$ a sofic approximation sequence for $G$. Let $\psi:(Z, G)\to(Y,G)$  a factor map of $\pi$ and $\phi:(X, G)\to(Z,G)$ a factor map with $\pi=\psi\circ\phi$. Then the image of the relative sofic mean dimension set of $Z$ relevant to $X$ under $\psi$ is a relative sofic mean dimension set of $Y$ relative to $X$ if it does not collapse to a single point.
	\end{thm}

	\begin{proof}
		Let $K$ be a relative sofic mean dimension set of $Z$ relevant to $X$. By Proposition \ref{ttttthm39}, for any distinct $n$ points $\{z_1,z_2,\dots,z_n\}\subset K,$ the tuple  $(z_i)_{i=1}^{n}\in \mathrm{D}^{\mathrm{md}}_{n}(Z|X,G,\Sigma)$. By Theorem \ref{liftY}, if $(z_i)_{i=1}^{n}\in \mathrm{D}^{\mathrm{md}}_{n}(Z|X,G,\Sigma)$ with $(\psi(z_i))_{i=1}^{n}\notin\Delta_n(Y)$, then $(\psi(z_i))_{i=1}^{n}\in \mathrm{D}^{\mathrm{md}}_{n}(Y|X,G,\Sigma),$ which implies $\psi(K)$ is a relative sofic mean dimension set of $Y$
		relative to $X.$
	\end{proof}
	
	\begin{defn}
		Let $\pi:(X,G)\to (Y,G)$ be a factor map between $G$-systems and $\Sigma$ a sofic approximation sequence for $G$. Denote by $\mathrm{D}^{\mathrm{md}}_s(X|\pi,G,\Sigma)$ the set of all sofic conditional  mean dimension sets of $X$ relevant to $\pi$. Denote by $\mathrm{D}^{\mathrm{md}}_{s}(Y|X,G,\Sigma)$ the set of all relative sofic mean dimension sets  of $Y$ relevant to $X$.
	\end{defn}
	
	\begin{thm}
		Let $\pi:(X,G)\to (Y,G)$ be a factor map between $G$-systems and $\Sigma$ a sofic approximation sequence for $G$. The following holds:
		\begin{itemize}	\item[(1)]$\mathrm{D}^{\mathrm{md}}_{s}(X|\pi,G,\Sigma)\neq \emptyset$ if and only if $\mathrm{mdim}_{\Sigma}(X|\pi)>0$.
			\item [(2)]$\mathrm{D}^{\mathrm{md}}_{s}(Y|X,G,\Sigma)\neq \emptyset$ if and only if $\mathrm{mdim}_{\Sigma}(Y|X)>0$.
			\item [(3)] Let $\psi:(Z, G)\to(Y,G)$  a factor map of $\pi$ and $\phi:(X, G)\to(Z,G)$ a factor map with $\pi=\psi\circ\phi$. Then
			\begin{align*}
				&\phi(\mathrm{D}^{\mathrm{md}}_{s}(X|\pi,G,\Sigma))\subset \mathrm{D}^{\mathrm{md}}_{s}(Z|\psi,G,\Sigma)\cup \Delta_n(Z),\\
				&\psi(\mathrm{D}^{\mathrm{md}}_{s}(Z|X,G,\Sigma))\subset \mathrm{D}^{\mathrm{md}}_{s}(Y|X,G,\Sigma)\cup \Delta_n(Y).
			\end{align*}
		\end{itemize}
	\end{thm}
	\begin{proof}
		(1) is from Theorem \ref{main3}, (2) is from Theorem \ref{main4}, (3)
		is from Theorems \ref{ttttthmm40} and \ref{ttttthmm41} .\end{proof}

	\section{Conditional mean dimension for amenable group actions}\label{sec8}
	In this section, we discuss what happens when $G$ is an infinite amenable group. Firstly, we review the notion of conditional mean topological dimension for amenable groups.
	
	\begin{defn}
		Let $G$ be a countable discrete infinite group and $\F(G)$ the set of all finite nonempty subsets
		of $G$.
		A sequence $\{F_n\}_{n \geq 1}$ of nonempty finite subsets of $G$ is called a \emph{F{\o}lner sequence} for $G$ if one has
		$$
		\lim_{n \to \infty} \frac{|F_n \setminus g F_n|}{|F_n|} = 0.
		$$
		We say that a countable group is \emph{amenable} if it admits a F{\o}lner sequence.
	\end{defn}
	
	\begin{defn}
		Let $\pi:(X,G)\to (Y,G)$ be a factor map between $G$-systems, where $G$ is an infinite amenable group. We define the \emph{conditional mean topological dimension of $(X,G)$ relative to $(Y,G)$} as
		$$
		\operatorname{mdim}(X | Y):=\sup _{\alpha} \lim _{F \in \mathcal{F}(G)} \frac{D\left(\mathcal{\alpha}^F | \pi\right)}{|F|}
		$$
		for $\alpha$ running over all finite open covers of $X$, where $\alpha^F$ is the finite open cover $\bigvee_{s \in F}s^{-1}\alpha.$
	\end{defn}

	\begin{defn}
		Let $\pi:(X,G)\to (Y,G)$ be a factor map between $G$-systems, where $G$ is an infinite amenable group. Consider $\rho_Y$ as a continuous pseudometric on $Y$ and fix $\theta>0$. We say a continuous map $f: X \rightarrow P$ is an ( $\varepsilon, \theta, \rho_X | \rho_Y$ )-embedding (with respect to $\pi$ ) if
		$$
		\rho_X\left(x, x^{\prime}\right)<\varepsilon
		$$
		for every $x, x^{\prime} \in X$ satisfying that $f(x)=f\left(x^{\prime}\right)$ and $\rho_Y\left(\pi(x), \pi\left(x^{\prime}\right)\right) \leq \theta$. We call $f$ is an $\left(\varepsilon, \rho_X | Y\right)$-embedding if $\rho_X\left(x, x^{\prime}\right)<\varepsilon$ holds for every $x, x^{\prime} \in X$ satisfying that $f(x)=f\left(x^{\prime}\right)$ and $\pi(x)=\pi\left(x^{\prime}\right)$.
		
		Denote by $\operatorname{Wdim}_{\varepsilon}\left(\rho_X | \rho_Y, \theta\right)$ the minimal (covering) dimension $\operatorname{dim}(P)$ of a compact  metric space $P$ that admits an ( $\varepsilon, \theta, \rho_X | \rho_Y$ )-embedding from $X$ to $P$. Similarly, denote by $\operatorname{Wdim}_{\varepsilon}\left(\rho_X | Y\right)$ the minimal (covering) dimension $\operatorname{dim}(Q)$ of a compact  metric space $Q$ that admits an $\left(\varepsilon, \rho_X | Y\right)$-embedding from $X$ to $Q$. 
		
		For each $F \in \mathcal{F}(G)$ it induces a continuous pseudometric $\rho_{X, F}$ on $X$ via
		$$
		\rho_{X, F}\left(x, x^{\prime}\right):=\max _{s \in F} \rho_X\left(s x, s x^{\prime}\right).
		$$
		Write
		$$
		\operatorname{mdim}_{\varepsilon}\left(\rho_X | \rho_Y\right):=\inf _{\theta>0} \lim _{F \in \mathcal{F}(G)} \frac{\operatorname{Wdim}_{\varepsilon}\left(\rho_{X, F} | \rho_{Y, F}, \theta\right)}{|F|},
		$$
		and 
		$$
		\operatorname{mdim}_{\varepsilon}\left(\rho_X | Y\right):=\lim _{F \in \mathcal{F}(G)}  \frac{\operatorname{Wdim}_{\varepsilon}\left(\rho_{X, F} | Y\right)}{|F|}.
		$$
		
		For any continuous pseudometrics $\rho_X$ and $\rho_Y$ on $X$ and $Y$ respectively, we define
		$$
		\operatorname{mdim}\left(\rho_X | \rho_Y\right):=\sup _{\varepsilon>0} \operatorname{mdim}_{\varepsilon}\left(\rho_X | \rho_Y\right)
		$$
		and
		$$
		\operatorname{mdim}\left(\rho_X | Y\right):=\sup _{\varepsilon>0} \operatorname{mdim}_{\varepsilon}\left(\rho_X | Y\right).
		$$
		
	\end{defn} 
	The following result is from \cite[Proposition 2.5]{LBB1} and \cite[Proposition 3.4]{LBB2}.
	\begin{prop}\label{popopo}
		Let $\pi: (X,G) \rightarrow (Y,G)$ be a factor map between $G$-systems, where $G$ is an infinite amenable group. Let $\rho_X$ and $\rho_Y$ be  compatible metrics on $X$ and $Y$ respectively, then we have
		$$\operatorname{mdim}(X | Y)=\operatorname{mdim}\left(\rho_X | \rho_Y\right)=\operatorname{mdim}\left(\rho_X | Y\right).$$
	\end{prop}
	
	By adopting the argument of \cite[Theorem 4.7]{LBB2}, we prove that the sofic conditional mean dimension (which is defined by finite open covers) generalizes the notion of conditional mean dimension. 
	
	\begin{thm}\label{ooooooq}
		Let $\pi: (X,G) \rightarrow (Y,G)$ be a factor map between $G$-systems, where $G$ is an infinite amenable group. Then $\operatorname{mdim}_{\Sigma}(X |\pi)=\operatorname{mdim}(X|Y)$.
	\end{thm} 
	
	Theorem \ref{ooooooq} follows directly from Lemma \ref{wowowowowaa} and Lemma \ref{oaoaoaoaw}.
	\begin{lem}\label{wowowowowaa}
		Let $\pi: (X,G) \rightarrow (Y,G)$ be a factor map between $G$-systems, where $G$ is an infinite amenable group. Then 
		$\operatorname{mdim}_{\Sigma}(X|\pi ) \leq \operatorname{mdim}(X | Y)$.
	\end{lem} 
	\begin{proof}
		Let $\rho_X$ and $\rho_Y$ be  compatible metrics on $X$ and $Y$ respectively. Consider the  sofic conditional mean dimension $\operatorname{mdim}_{\Sigma}(\rho_X | \rho_Y)$(see \cite[Definition 4.2]{LBB2}),  it's obvious that 
		for any $\varepsilon>0$, we have 
		$$
		\operatorname{mdim}_{\Sigma,\varepsilon}(\rho_X | Y) \leq \operatorname{mdim}_{\Sigma}^{\varepsilon }\left(\rho_X | \rho_Y\right).
		$$
		By \cite[Lemma 4.9]{LBB2}, we have
		$$\operatorname{mdim}_{\Sigma}^{\varepsilon }\left(\rho_X | \rho_Y\right) \leq \operatorname{mdim}_{\varepsilon /3}\left(\rho_X | \rho_Y\right).$$
		Then 
		$$\operatorname{mdim}_{\Sigma,\varepsilon}(\rho_X | Y) \leq \operatorname{mdim}_{\varepsilon /3}\left(\rho_X | \rho_Y\right).$$
		Combining Proposition \ref{xoxoxoxo} and Proposition \ref{popopo}, we have
		$$\operatorname{mdim}_{\Sigma}(X| \pi)=\operatorname{mdim}_{\Sigma}(\rho_X | Y) \leq \operatorname{mdim}\left(\rho_X | \rho_Y\right)=\operatorname{mdim}\left(X |Y\right),$$
		which implies the desired result.
	\end{proof}
	
	Using the idea in the proof of \cite[Lemma 4.8]{LBB2}, one can get the following result.
	For the sake of completeness, we provide the proofs.

	\begin{lem}\label{oaoaoaoaw}
		Let $\pi: (X,G) \rightarrow (Y,G)$ be a factor map  where $G$ is an infinite amenable group. Let $\rho_X$ be a compatible metric on $X$. Then for every $\varepsilon>0$, we have $\operatorname{mdim}_{\Sigma,\varepsilon}\left(\rho_X |Y\right) \geq \operatorname{mdim}_{\varepsilon}\left(\rho_X |Y\right)$. Moreover, $\operatorname{mdim}_{\Sigma}(X |\pi) \geq \operatorname{mdim}(X | Y)$.
	\end{lem}
	\begin{proof}
		
		Without loss of generality, assume $\mdim_{\varepsilon}(\rho_X | Y) > 0$. Fix $0 < \beta < \mdim_{\varepsilon}(\rho_X | Y)$. To establish the desired inequality, it is sufficient to prove:
		$$
		\operatorname{mdim}_{\Sigma}^{\varepsilon}\left(\rho_X |Y\right) \geq \operatorname{mdim}_{\varepsilon}\left(\rho_X |Y\right)-2 \beta.
		$$
		Let $F \in \mathcal{F}(G)$ and $\delta>0$. We will show that for sufficiently good $\sigma$, the following holds:
		$$
		\mathrm{Wdim}_{\varepsilon}\left(\operatorname{Map}\left(\rho_X, F, \delta, \sigma\right),  \rho_{X, \infty} |Y\right) \geq\left(\operatorname{mdim}_{\varepsilon}\left(\rho_X |Y\right)-2 \beta\right)d.
		$$
		
		Choose a finite set $K \subset G$ containing $F$ and $\varepsilon_0>0$ such that for any $F^{\prime} \in \mathcal{F}(G)$ satisfying $|KF' \setminus F'| < \varepsilon_0 |F'|$, we have:
		\begin{equation}\label{qqqqqqqqqqq}
			\mathrm{Wdim}_{\varepsilon}\left(\rho_{X, F^{\prime}}|Y\right) \geq\left(\operatorname{mim}_{\varepsilon}\left(\rho_X |Y\right)-\beta\right)\left|F^{\prime}\right|.  
		\end{equation}
		Since $\operatorname{mdim}_{\varepsilon}\left(\rho_X |Y\right)-\beta>0$, there exists $\tau > 0$ such that $\sqrt{\tau} \operatorname{diam}\left(X, \rho_X\right)<\delta / 2$ and
		\begin{equation}\label{wwwwwww}
			(1-\tau)\left(\operatorname{mdim}_{\varepsilon}\left(\rho_X |Y\right)-\beta\right) \geq \operatorname{mdim}_{\varepsilon}\left(\rho_X |Y\right)-2 \beta.
		\end{equation}
		Applying \cite[Lemma 3.2]{LH}, there exist $\ell \in \mathbb{N}$ and finite $F_1, \ldots, F_{\ell} \subset G$ with $|K F_k \backslash$ $F_k|<\min \left(\varepsilon_0, \tau\right)| F_k |$ for every $k=1, \ldots, \ell$ such that for any good enough sofic approximation map $\sigma: G \rightarrow \operatorname{Sym}(d)$ and any $\mathcal{W} \subset[d]$ with $|\mathcal{W}| \geq(1-\tau / 2) d$, there exist $\mathcal{C}_1, \ldots, \mathcal{C}_{\ell} \subset \mathcal{W}$ satisfying the following:
		\begin{itemize}
			\item [(i)]the map $F_k \times \mathcal{C}_k \rightarrow \sigma(F_k) \mathcal{C}_k$ sending $(s, c)$ to $\sigma_s(c)$ is bijective for every $k=1, \ldots, \ell ;$
			\item[(ii)] the sets $\left\{\sigma(F_k) \mathcal{C}_k\right\}_{k=1}^{\ell}$ are pairwise disjoint with $\left|\bigsqcup_{k=1}^{\ell} \sigma(F_k) \mathcal{C}_k\right| \geq(1-\tau) d$.
		\end{itemize}
		Let $\sigma: G \rightarrow \operatorname{Sym}(d)$ be a map. Consider the subset
		$$
		\mathcal{W}=\left\{v \in[d]: \sigma_t \sigma_s(v)=\sigma_{t s}(v) \text { for every } t \in F, s \in \cup_{k=1}^{\ell} F_k\right\}.
		$$
		By the soficity of $G$, as $\sigma$ is a good enough sofic approximation, we have $|\mathcal{W}| \geq$ $(1-\tau / 2) d$ and there exist $\mathcal{C}_1, \ldots, \mathcal{C}_{\ell} \subset \mathcal{W}$ as above.
		
		Since $G$ is infinite, there exist maps $\psi_k: \mathcal{C}_k \rightarrow G$ for $k=1, \ldots, \ell$ such that the map $ \Psi: \bigsqcup_{k=1}^{\ell} F_k \times \mathcal{C}_k \rightarrow G$ sending $(s, c) \in F_k \times \mathcal{C}_k$ to $s \psi_k(c)$ is injective. Denote by $\tilde{F}$ the image of $\Psi$. Since $\left|K F_k \backslash F_k\right|<\varepsilon_0\left|F_k\right|$ for every $k=1, \ldots, \ell$, we have that
		\begin{equation}\label{eeeeeee}
			|K \tilde{F} \backslash \tilde{F}|<\varepsilon_0|\tilde{F}| \text { and }|\tilde{F}|=|\bigsqcup_{k=1}^{\ell} \sigma(F_k) \mathcal{C}_k| \geq(1-\tau) d.
		\end{equation}
		
		Fix $x_0 \in X$. For every $x \in X$ define a map $\xi_x:[d] \rightarrow X$ by $\xi_x(v)=x_0$ if $v \in[d] \backslash \bigsqcup_{k=1}^{\ell} \sigma\left(F_k\right) \mathrm{C}_k$, and
		$$
		\xi_x(v)=s \psi_k(c) x
		$$
		if $v=\sigma_s(c)$ for some $k=1, \ldots, \ell, s \in F_k$, and $c \in \mathcal{C}_k$. By \cite[Lemma 3.4]{LH}, one has $\xi_x \in \operatorname{Map}\left(\rho_X, F, \delta, \sigma\right)$.
		
		Let $\Phi: \operatorname{Map}\left(\rho_X, F, \delta, \sigma\right) \rightarrow P$ be an $\left(\varepsilon|Y, \rho_{X, \infty} \right)$-embedding. Since the map $X \rightarrow \operatorname{Map}\left(\rho_X, F, \delta, \sigma\right)$ sending $x$ to $\xi_x$ is continuous, the map $f: X \rightarrow P$ sending $x$ to $\Phi\left(\xi_x\right)$ is also continuous.
		
		\begin{claim}\label{end}
			$f$ is an $(\varepsilon, \rho_{X, \tilde{F}}|Y)$ -embedding.
		\end{claim}
		Therefore, combining (\ref{qqqqqqqqqqq}), (\ref{wwwwwww}) with (\ref{eeeeeee}), we obtain
		\begin{align*}
			\operatorname{Wdim}_{\varepsilon}(\operatorname{Map}{(\rho_X, F, \delta, \sigma), \rho_{X, \infty} |Y)}) & \geq \operatorname{Wdim}_{\varepsilon}(\rho_{X, \tilde{F}} | Y) \\
			&\geq|\tilde{F}|\left(\operatorname{mdim}_{\varepsilon}\left(\rho_X |Y\right)-\beta\right) \\
			& \geq(1-\tau) d\left(\operatorname{mdim}_{\varepsilon}\left(\rho_X |Y\right)-\beta\right) \\
			& \geq d\left(\operatorname{mdim}_{\varepsilon}\left(\rho_X |Y\right)-2\beta\right).
		\end{align*}
		This implies that
		$$\operatorname{mdim}_{\Sigma}\left(\rho_X |Y\right) \geq \operatorname{mdim}\left(\rho_X |Y\right).$$
		Combining Proposition \ref{xoxoxoxo} and Proposition \ref{popopo}, we get $\operatorname{mdim}_{\Sigma}(X | Y) \geq \operatorname{mdim}(X | Y)$.
		
		\emph{Proof of Claim \ref{end}}: Suppose that $f(x)=f\left(x^{\prime}\right)$ and $\pi(x)=\pi(x')$. We need to show $\rho_X\left(t x, t x^{\prime}\right)<\varepsilon$ for every $t \in \tilde{F}$.
		
		Since $\pi(x)=\pi(x')$, by the construction of $\xi_x$ and $\xi_{x'}$, we have 
		$$\pi \circ \xi_x=\pi \circ \xi_{x'}.$$
		
		Since $\Phi: \operatorname{Map}\left(\rho_X, F, \delta, \sigma\right) \rightarrow P$ is an $\left(\varepsilon|Y, \rho_{X, \infty}\right)$-embedding, from $f(x)=f(x^{\prime})$, we have $\Phi(\xi_x)=\Phi(\xi_{x'})$. Then
		$$
		\rho_{X, \infty}\left(\xi_x, \xi_{x^{\prime}}\right)<\varepsilon.
		$$
		
		For every $t \in \tilde{F}$ write $t=s \psi_k(c)$ for some $k=1, \ldots, \ell, s \in F_k$, and $c \in \mathcal{C}_k$. Setting $v=\sigma_s(c)$, we know that $\xi_x(v)=t x$ and $\xi_{x^{\prime}}(v)=t x^{\prime}$. Thus
		$$
		\rho_X\left(t x, t x^{\prime}\right)=\rho_X\left(\xi_x(v), \xi_{x^{\prime}}(v)\right)<\varepsilon
		$$
		as desired.
	\end{proof}
	In \cite{LQ}, Jin and Qiao proved a  product formula for mean dimension under the amenable group action. Next we show that 
	conditional mean dimension under the amenable group action also has the product formula. Before this, we need some preparations.

	Let $G$ be an infinite amenable group which acts continuously on a compact  metric space $X$, $(X^n,G)$ define the product action of $G$ on the product space $X^n$. Let
	$$\Sigma:\{\sigma:G\to \text{Sym}(d_i)\}_{n=1}^\infty$$
	be a sofic approximation sequence for $G$.  Jin and Qiao \cite[Lemma 3.1]{LQ} showed that 
	$$\Sigma^{(n)}:\{\sigma^{(n)}_i:G\to \text{Sym}(nd_i)\}_{n=1}^\infty$$
	is also a sofic approximation sequence for $G,$ where for every $i\in \N,$ the map  $$\sigma^{(n)}_i:G\to \text{Sym}(nd_i)$$
	is defined by 
	$$\sigma_{i}^{(n)}(g)((j-1)n+l)=(\sigma_i(g)(j)-1)n+l,\forall g\in G,\forall j\in [d_i],\forall l\in [n].$$
	
	\begin{prop}\label{llllllllem21}
		Let $\pi:(X,G)\to (Y,G)$ be a factor map where $G$ is a countable sofic group. Then we have 
		\begin{align*}
			\mathrm{mdim}_{\Sigma^{(n)}}(X|\pi)&=\frac1n\mathrm{mdim}_{\Sigma}(X^n|\pi^n);\\
			\mathrm{mdim}_{\Sigma^{(n)}}(Y|X)&=\frac1n\mathrm{mdim}_{\Sigma}(Y^n|X^n).
		\end{align*}
	\end{prop}
	\begin{proof}
		Let $\rho_X$ and $\rho_Y$ be compatible metrics on $X$ and $Y$ respectively. Let $\rho_X^{n}$ and $\rho_Y^{n}$ be compatible metrics on $X^n$ and $Y^n$ respectively. According to the proof of \cite[Lemma 3.2]{LQ}, we have \begin{equation}\label{eeeeeea1}\text{Map}(\rho_X^{n},F,\delta,\sigma)\subset \text{Map}(\rho_X,F,\delta,\sigma^{(n)})\subset \text{Map}(\rho_X^{n},F,\sqrt{n}\delta,\sigma)\end{equation}
		for any  $F \in \mathcal{F}(G),\delta>0$ and any map $\sigma$ from $G$ to $\text{Sym}(d).$
		
		Next we will show that
		\begin{align*}
			\mathrm{Wdim}_{\varepsilon}(\rho_X^{n},F,\delta,\sigma),\rho^{n}_{{X},\infty}|Y^n)&\leq \mathrm{Wdim}_{\varepsilon}(\rho_X,F,\delta,\sigma^{(n)}),\rho_{X,\infty}|Y)\\
			&\leq \mathrm{Wdim}_{\varepsilon}(\rho_X^{n},F,\sqrt{n}\delta,\sigma),\rho^{n}_{X,\infty}|Y^n),
		\end{align*}
		which implies that
		$$\text{mdim}_{\Sigma^{(n)}}(X|\pi)=\frac1n\text{mdim}_{\Sigma}(X^{n}|\pi^{n}).$$.
		
		There exists a continuous map $f:\text{Map}(\rho_X,F,\delta,\sigma^{(n)})\to P$ satisfying $f$ is an $(\varepsilon|Y,\rho_{X,\infty})$-embedding with $$\text{Wdim}_{\varepsilon}(\text{Map}(\rho_X,F,\delta,\sigma),\rho_{X,\infty}|Y))=\text{dim}(P).$$ 
		Therefore for every $\xi,\xi'\in \text{Map}(\rho_X,F,\delta,\sigma^{(n)})$ satisfying $f(\xi)=f(\xi')$ and $\pi\circ \varphi=\pi\circ \varphi'$, we have $\rho_{X,\infty}(\xi,\xi')<\varepsilon$. 
		
		Since $\text{Map}(\rho_X^{n},F,\delta,\sigma)\subset \text{Map}(\rho_X,F,\delta,\sigma^{(n)})$ and 
		$$\rho_{X,\infty}(\xi,\xi')=\rho_{\infty}^{n}(\xi,\xi'),\forall \xi,\xi'\in X^{nd}=(X^n)^{d},$$
		then for every $\xi,\xi'\in \text{Map}(\rho^{n}_X,F,\delta,\sigma)$ satisfying $f(\xi)=f(\xi')$ and $\pi^n\circ \xi=\pi^n\circ \xi'$, i.e. $\pi\circ \xi=\pi\circ \xi'$, we have 
		$$\rho_{X,\infty}^{n}(\xi,\xi')=\rho_{X,\infty}(\xi,\xi')<\varepsilon,$$ 
		
		which implies that $f$ is a $(\varepsilon|Y^n,\rho_{X,\infty}^{n})$-embedding.
		Therefore, \begin{align*}
			\text{Wdim}_{\varepsilon}(\text{Map}(\rho^{n}_X,F,\delta,\sigma),\rho^{n}_{X,\infty}|Y^n)&\leq \text{dim}(P)\\
			&= \text{Wdim}_{\varepsilon}(\text{Map}(\rho_X,F,\delta,\sigma),\rho_{X,\infty}|Y).
		\end{align*}
		Similarly, we can prove that $$\text{Wdim}_{\varepsilon}(\text{Map}(\rho_X,F,\delta,\sigma),\rho_{X,\infty}|Y)\leq \text{Wdim}_{\varepsilon}(\text{Map}(\rho_X^{n},F,\sqrt{n}\delta,\sigma),\rho^{n}_{X,\infty}|Y^n).$$
		This implies the desired result.
		
		On the other hand, by (\ref{eeeeeea1}), we have 
		$$(\pi^n)^d(\text{Map}(\rho_X^{n},F,\delta,\sigma))\subset \pi^{nd}(\text{Map}(\rho_X,F,\delta,\sigma^{(n)}))\subset (\pi^n)^d(\text{Map}(\rho^{n}_X,F,\sqrt{n}\delta,\sigma)).$$
		This implies that \begin{align*}
			\text{Wdim}_{\varepsilon}((\pi^n)^{d}(\text{Map}(\rho_X^{n},F,\delta,\sigma)),\rho^n_{Y,\infty})&\leq \text{Wdim}_{\varepsilon}(\pi^{nd}(\text{Map}(\rho_X,F,\delta,\sigma^{(n)})),\rho_{Y,\infty})\\
			&\leq \text{Wdim}_{\varepsilon}((\pi^n)^{d}(\text{Map}(\rho_X^{n},F,\sqrt{n}\delta,\sigma)),\rho^{n}_{Y,\infty}).
		\end{align*}
		Therefore, $$\text{mdim}_{\Sigma^{(n)}}(Y|X)=\frac1n\text{mdim}_{\Sigma}(Y^n|X^{n}).$$
		This proof is complete.
	\end{proof}
	
	\begin{cor}
		Let $\pi:(X,G)\to (Y,G)$ be a factor map where $G$ is an infinite amenable group. Then 
		\begin{align*}
			\mathrm{mdim}(X|\pi)&=\frac1n\mathrm{mdim}(X^n|\pi^n),\\
			\mathrm{mdim}(Y|X)&=\frac1n\mathrm{mdim}(Y^n|X^n).
		\end{align*}
	\end{cor}
	\begin{proof}
		Combining  Theorem \ref{ooooooq} and Proposition \ref{llllllllem21}, we get the desired result.
	\end{proof}

	\bibliographystyle{amsalpha}

\end{document}